\let\myfrac=\frac
\input eplain
\let\frac=\myfrac
\input amstex
\input epsf




\loadeufm \loadmsam \loadmsbm
\message{symbol names}\UseAMSsymbols\message{,}

\font\myfontdefault=cmr10

\font\mytdmchapfont=cmb10 at 14pt

\magnification 1200
\newif\ifinappendices
\newif\ifundefinedreferences
\newif\ifchangedreferences
\newif\ifloadreferences
\newif\ifmakebiblio
\newif\ifmaketdm
\newif\ifloadtdm
\newif\ifmakeglossary
\newif\ifloadglossary

\undefinedreferencesfalse
\changedreferencesfalse


\loadreferencestrue
\makebibliofalse
\maketdmfalse
\makeglossarytrue

\ifmaketdm%
\loadtdmfalse%
\else%
\loadtdmtrue%
\fi%


%
%

\def\turnintolatin#1{\ifcase #1 _\or i\or ii\or iii\or iv\or v\or vi\or vii\or viii\or ix\or x\or xi\or xii\or xiii\or xiv\or xv\or xvi\or xvii\or xviii\or xix\or xx\or xxi\or xxii\or xxiii\or xxiv\or xxv\or xxvi\fi}
\def\alphanum#1{\ifcase #1 _\or A\or B\or C\or D\or E\or F\or G\or H\or I\or J\or K\or L\or M\or N\or O\or P\or Q\or R\or S\or T\or U\or V\or W\or X\or Y\or Z\fi}
\def\gobbleeight#1#2#3#4#5#6#7#8{}

\newwrite\references
\newwrite\tdm
\newwrite\biblio
\newwrite\glossary

\newcount\chapno
\newcount\headno
\newcount\subheadno
\newcount\procno
\newcount\figno
\newcount\citationno
\newcount\eqnno

\def\setcatcodes{%
\catcode`\!=0 \catcode`\\=11}%

\ifloadreferences
    {\catcode`\@=11 \catcode`\_=11%
    \global\def\_@citation@Brezis{1}
\global\def\_@citation@CaffNirSpruckI{2}
\global\def\_@citation@CaffNirSpruckII{3}
\global\def\_@citation@CaffNirSpruckIII{4}
\global\def\_@citation@Calabi{5}
\global\def\_@citation@doCarmo{6}
\global\def\_@citation@doCarmoII{7}
\global\def\_@citation@CrandhallIshiiLyons{8}
\global\def\_@citation@Dold{9}
\global\def\_@citation@FriedlanderJoshi{10}
\global\def\_@citation@GilbTrud{11}
\global\def\_@citation@GuanSpruckII{12}
\global\def\_@citation@GuillemanPollack{13}
\global\def\_@citation@Guttierez{14}
\global\def\_@citation@Milnor{15}
\global\def\_@citation@Nirenberg{16}
\global\def\_@citation@Pog{17}
\global\def\_@citation@RudinOld{18}
\global\def\_@citation@Rudin{19}
\global\def\_@citation@RudinII{20}
\global\def\_@citation@ShengUrbasWang{21}
\global\def\_@citation@Simon{22}
\global\def\_@citation@Smale{23}
\global\def\_@citation@SmithCGC{24}
\global\def\_@citation@SmithPPG{25}
\global\def\_@citation@TrudWang{26}
\global\def\_@proc@ThmGeneralPlateauViscositySenseWithSingularities{1.1}
\global\def\_@proc@MostGeneralPlateau{1.2}
\global\def\_@head@TheCNSMethod{2}
\global\def\_@subhead@TheFramework{2.1}
\global\def\_@eqn@GaussCurvatureEquation{(A)}
\global\def\_@eqn@BasicEquation{(B)}
\global\def\_@proc@ZeroethAndFirstOrderBounds{2.1}
\global\def\_@subhead@BasicPropertiesOfF{2.2}
\global\def\_@proc@FirstDerivativeOfF{2.2}
\global\def\_@proc@HomogeneityOfF{2.3}
\global\def\_@proc@SecondDerivativeOfF{2.4}
\global\def\_@proc@ConcavityOfF{2.5}
\global\def\_@proc@LowerBoundsOfTraceOfHessian{2.6}
\global\def\_@subhead@Linearisation{2.3}
\global\def\_@proc@Ellipticity{2.7}
\global\def\_@proc@SuperHarmonicFunctions{2.8}
\global\def\_@proc@StrongMaximumPrinciple{2.9}
\global\def\_@subhead@TheCNSTechnique{2.4}
\global\def\_@proc@LambdaIsGreaterThanOne{2.10}
\global\def\_@proc@OffNormalDerivatives{2.11}
\global\def\_@proc@FirstBarrierFunction{2.12}
\global\def\_@proc@SecondBarrierFunction{2.13}
\global\def\_@proc@OffNormalSecondOrderBounds{2.14}
\global\def\_@proc@CorOffNormalSecondOrderBounds{2.15}
\global\def\_@proc@FinalStep{2.16}
\global\def\_@proc@OffNormalSecondOrderLowerBounds{2.17}
\global\def\_@proc@FirstOrderLowerBounds{2.18}
\global\def\_@proc@SecondOrderBounds{2.19}
\global\def\_@subhead@BoundaryToGlobal{2.6}
\global\def\_@proc@DefinitionWeakSolutions{2.20}
\global\def\_@proc@DerivativeOfLogarithm{2.21}
\global\def\_@proc@ThirdBarrierFunction{2.22}
\global\def\_@subhead@HigherOrderBounds{2.7}
\global\def\_@proc@Schauder{2.24}
\global\def\_@proc@Krylov{2.25}
\global\def\_@proc@Compactness{2.26}
\global\def\_@head@HolderSpaces{3}
\global\def\_@subhead@SmoothFunctionals{3.1}
\global\def\_@subhead@BanachSpaces{3.2}
\global\def\_@proc@InverseFunctionTheorem{3.1}
\global\def\_@proc@SubmanifoldIsManifold{3.2}
\global\def\_@proc@SubmersionTheorem{3.3}
\global\def\_@subhead@DegreeTheory{3.3}
\global\def\_@proc@DifferentialTopologicalDegree{3.4}
\global\def\_@subhead@HolderSpaces{3.4}
\global\def\_@proc@ArzelaAscoli{3.6}
\global\def\_@proc@SmoothFunctions{3.8}
\global\def\_@proc@SecondOrderFunctionals{3.9}
\global\def\_@proc@EllipticSecondOrderFunctionals{3.10}
\global\def\_@proc@Regularity{3.11}
\global\def\_@proc@SecondOrderFredholmFunctionals{3.12}
\global\def\_@proc@FunctionalIsFredholm{3.13}
\global\def\_@proc@PropernessIfPiInContextOfCurvature{3.14}
\global\def\_@proc@CorrectChoiceOfE{3.15}
\global\def\_@proc@FirstExistenceTheorem{3.16}
\global\def\_@head@TheStructureOfSingularities{4}
\global\def\_@subhead@TheHausdorffTopology{4.1}
\global\def\_@proc@NestedIntersectionIsHausdorffLimit{4.1}
\global\def\_@proc@CompactnessOfCompactSets{4.2}
\global\def\_@proc@ConvexSetsIsAClosedSet{4.3}
\global\def\_@subhead@SupportingNormals{4.2}
\global\def\_@proc@ConvergenceOfSupportingNormals{4.4}
\global\def\_@proc@CharacterisationOfClosestPoints{4.6}
\global\def\_@proc@ExistenceOfSupportingNormals{4.7}
\global\def\_@proc@ContinuityOfSupportingNormals{4.8}
\global\def\_@proc@SupportingNormalIsLocalProperty{4.9}
\global\def\_@proc@SupportingNormalForClosedConvexSets{4.10}
\global\def\_@proc@LemmaGraphFunctionIsLipschitz{4.11}
\global\def\_@proc@ExistenceOfGraphFunction{4.12}
\global\def\_@proc@ClosureOfConvexIsConvex{4.13}
\global\def\_@proc@ConvexHullOfCompactIsCompact{4.14}
\global\def\_@proc@ConvexHullOfOpenIsOpen{4.15}
\global\def\_@proc@ConvexHullOfConvexIsTheSame{4.16}
\global\def\_@subhead@LocalGeodesicProperty{4.5}
\global\def\_@proc@GBPLineIsContainedInBoundary{4.17}
\global\def\_@proc@GBPBdyPointsAreConvexHull{4.18}
\global\def\_@proc@ConvexHullHasLGP{4.19}
\global\def\_@proc@GBPMeansNotStrictlyContainedInAHemisphere{4.20}
\global\def\_@proc@ZeroIsInConvexHull{4.21}
\global\def\_@proc@NoLGPMeansPointyBit{4.22}
\global\def\_@proc@FourthBarrierFunction{4.23}
\global\def\_@proc@FifthBarrierFunction{4.24}
\global\def\_@proc@PogorelovsInteriorSecondOrderBound{4.25}
\global\def\_@proc@NonTrivialInteriorMeansSupportingNormalsContainedInHypersphere{4.26}
\global\def\_@proc@OptimalSupportingNormal{4.27}
\global\def\_@proc@StructureOfSingularities{4.28}
\global\def\_@head@WeakBarriers{5}
\global\def\_@subhead@OpenHalfSpacesAndConvexHulls{5.1}
\global\def\_@proc@InteriorOfClosureIsSame{5.1}
\global\def\_@proc@ConvexHullIsIntersectionOfOpenHalfSpaces{5.2}
\global\def\_@proc@ClosureOfInteriorIsSame{5.3}
\global\def\_@subhead@ConvexSubsetsOfTheSphere{5.2}
\global\def\_@proc@GreatCircularArcs{5.4}
\global\def\_@proc@PIsADiffeomorphism{5.5}
\global\def\_@proc@PMapsGreatCircularArcsToLineSegments{5.6}
\global\def\_@proc@PMapsConvexToConvex{5.7}
\global\def\_@proc@LemmaBasicPropertiesOfConvexSubsetsOfSpheres{5.8}
\global\def\_@proc@PMapsOpenHemispheresToOpenHalfSpaces{5.9}
\global\def\_@proc@ConvIsIntersectionOfOpenHemispheres{5.10}
\global\def\_@proc@NonEmptyToStrictlyContainedInAHemisphere{5.11}
\global\def\_@proc@OpenToClosed{5.12}
\global\def\_@proc@NonEmptyToConvex{5.13}
\global\def\_@proc@DoubleDualOfClosedIsConvexHull{5.14}
\global\def\_@proc@DualOfUnionIsIntersectionOfDual{5.15}
\global\def\_@proc@DualOfIntersectionIsConvexHullOfUnionOfDual{5.16}
\global\def\_@subhead@Links{5.4}
\global\def\_@proc@SupportingNormalsIsDualOfLink{5.18}
\global\def\_@proc@SupportingNormalsOfIntersection{5.20}
\global\def\_@subhead@DistanceFunctions{6.1}
\global\def\_@head@TheTheoryOfWeakSolutions{6}
\global\def\_@proc@PointMappingToXIsUnique{6.1}
\global\def\_@proc@ClosestPointIsUnique{6.2}
\global\def\_@proc@FirstDerivativeOfDistance{6.3}
\global\def\_@proc@TwiceDifferentiableWheneverPiIsDifferentiable{6.4}
\global\def\_@proc@ProjectionIsLipschitz{6.5}
\global\def\_@proc@PiIsDifferentiable{6.6}
\global\def\_@proc@DIsTwiceDifferentiable{6.7}
\global\def\_@proc@SecondDerivativeOfDIsSymmetric{6.8}
\global\def\_@proc@DistanceAndProjectionAreSmooth{6.10}
\global\def\_@proc@PiIsCloseToOrthogonalProjection{6.11}
\global\def\_@proc@LowerBoundsOnOrthogonalComponentOfDeterminantNearHypersurface{6.12}
\global\def\_@subhead@IntersectingConvexSets{6.3}
\global\def\_@proc@SecondDerivativeWhenProjectionIsInFirstConvex{6.13}
\global\def\_@proc@SecondDerivativeWhenProjectionIsInSecondConvex{6.14}
\global\def\_@proc@SecondDerivativeLowerBoundWhenNormalsMatch{6.15}
\global\def\_@proc@NormalLiesBetweenTheNormalsOfEachHypersurface{6.16}
\global\def\_@proc@SecondDerivativesCaseIVA{6.17}
\global\def\_@proc@SecondDerivativesCaseIVB{6.18}
\global\def\_@proc@LowerBoundOfL{6.19}
\global\def\_@proc@DerivativeOfClosestPointProjectionIsLikeNormalProjection{6.20}
\global\def\_@proc@SecondDerivativeIsInTheCorrectSet{6.21}
\global\def\_@subhead@SmoothingOperators{6.4}
\global\def\_@proc@MollifiedFunctionIsContinuous{6.22}
\global\def\_@proc@MollifiedFunctionsConvergeUniformly{6.23}
\global\def\_@proc@MollifierOfDerivative{6.24}
\global\def\_@proc@SmoothingPreservesConvexity{6.26}
\global\def\_@proc@MollifiedFunctionRemainsWithinConvexSet{6.27}
\global\def\_@proc@MollifierOfDistanceFunction{6.28}
\global\def\_@subhead@SmoothingTheIntersection{6.5}
\global\def\_@proc@SubmersionTheoremWithShapeOperator{6.29}
\global\def\_@proc@ThatSetKappaIsConvex{6.30}
\global\def\_@proc@DistanceFunctionHasHessianInCorrectSet{6.31}
\global\def\_@proc@DerivativeOfMollifiedFunctionIsInCorrectSet{6.32}
\global\def\_@proc@SmoothIntersectionHasHighCurvature{6.33}
\global\def\_@proc@WeakBarriersClosedUnderLimits{6.34}
\global\def\_@proc@WeakBarriersClosedUnderIntersection{6.35}
\global\def\_@proc@WeakBarriersClosedUnderExcision{6.37}
\global\def\_@eqn@EqnKLooksTheSameInSmallBalls{(C)}
\global\def\_@proc@EveryElementHasNonTrivialInterior{6.38}
\global\def\_@proc@SmoothPointsHaveCurvatureBoundedBelow{6.39}
\global\def\_@proc@InfimumOfVolumeIsPositive{6.40}
\global\def\_@proc@ContainedMeansLessVolume{6.41}
\global\def\_@proc@VolumeMinimiserExists{6.42}
\global\def\_@subhead@SingularitiesAndSmoothness{6.8}
\global\def\_@proc@VolumeMinimiserIsSmoothWhenNotLGP{6.43}
\global\def\_@proc@ThmDescriptionOfSingularSet{6.44}
\global\def\_@proc@InteriorBallCondition{6.45}
\global\def\_@head@Terminology{A}
    }%
\else
    \openout\references=references.tex
\fi

\ifmaketdm
    \openout\tdm=tdm.tex
\fi

%
%
\ifmakeglossary%
    \openout\glossary=glossary.tex%
\fi%
\def\gloss#1{%
\def\temp{#1}%
\edef\word{\expandafter\gobbleeight\meaning\temp}%
\edef\textoutput{\write\glossary{\word\noexpand\quad\noexpand\folio\noexpand\par}}%
\textoutput%
{\bf #1}%
}%
%
%

\newcount\newchapflag 
\newcount\showpagenumflag 

\global\eqnno = 0
\global\chapno = -1 
\global\citationno=0
\global\headno = 0
\global\subheadno = 0
\global\procno = 0
\global\figno = 0

\def\resetcounters{%
\global\eqnno = 0%
\global\headno = 0%
\global\subheadno = 0%
\global\procno = 0%
\global\figno = 0%
}

\global\newchapflag=0 
\global\showpagenumflag=0 

\def\eqninfo{(\alphanum{\the\eqnno})}
\def\chinfo{\ifinappendices\alphanum\chapno\else\the\chapno\fi}%
\def\headinfo{\ifinappendices\alphanum\headno\else\the\headno\fi}%
\def\subheadinfo{\headinfo.\the\subheadno}
\def\procinfo{\headinfo.\the\procno}
\def\figinfo{\the\figno}        
\def\citationinfo{\the\citationno}%
\def\nexteqnno{\global\advance\eqnno by 1 \eqninfo}
\def\nextheadno{\global\advance\headno by 1 \global\subheadno = 0 \global\procno = 0}
\def\nextsubheadno{\global\advance\subheadno by 1}
\def\nextprocno{\global\advance\procno by 1 \procinfo}
\def\nextfigno{\global\advance\figno by 1 \figinfo}

{\global\let\noe=\noexpand%
%
%
\catcode`\@=11%
\catcode`\_=11%
\setcatcodes%
!global!def!_@@internal@@makeref#1{%
!global!expandafter!def!csname #1ref!endcsname##1{%
!csname _@#1@##1!endcsname%
!expandafter!ifx!csname _@#1@##1!endcsname!relax%
    !write16{#1 ##1 not defined - run saving references}%
    !undefinedreferencestrue%
!fi}}%
!global!def!_@@internal@@makelabel#1{%
!global!expandafter!def!csname #1label!endcsname##1{%
!edef!temptoken{!csname #1info!endcsname}%
!ifloadreferences%
    !expandafter!ifx!csname _@#1@##1!endcsname!relax%
        !write16{#1 ##1 not hitherto defined - rerun saving references}%
        !changedreferencestrue%
    !else%
        !expandafter!ifx!csname _@#1@##1!endcsname!temptoken%
        !else%
            !write16{#1 ##1 reference has changed - rerun saving references}%
            !changedreferencestrue%
        !fi%
    !fi%
!else%
    !expandafter!edef!csname _@#1@##1!endcsname{!temptoken}%
    !edef!textoutput{!write!references{\global\def\_@#1@##1{!temptoken}}}%
    !textoutput%
!fi}}%
!global!def!makecounter#1{!_@@internal@@makelabel{#1}!_@@internal@@makeref{#1}}%
!unsetcatcodes%
}
\makecounter{eqn}%
\makecounter{ch}%
\makecounter{head}%
\makecounter{subhead}%
\makecounter{proc}%
\makecounter{fig}%
\makecounter{citation}%
\def\newref#1#2{%
\def\temptext{#2}%
\edef\bibliotextoutput{\expandafter\gobbleeight\meaning\temptext}%
\global\advance\citationno by 1\citationlabel{#1}%
\ifmakebiblio%
    \edef\fileoutput{\write\biblio{\noindent\hbox to 0pt{\hss$[\the\citationno]$}\hskip 0.2em\bibliotextoutput\medskip}}%
    \fileoutput%
\fi}%
\def\cite#1{%
$[\citationref{#1}]$%
\ifmakebiblio%
    \edef\fileoutput{\write\biblio{#1}}%
    \fileoutput%
\fi%
}%
%
%
%

\let\mypar=\par


\def\raggedleft{\leftskip=0pt plus 1fil \parfillskip=0pt}


\font\lettrinefont=cmr10 at 28pt
\def\lettrine #1[#2][#3]#4%
{\hangafter -#1 \hangindent #2
\noindent\hskip -#2 \vtop to 0pt{
\kern #3 \hbox to #2 {\lettrinefont #4\hss}\vss}}

\font\mylettrinefont=cmr10 at 28pt
\def\mylettrine #1[#2][#3][#4]#5%
{\hangafter -#1 \hangindent #2
\noindent\hskip -#2 \vtop to 0pt{
\kern #3 \hbox to #2 {\mylettrinefont #5\hss}\vss}}


\edef\Pagetitle={Blank}

\headline={\hfil\Pagetitle\hfil}

\edef\Pagefooter={Blank}
\footline={\hfil\Pagefooter\hfil}

\def\nextoddpage
{
\newpage%
\ifodd\pageno%
\else%
    \global\showpagenumflag = 0%
    \null%
    \vfil%
    \eject%
    \global\showpagenumflag = 1%
\fi%
}


\def\newchap#1#2%
{%
%
%
\global\advance\chapno by 1%
\resetcounters%
%
%
\newpage%
\ifodd\pageno%
\else%
    \global\showpagenumflag = 0%
    \null%
    \vfil%
    \eject%
    \global\showpagenumflag = 1%
\fi%
\global\newchapflag = 1%
\global\showpagenumflag = 1%
%
%
{\font\chapfontA=cmsl10 at 30pt%
\font\chapfontB=cmsl10 at 25pt%
\null\vskip 5cm%
{\chapfontA\raggedleft\hfil%
{%
\ifnum\chapno=0
    \phantom{%
    \ifinappendices%
        Annexe \alphanum\chapno%
    \else%
        \the\chapno%
    \fi}%
\else%
    \ifinappendices%
        Annexe \alphanum\chapno%
    \else%
        \the\chapno%
    \fi%
\fi%
}%
\par}%
\vskip 2cm%
{\chapfontB\raggedleft%
\lineskiplimit=0pt%
\lineskip=0.8ex%
\hfil #1\par}%
\vskip 2cm%
}%
\edef\Pagetitle{#2}%
%
%
\ifmaketdm%
    \def\temp{#2}%
    \def\tempbis{\nobreak}%
    \edef\chaptitle{\expandafter\gobbleeight\meaning\temp}%
    \edef\mynobreak{\expandafter\gobbleeight\meaning\tempbis}%
    \edef\textoutput{\write\tdm{\bigskip{\noexpand\mytdmchapfont\noindent\chinfo\ - \chaptitle\hfill\noexpand\folio}\par\mynobreak}}%
\fi%
\textoutput%
}


\def\newhead#1%
{%
%
%
\vfil
\eject%
\ifodd\pageno%
\else%
\null\vfil\eject%
\fi%
%
%
\nextheadno%
%
%
\font\headfontA=cmb10 at 48pt
\font\headfontB=cmb10 at 24pt
\null%
\vskip 15ex%
{\headfontA\raggedleft%
\headinfo\par%
\headfontB\vskip 6ex%
#1\par}%
\vskip 15ex
%
%
\ifmaketdm%
\def\temp{#1}%
\edef\sectiontitle{\expandafter\gobbleeight\meaning\temp}%
\ifnum\headno=1%
\ifinappendices%
\edef\textoutput{\write\tdm{\noindent\noexpand\medskip}}%
\textoutput%
\fi%
\else%
\edef\textoutput{\write\tdm{\noindent\noexpand\medskip}}%
\textoutput%
\fi%
\edef\textoutput{\write\tdm{\noindent{\noexpand\bf\headinfo\ -\ \sectiontitle\noexpand\hfill\noexpand\folio}\par}}%
\textoutput%
\fi%
}%


\def\newsubhead#1%
{%
\ifhmode%
    \mypar%
\fi%
\ifnum\subheadno=0%
\else%
    \goodbreak\medskip%
\fi%
\nextsubheadno%
{\goodbreak\medskip%
{\bf\noindent\subheadinfo\ #1.\ }
\nobreak\medskip}%
%
%
\ifmaketdm%
\def\temp{#1}%
\edef\subsectiontitle{\expandafter\gobbleeight\meaning\temp}%
\edef\textoutput{\write\tdm{\noindent{\noindent\quad\subheadinfo\ -\ \subsectiontitle\noexpand\hfill\noexpand\folio}\par}}%
\textoutput%
\fi%
}%

%
%


\font\mathromanten=cmr10
\font\mathromanseven=cmr7
\font\mathromanfive=cmr5
\newfam\mathromanfam
\textfont\mathromanfam=\mathromanten
\scriptfont\mathromanfam=\mathromanseven
\scriptscriptfont\mathromanfam=\mathromanfive
\def\mathroman{\fam\mathromanfam}


\font\sf=cmss12

\font\sansseriften=cmss10
\font\sansserifseven=cmss7
\font\sansseriffive=cmss5
\newfam\sansseriffam
\textfont\sansseriffam=\sansseriften
\scriptfont\sansseriffam=\sansserifseven
\scriptscriptfont\sansseriffam=\sansseriffive
\def\mathsf{\fam\sansseriffam}


\font\bftwelve=cmb12

\font\boldten=cmb10
\font\boldseven=cmb7
\font\boldfive=cmb5
\newfam\mathboldfam
\textfont\mathboldfam=\boldten
\scriptfont\mathboldfam=\boldseven
\scriptscriptfont\mathboldfam=\boldfive
\def\mathbf{\fam\mathboldfam}


\font\mycmmiten=cmmi10
\font\mycmmiseven=cmmi7
\font\mycmmifive=cmmi5
\newfam\mycmmifam
\textfont\mycmmifam=\mycmmiten
\scriptfont\mycmmifam=\mycmmiseven
\scriptscriptfont\mycmmifam=\mycmmifive

\def\hexa#1{\ifcase #1 0\or 1\or 2\or 3\or 4\or 5\or 6\or 7\or 8\or 9\or A\or B\or C\or D\or E\or F\fi}
\mathchardef\mathi="7\hexa\mycmmifam7B
\mathchardef\mathj="7\hexa\mycmmifam7C


\font\mymsbmten=msbm10 at 8pt
\font\mymsbmseven=msbm7 at 5.6pt
\font\mymsbmfive=msbm5 at 4pt
\newfam\mymsbmfam
\textfont\mymsbmfam=\mymsbmten
\scriptfont\mymsbmfam=\mymsbmseven
\scriptscriptfont\mymsbmfam=\mymsbmfive

\mathchardef\mybeth="7\hexa\mymsbmfam69
\mathchardef\mygimmel="7\hexa\mymsbmfam6A
\mathchardef\mydaleth="7\hexa\mymsbmfam6B


\def\placelabel[#1][#2]#3{{%
\setbox10=\hbox{\raise #2cm \hbox{\hskip #1cm #3}}%
\ht10=0pt%
\dp10=0pt%
\wd10=0pt%
\box10}}%


\newif\ifinproclaim%
\global\inproclaimfalse%
\def\proclaim#1{%
\medskip%
%
%
\bgroup%
\inproclaimtrue%
\setbox10=\vbox\bgroup\leftskip=0.8em\noindent{\bftwelve #1}\sf%
}

\def\endproclaim{%
\egroup%
\setbox11=\vtop{\noindent\vrule height \ht10 depth \dp10 width 0.1em}%
\wd11=0pt%
\setbox12=\hbox{\copy11\kern 0.3em\copy11\kern 0.3em}%
\wd12=0pt%
\setbox13=\hbox{\noindent\box12\box10}%
\noindent\unhbox13%
\egroup%
\medskip\ignorespaces%
}

\def\proclaim#1{%
\goodbreak\medskip
\bgroup%
\inproclaimtrue%
\noindent{\bf #1}%
\nobreak\medskip%
\sl%
}%

\def\endproclaim{%
\mypar\egroup\nobreak\medskip\ignorespaces%
}

\def\noskipproclaim#1{%
\goodbreak\medskip%
\bgroup%
\inproclaimtrue%
\noindent{\bf #1}\nobreak\sl%
}

\def\endnoskipproclaim{%
\mypar\egroup\nobreak\medskip\ignorespaces%
}


\def\ninn{{n\in\Bbb{N}}}
\def\minn{{m\in\Bbb{N}}}

\def\proof{{\noindent\bf Proof:\ }}

\def\mlim{\mathop{{\mathroman Lim}}}

\def\msup{\mathop{{\mathroman Sup}}}
\def\minf{\mathop{{\mathroman Inf}}}
\def\msf#1{{\mathsf #1}}

\def\qed{~$\square$}
\def\munion{\mathop{\cup}}
\def\minter{\mathop{\cap}}
\def\myitem#1{%
    \noindent\hbox to .5cm{\hfill#1\hss}
}

\catcode`\@=11
\def\Eqalign#1{\null\,\vcenter{\openup\jot\m@th\ialign{%
\strut\hfil$\displaystyle{##}$&$\displaystyle{{}##}$\hfil%
&&\quad\strut\hfil$\displaystyle{##}$&$\displaystyle{{}##}$%
\hfil\crcr #1\crcr}}\,}
\catcode`\@=12

\def\makeop#1{%
\global\expandafter\def\csname op#1\endcsname{{\mathroman #1}}}%

\def\makeopsmall#1{%
\global\expandafter\def\csname op#1\endcsname{{\mathroman{\lowercase{#1}}}}}%

\makeopsmall{ArcTan}%
\makeopsmall{ArcCos}%
\makeop{Arg}%
\makeop{Det}%
\makeop{Log}%
\makeop{Re}%
\makeop{Im}%
\makeop{Dim}%
\makeopsmall{Tan}%
\makeop{Mult}
\makeop{Cot}
\makeop{low}
\makeop{O}
\makeop{high}
\makeop{Ker}%
\makeopsmall{Cos}%
\makeopsmall{Sin}%
\makeop{Exp}%
\makeop{Hull}
\makeop{Conv}
\makeop{sig}
\makeopsmall{Tanh}%
\makeop{Tr}%
\makeop{End}%
\makeop{Long}%
\makeop{Ch}%
\makeop{Exp}%
\makeopsmall{Tanh}
\makeop{Eval}%
\makeop{Lift}%
\makeop{index}
\makeop{emb}%
\makeop{Int}%
\makeop{Ext}%
\makeop{Lin}
\makeop{Aire}%
\makeop{Im}%
\makeop{Conf}%
\makeop{Hom}
\makeop{Min}
\makeop{Exp}%
\makeop{Mod}%
\makeop{Log}%
\makeop{Outrad}%
\makeop{Ext}%
\makeop{Int}%
\makeop{Dist}%
\makeop{Aut}%
\makeop{Id}%
\makeop{GL}%
\makeop{Supp}
\makeop{SO}%
\makeop{Homeo}%
\makeop{Vol}%
\makeop{Ric}%
\makeop{Hess}%
\makeop{Euc}%
\makeop{Isom}%
\makeop{Max}%
\makeop{SW}%
\makeop{Coker}%
\makeop{Ind}%
\makeop{Dim}%
\makeop{even}
\makeop{Diff}%
\makeop{SL}%
\makeop{Lip}
\makeop{Emb}
\makeop{gl}
\makeop{sl}
\makeop{in}
\makeop{Imm}%
\makeop{Long}%
\makeop{Fix}%
\makeop{Smooth}%
\makeop{Wind}%
\makeop{Adj}
\makeop{Cone}
\makeop{imm}%
\makeop{Mush}%
\makeop{dVol}%
\makeop{bdd}%
\makeop{Deg}
\makeop{Struct}
\makeop{Sym}
\makeop{Quad}
\makeop{min}
\makeop{Gr}
\makeop{max}
\makeop{Ad}%
\makeop{Scal}
\makeop{equiv}
\makeop{Seq}
\makeop{Asymp}
\makeop{fibre}
\makeop{Symm}
\makeop{Coth}
\makeop{loc}%
\makeop{Diam}%
\makeop{Sinh}
\makeop{AC}
\makeopsmall{Cosh}
\makeop{Sing}
\makeop{Cyl}
\makeop{Gr}

\makeop{IDO}%
\makeop{Len}%
\makeop{Area}%
\makeop{Spec}%
\font\mycirclefont=cmsy7
\def\textcircle{{\raise 0.3ex \hbox{\mycirclefont\char'015}}}

\let\emph=\bf

\hyphenation{quasi-con-formal}

%
%

\ifmakebiblio%
    \openout\biblio=biblio.tex %
    {%
        \edef\fileoutput{\write\biblio{\bgroup\leftskip=2em}}%
        \fileoutput
    }%
\fi%

\newref{Brezis}{Brezis H., {\sl Functional analysis, Sobolev spaces and partial differential equations}, Universitext, Springer, New York, (2011)}
\newref{CaffNirSpruckI}{Caffarelli L., Nirenberg L., Spruck J., The Dirichlet problem for nonlinear second-order elliptic equations. I. Monge Amp\`ere equation, {\it Comm. Pure Appl. Math.} {\bf 37} (1984), no. 3, 369--402}
\newref{CaffNirSpruckII}{Caffarelli L., Kohn J. J., Nirenberg L., Spruck J., The Dirichlet problem for nonlinear second-order elliptic equations. II. Complex Monge Amp\`ere, and uniformly elliptic, equations, {\it Comm. Pure Appl. Math.} {\bf 38} (1985), no. 2, 209\96252}
\newref{CaffNirSpruckIII}{Caffarelli L., Nirenberg L., Spruck J., Nonlinear second-order elliptic equations. V. The Dirichlet problem for Weingarten hypersurfaces, {\sl Comm. Pure Appl. Math.} {\bf 41} (1988), no. 1, 47\9670}
\newref{Calabi}{Calabi E., Improper affine hyperspheres of convex type and a generalization of a theorem by K. J\"orgens, {\it Michigan Math. J.} {\bf 5} (1958), 105\96126}
\newref{doCarmo}{do Carmo M. P., {\it Differential Geometry of Curves and Surfaces}, Pearson, (1976)}
\newref{doCarmoII}{do Carmo M. P., {\it Riemannian Geometry}, Birkha\"user, Boston-Basel-Berlin, (1992)}
\newref{CrandhallIshiiLyons}{Crandall M. G., Ishii H., Lions P.-L., User\92s guide to viscosity solutions of second
order partial differential equations, {\it Bull. Amer. Math. Soc.} {\bf 27} (1992), no. 1, 1--67}
\newref{Dold}{Dold A., {\it Lectures on Algebraic Topology}, Classics in Mathematics, Springer-Verlag, Berlin-Heidelberg-New York, (1995)}
\newref{FriedlanderJoshi}{Friedlander F. G., Joshi M., {\it Introduction to the Theory of Distributions}, Cambridge University Press, Cambridge, (1998)}
\newref{GilbTrud}{Gilbarg D., Trudinger N. S., {\it Elliptic partical differential equations of second order}, Die Grundlehren der mathemathischen Wissenschaften, {\bf 224}, Springer-Verlag, Berlin, New York (1977)}
\newref{GuanSpruckII}{Guan B., Spruck J., The existence of hypersurfaces of constant Gauss curvature with
prescribed boundary, {\it J. Differential Geom.} {\bf 62} (2002), no. 2, 259--287}
\newref{GuillemanPollack}{Guillemin V., Pollack A., {\it Differential Topology}, Prentice-Hall, Englewood Cliffs, N.J.,
(1974)}
\newref{Guttierez}{Gutierrez C. E., {\it The Monge-Ampere Equation}, Progress in Nonlinear Differential Equations and Their Applications, Birkha\"user, (2001)}
\newref{Milnor}{Milnor J. W., {\it Topology from the differential viewpoint}, Princeton Landmarks in Mathematics, Princeton, (1997)} 
\newref{Nirenberg}{Nirenberg L., {\it Topics in non-linear functional analysis}, Courant Lecture Notes, {\bf 6}, AMS, (2001)}
\newref{Pog}{Pogorelov A. V., On the improper convex affine hyperspheres, {\it Geometriae Dedicata} {\bf 1} (1972), no. 1, 33--46}
\newref{RudinOld}{Rudin W., {\it Principles of Mathematical Analysis}, International Series in Pure \& Applied Mathematics, McGraw-Hill, (1976)}
\newref{Rudin}{Rudin W., {\it Real \& Complex Analysis}, McGraw-Hill, (1987)}
\newref{RudinII}{Rudin W., {\it Functional Analysis}, International Series in Pure \& Applied Mathematics, McGraw-Hill, (1990)}
\newref{ShengUrbasWang}{Sheng W., Urbas J., Wang X., Interior curvature bounds for a class of curvature equations. (English summary), {\it Duke Math. J.} {\bf 123} (2004), no. 2, 235--264}
\newref{Simon}{Simon L., {\it Lectures on geometric measure theory}, Centre for Mathematical Analysis, Australian National University (1984)}
\newref{Smale}{Smale S., An infinite dimensional version of Sard\92s theorem, {\it Amer. J. Math.}, {\bf 87}, (1965), 861--866}
\newref{SmithCGC}{Smith G., Compactness results for immersions of prescribed Gaussian curvature II - geometric aspects, {\it Geom. Dedicata},{\bf 172}, no. 1, (2014), 303--350}
\newref{SmithPPG}{Smith G., The Plateau problem for convex curvature functions, arXiv:1008.3545}
\newref{TrudWang}{Trudinger N. S., Wang X., On locally locally convex hypersurfaces with boundary, {\it J.
Reine Angew. Math.} {\bf 551} (2002), 11--32}

\ifmakebiblio%
    {\edef\fileoutput{\write\biblio{\egroup}}%
    \fileoutput}%
\fi%

\makeop{AC}
\makeop{pAC}

%
%
%
%
\document
\myfontdefault
\global\chapno=1
\global\showpagenumflag=1
\def\Pagetitle{\null\hfill}
\def\Pagefooter{\null\hfill}
\font\titlefontA=cmb10 at 32pt
\font\namefontA=cmb10 at 18pt%
\null
\vskip 15ex
{\raggedleft
{\titlefontA
Global\par
\vskip 2ex
Singularity Theory\par
\vskip 2ex
for the\par
\vskip 2ex
Gauss Curvature\par
\vskip 2ex
Equation\par
\vskip 4ex}
{\namefontA
Graham A. C. Smith\par}}
\nextoddpage
\global\pageno=1
\def\Pagetitle{%
\hfil\ifodd\pageno%
{\sl for the Gauss Curvature Equation}%
\else%
{\sl Global Singularity Theory}%
\fi\hfil}%

\def\Pagefooter{{\rm\turnintolatin{\the\pageno}}}
\ifloadtdm
\font\tdmfont=cmb10 at 24pt%
\null%
\vskip 15ex%
{\noindent\tdmfont Contents}
\vskip 10ex%
\noindent {\relax \tenbf 1\ -\ Introduction\hfill 1}\par 
\noindent {\noindent \hskip 1em\relax 1.1\ -\ Singularities and the Plateau problem\hfill 1}\par 
\noindent {\noindent \hskip 1em\relax 1.2\ -\ Overview and acknowledgements\hfill 4}\par 
\noindent \vskip \medskipamount 
\noindent {\relax \tenbf 2\ -\ The CNS Method\hfill 7}\par 
\noindent {\noindent \hskip 1em\relax 2.1\ -\ The framework\hfill 8}\par 
\noindent {\noindent \hskip 1em\relax 2.2\ -\ Basic properties of $F$\hfill 10}\par 
\noindent {\noindent \hskip 1em\relax 2.3\ -\ Linearisation\hfill 12}\par 
\noindent {\noindent \hskip 1em\relax 2.4\ -\ The CNS technique\hfill 14}\par 
\noindent {\noindent \hskip 1em\relax 2.5\ -\ The double normal derivative\hfill 19}\par 
\noindent {\noindent \hskip 1em\relax 2.6\ -\ Boundary to global\hfill 21}\par 
\noindent {\noindent \hskip 1em\relax 2.7\ -\ Higher order bounds\hfill 24}\par 
\noindent \vskip \medskipamount 
\noindent {\relax \tenbf 3\ -\ Degree Theory\hfill 27}\par 
\noindent {\noindent \hskip 1em\relax 3.1\ -\ Smooth mappings and differential operators\hfill 28}\par 
\noindent {\noindent \hskip 1em\relax 3.2\ -\ Banach spaces\hfill 29}\par 
\noindent {\noindent \hskip 1em\relax 3.3\ -\ Degree theory\hfill 30}\par 
\noindent {\noindent \hskip 1em\relax 3.4\ -\ H\"older spaces and H\"older norms\hfill 31}\par 
\noindent {\noindent \hskip 1em\relax 3.5\ -\ Smooth mappings of H\"older spaces\hfill 33}\par 
\noindent {\noindent \hskip 1em\relax 3.6\ -\ Existence\hfill 35}\par 
\noindent \vskip \medskipamount 
\noindent {\relax \tenbf 4\ -\ Singularities\hfill 39}\par 
\noindent {\noindent \hskip 1em\relax 4.1\ -\ The Hausdorff topology\hfill 39}\par 
\noindent {\noindent \hskip 1em\relax 4.2\ -\ Supporting normals\hfill 41}\par 
\noindent {\noindent \hskip 1em\relax 4.3\ -\ Convex sets as graphs\hfill 44}\par 
\noindent {\noindent \hskip 1em\relax 4.4\ -\ Convex hulls\hfill 47}\par 
\noindent {\noindent \hskip 1em\relax 4.5\ -\ The local geodesic property\hfill 48}\par 
\noindent {\noindent \hskip 1em\relax 4.6\ -\ Interior a-priori bounds\hfill 52}\par 
\noindent {\noindent \hskip 1em\relax 4.7\ -\ The structure of singularities\hfill 55}\par 
\noindent \vskip \medskipamount 
\noindent {\relax \tenbf 5\ -\ Duality of Convex Sets\hfill 59}\par 
\noindent {\noindent \hskip 1em\relax 5.1\ -\ Open half spaces and convex hulls\hfill 59}\par 
\noindent {\noindent \hskip 1em\relax 5.2\ -\ Convex subsets of the sphere\hfill 61}\par 
\noindent {\noindent \hskip 1em\relax 5.3\ -\ Duality\hfill 63}\par 
\noindent {\noindent \hskip 1em\relax 5.4\ -\ Links\hfill 66}\par 
\noindent \vskip \medskipamount 
\noindent {\relax \tenbf 6\ -\ Weak Barriers\hfill 69}\par 
\noindent {\noindent \hskip 1em\relax 6.1\ -\ Distance functions\hfill 70}\par 
\noindent {\noindent \hskip 1em\relax 6.2\ -\ Convex sets with smooth boundary\hfill 74}\par 
\noindent {\noindent \hskip 1em\relax 6.3\ -\ Intersecting convex sets\hfill 76}\par 
\noindent {\noindent \hskip 1em\relax 6.4\ -\ Smoothing functions and convexity\hfill 83}\par 
\noindent {\noindent \hskip 1em\relax 6.5\ -\ Smoothing the intersection\hfill 86}\par 
\noindent {\noindent \hskip 1em\relax 6.6\ -\ Weak barriers\hfill 89}\par 
\noindent {\noindent \hskip 1em\relax 6.7\ -\ The Plateau problem\hfill 93}\par 
\noindent {\noindent \hskip 1em\relax 6.8\ -\ Singularities and smoothness\hfill 95}\par 
\noindent \vskip \medskipamount 
\noindent {\relax \tenbf A\ -\ Terminology\hfill 101}\par 
\noindent \vskip \medskipamount 
\noindent {\relax \tenbf B\ -\ Index\hfill 105}\par 
\noindent \vskip \medskipamount 
\noindent {\relax \tenbf C\ -\ Bibliography\hfill 107}\par 

\nextoddpage
\global\pageno=1
\fi
\def\Pagefooter{{\rm\folio}}
\newhead{Introduction}
\newsubhead{Singularities and the Plateau problem}
The theory of singularities of solutions of totally non-linear partial differential equations presents a vast and fascinating field of mathematics about which much remains to be learnt. In this text, we study the singularities of otherwise smooth solutions of operators of Hessian type. Although little is known in the general case, when the operator is also of convex type and the ambient space is flat, a complete and satisfying description of the singularity set of any solution becomes possible. Indeed, this set decomposes as a union of convex hulls, thereby presenting a nice analogy with the linear case, where H\"ormander showed (c.f. \cite{FriedlanderJoshi}) that the wave-front set of any solution of a linear partial differential equation is a union of complete orbits of the Hamiltonian vector field of its principal symbol.
\par
The case of gaussian curvature presents a nice geometric framework within which to present this theory, though it should be borne in mind that the techniques developed in the sequel apply equally well in a far wider context. We first recall some basic definitions of riemannian geometry (c.f. \cite{doCarmoII}). Let $S$ be a smooth, oriented, embedded hypersurface in $\Bbb{R}^{n+1}$. Let $\msf{N}$ be the unit normal vector field over $S$ which is compatible with the orientation. Let $A$ be its \gloss{shape operator} (also known as the \gloss{Weingarten operator}), which is defined to be the derivative of $\msf{N}$. That is, for all $x\in S$, and for every tangent vector $V$ to $S$ at $x$,
$$
A(x)V := D\msf{N}(x)V.
$$
We recall that, for all $x$, $A(x)$ defines a linear map from the tangent space of $S$ at $x$ to itself. In particular, it always has a well-defined determinant, and we therefore define the mapping $\kappa:S\rightarrow\Bbb{R}$ by
$$
\kappa(x) := \opDet(A(x)).
$$
We call this function the \gloss{gaussian curvature} (or \gloss{extrinsic curvature}) of $S$. It is one of an immense family of possible scalar curvatures which includes the mean curvature, the so-called ``scalar curvature'', and so on. Within this family, the gaussian curvature itself is of particular interest since, after the mean curvature, it is often the most analytically tractable.
\par
We will study the Plateau problem for gaussian curvature, which asks for constant curvature hypersurfaces with prescribed boundary. Before stating the result, we believe it is worth reviewing certain geometric features of the gaussian curvature. The first concerns its relationship with convexity. Let $\opSymm$ denote the space of real, symmetric $n\times n$ matrices and consider the determinant function $\opDet:\opSymm\rightarrow\Bbb{R}$. The set $Z:=\opDet^{-1}(\left\{0\right\})$ divides $\opSymm$ into $n+1$ connected components. Indeed, for $0\leq k\leq n$, we denote by $\opSymm_k$ the subset of $\opSymm$ consisting of all invertible, symmetric matrices with exactly $k$ positive eigenvalues. The complement of $Z$ coincides with the union of all the $\opSymm_k$, each of which is connected. Significantly, $\opSymm_n$ coincides with the set of positive-definite matrices. In particular, if $\opDet(A)$ is positive and if $A$ lies in the correct connected component of the complement of $Z$, then $A$ is positive definite.
\par
Now consider the embedded hypersurface $S$. We recall that $S$ is said to be \gloss{strictly convex} whenever the matrix $A(x)$ is positive definite at every point. However, if the gaussian curvature of $S$ is everywhere strictly positive, then it follows by connectedness that $S$ is strictly convex whenever $A(x)$ is an element of $\opSymm_n$ for one single $x$. In other words, when the gaussian curvature is strictly positive, the condition of strict convexity of $S$ reduces to a single topological datum which may take one of only $n+1$ possible values.
\par
Now suppose that $(S_m)_\minn$ is a sequence of embedded hypersurfaces converging smoothly (in some reasonable sense) to $S$. If $S_m$ is strictly convex for all $m$, then $S$ is also convex, though not necessarily strictly so. However, if the gaussian curvature of $S$ is everywhere strictly positive, then $S$ is also strictly convex. That is, strict convexity, which is a-priori an open condition, becomes also a closed condition, provided again that the gaussian curvature is assumed to be strictly positive.
\par
This describes the relationship between gaussian curvature and strict convexity. It is of particular significance to us as strict convexity will play an important role in the development of the singularity theory presented in the sequel. It is a particular property of the gaussian curvature, which is not possessed, for example, by the mean curvature or the so-called ``scalar curvature''. Nonetheless, there is still an immense family of scalar curvatures which do possess this property, and which we refer to collectively as curvatures of convex type. We will not discuss this further here, but we refer the interested reader to \cite{SmithPPG} for a complete treatment.
\par
The second feature of the Plateau problem for gaussian curvature is the importance of outer barriers. This is a more subtle feature arising from the totally non-linear nature of the problem. The situation is best illustrated by the case of a circle, $C$, of unit radius in the plane. We furnish $C$ with the canonical orientation, and for $k>0$, we look for compact, oriented, embedded surfaces, $S$, in $\Bbb{R}^3$ of constant gaussian curvature $k^2$, with boundary $C$, and which lie locally to the right of this boundary curve. For $k\in]0,1[$, it is easy to construct two distinct solutions to this problem. Indeed, there are exactly two spheres of radius $1/k$ containing $C$, the first, which we denote by $S_k^+$, lying mostly above the plane, and second, which we denote by $S^-_k$, lying mostly below it. If $H$ now denotes the closed lower half-space $\left\{(x,y,z)\ |\ z\leq 0\right\}$, then the intersections $S_k^\pm\minter H$ are the desired solutions.
\par
Observe that, as $k$ tends to $1$, the pair $S_k^\pm$ degenerates to the single sphere $S_1$, which is the unique sphere with equator $C$. For $k<1$, the two solutions are then distinguished geometrically by their position with respect to $S_1\minter H$. Indeed, the ``small solution'', $S^+_k\minter H$, lies on the inside, that is to say, the concave side, of $S_1\minter H$, whilst the ``big solution'', $S^-_k\minter H$, lies on the outside, that is to say, the convex side, of this surface. In technical terms, we say that $S_1\minter H$ serves as an \gloss{outer barrier} for $S^-_k\minter H$, but not for $S^+_k\minter H$.
\par
This qualitative difference influences in a fundamental manner the behaviour of nearby solutions which follow perturbations of the boundary curve. Indeed, if $(C_t)_{\left|t\right|<\epsilon}$ is a smooth family of curves such that $C_0=C$, and if $(S^\pm_{k,t})_{\left|t\right|<\epsilon}$ are smooth families of surfaces such that $S^\pm_{k,0}=S^\pm_k\minter H$, and that, for all $t$, $S^\pm_{k,t}$ is a solution to the Plateau problem with gaussian curvature equal to $k$ and boundary curve $C_t$, then, although the existence of an outer barrier makes it relatively easy to ensure that the family $(S^+_{k,t})_{\left|t\right|<\epsilon}$ of small solutions remains within some fixed compact set, the same cannot be said for the family $(S^-_{k,t})_{\left|t\right|<\epsilon}$ of big solutions. Indeed, depending on $(C_t)_{\left|t\right|<\epsilon}$, they may diverge in arbitrarily short time. It is for this reason that the assumption of existence of an outer barrier is often indispensable in the statement of our results.
\par
We are now in a position to state the main result of this text. Here the role of the outer barrier is played by the boundary, $\partial K$, of the convex set, $K$.
\proclaim{Theorem \nextprocno, {\bf Existence and Singularities}}
\noindent Choose $k>0$. Let $K$ be a compact, convex subset of $\Bbb{R}^{n+1}$ with non-trivial interior. Let $X$ be a closed subset of $\partial K$ whose convex hull also has non-trivial interior. If $\partial K$ is smooth with gaussian curvature greater than $k$ at every point of $(\partial K)\setminus X$, then there exists a compact, convex subset $K_0\subseteq K$ with non-trivial interior such that
\medskip
\myitem{(1)} $K_0\minter\partial K=X$; and
\medskip
\myitem{(2)} $\partial K_0\minter K^o$ has constant gaussian curvature equal to $k$ in the viscosity sense.
\medskip
\noindent Furthermore, if we denote by $\opSing(K_0)$ the set of all points in $\partial K_0$ near which $\partial K_0$ is not smooth, then there exists a family $(X_\alpha)_{\alpha\in A}$ of subsets of $X$ such that
$$
\opSing(K_0) = \munion_{\alpha\in A}\opConv(X_\alpha),
$$
where $\opConv(Y)$ here denotes the convex hull of $Y$ for any set $Y$.
\endproclaim
\proclabel{ThmGeneralPlateauViscositySenseWithSingularities}
When the set $X$ has more structure, straightforward geometric arguments may often be applied to ensure that the singularity set is empty. For example, we obtain the following more classical version of the Plateau problem (c.f. \cite{GuanSpruckII} and \cite{TrudWang}).
\proclaim{Theorem \nextprocno}
\noindent Choose $k>0$. Let $K$ be a compact, convex subset of $\Bbb{R}^{n+1}$ with smooth boundary. Let $X$ be a closed subset of $\partial K$ with $C^2$ boundary $C:=\partial X$. If $\partial K$ has gaussian curvature bounded below by $k$ at every point of $(\partial K)\setminus X$, then there exists a compact, strictly convex, $C^{0,1}$ embedded hypersurface $S\subseteq\Bbb{R}^{n+1}$ such that
\medskip
\myitem{(1)} $S\subseteq K$;
\medskip
\myitem{(2)} $\partial S = C$; and
\medskip
\myitem{(3)} $S\setminus C$ is smooth and has constant gaussian curvature equal to $k$.
\endproclaim
\proclabel{MostGeneralPlateau}
{\sl\noindent{\bf Remark:\ }In fact, if $\partial X$ is smooth, then the techniques of Chapter \headref{TheCNSMethod} may readily be adapted to show that $S$ is smooth up to the boundary (and not just over its interior). We shall not study this here, although we refer the interested reader to \cite{SmithCGC} for a proof of this result following a slightly different approach.}
\newsubhead{Overview and acknowledgements}
The proof of Theorem \procref{ThmGeneralPlateauViscositySenseWithSingularities} leads us on a grand tour of various geometric and analytic aspects of the theory of non-linear partial differential equations of Hessian type. First, in Chapters \headref{TheCNSMethod} and \headref{HolderSpaces}, we prove the existence of solutions to the Plateau problem for the classical case of graphs over compact, convex subsets of $\Bbb{R}^n$ with smooth boundary. We carry this out in two stages. First, in Chapter \headref{TheCNSMethod}, we obtain compactness results for families of smooth solutions to the gaussian curvature equation, which we achieve using the technique of Caffarelli, Nirenberg and Spruck (c.f. \cite{CaffNirSpruckI}, \cite{CaffNirSpruckII} and \cite{CaffNirSpruckIII}). This general technique applies to solutions of any non-linear partial differential equation of Hessian type, and presents a fascinating object of study in its own right. With this compactness result in hand, we then use a differential-topological argument to prove existence. The formal development of this argument, which requires a considerable detour through basic functional analysis, forms the content of Chapter \headref{HolderSpaces}
\par
The local theory of singularities is introduced in Chapter \headref{TheStructureOfSingularities}. We describe the possible singularities that appear in $C^0$ limits of families of those solutions to the classical Plateau problem which we constructed in Chapter \headref{HolderSpaces}. To this end, we first present an in-depth study of the elementary geometry of convex subsets of $\Bbb{R}^{n+1}$. With a firm understanding of this theory in hand, we then apply a straightforward barrier argument dating back to the work \cite{Pog} of Pogorelov to show that singularities always lie along open straight-line segments contained within the limiting surface. This property, which we call the local geodesic property, directly implies the global structure of singularities described in Theorem \procref{ThmGeneralPlateauViscositySenseWithSingularities}.
\par
In order to apply this singularity theory, we introduce in Chapter \headref{TheTheoryOfWeakSolutions} a concept of weak supersolutions to the gaussian curvature equation, which we refer to as weak barriers. This concept, which is a more sophisticated variant of that of viscosity supersolutions, requires considerable technical work in order to establish its basic properties. This forms the content of Sections \subheadref{DistanceFunctions} to \subheadref{SmoothingTheIntersection} inclusive, which also makes use of additional properties of convex sets studied in the parenthetical Chapter \headref{WeakBarriers}. Once fully developed, however, this theory allows us to easily construct weak solutions to the Plateau problem via the Perron method, that is, by constructing a unique minimiser of a certain functional - in this case, the volume functional - amongst the set of weak barriers. Finally, in Section \subheadref{SingularitiesAndSmoothness}, we show that weak solutions are viscosity solutions, and, using the existence result of Chapter \headref{HolderSpaces}, which here serves as a local regularising operation, together with the local singularity theory developed in Chapter \headref{TheStructureOfSingularities}, we obtain the complete description of the structure of the singular parts of these solutions, thereby completing the proof of Theorem \procref{ThmGeneralPlateauViscositySenseWithSingularities}.
\par
A brief overview of the notations and terminology used throughout the text is presented in Appendix \headref{Terminology}. This text is an expanded and revised version of a mini-course presented to students in the XVII Escola de Geometria Diferencial, held in July 2012, in Manaus, Brazil. The author is grateful to L\'ucio Rodriguez and Olivier Druet for many helpful suggestions and comments.
\input GST_Chap1_150213.tex
\newhead{Degree Theory}
We develop a differential-topological degree for smooth mappings between open subsets of Banach spaces. Together with Theorem \procref{Compactness}, this yields the main result of this chapter, namely Theorem \procref{FirstExistenceTheorem}, which proves the existence of unique solutions to the classical Plateau problem for gaussian curvature in the case of graphs. This result itself constitutes an important component of the proof of Theorem \procref{MostGeneralPlateau}, and we will see in Chapter \headref{TheTheoryOfWeakSolutions}, below, how it is used construct a local regularisation operation for weak barriers.
\par
The topological degree theory we use dates back to Smale's infinite-dimensional adaptation (c.f. \cite{Smale}) of the classical finite-dimensional theory (c.f. \cite{GuillemanPollack} and \cite{Milnor}) and requires a fairly in-depth detour into functional analysis. A complete exposition of the required background material would take us too far afield, and we therefore quote a number of results without proof. We hope that this will not obscure too much the main ideas, and we refer the interested reader to the numerous excellent introductions to functional analysis (c.f. for example \cite{Brezis}, \cite{Nirenberg}, \cite{Rudin} and \cite{RudinII}) for more information.
\par
The key result is Theorem \procref{DifferentialTopologicalDegree}, which constructs a $\Bbb{Z}_2$-valued diff\-erential-topological degree for the zero set of a given smooth function between Banach spaces. The main step in our argument, encapsulated in Lemma \procref{CorrectChoiceOfE}, uses the classical Sard Theorem together with finite-dimensional reduction. In particular, even though we essentially follow Smale's reasoning (c.f. \cite{Smale}), we do not directly use the Sard-Smale Theorem. We hope that this approach will be of use to the novice reader, partly as we believe it clarifies the main ideas of Smale's result, but also because Smale's result is sometimes too specific as stated to be applied in many settings of interest in present-day mathematics.
\par
The content of this chapter is independent of the rest of the text, and the reader only interested in understanding the theory of singularities of the Gauss curvature equation may skip it if he so wishes.
\goodbreak
\newsubhead{Smooth mappings and differential operators}
Let $E$ and $F$ be normed vector spaces. Denote by $\opLin(E,F)$ the space of {\it bounded\/} linear maps from $E$ into $F$. Observe that $\opLin(E,F)$ is also a normed vector space with norm given by
\headlabel{HolderSpaces}
\subheadlabel{SmoothFunctionals}
$$
\|A\| := \msup_{x\in E\setminus\left\{0\right\}}\frac{\|Ax\|}{\|x\|}.
$$
Let $U$ be an open subset of $E$ and let $\Phi$ be a mapping from $U$ into $F$. For $x\in U$, we say that $\Phi$ is \gloss{differentiable} at $x$ whenever there exists a {\sl bounded\/} linear map $A:E\rightarrow F$ such that
$$
\mlim_{y\rightarrow 0}\frac{1}{\|y\|}\|\Phi(x+y) - \Phi(x) - A(y)\| = 0,
$$
where $y$ varies over all vectors in $E\setminus\left\{0\right\}$ having the property that $x+y\in U$. We refer to $A$ as the \gloss{derivative} of $\Phi$ at $x$. Whenever it exists, the derivative is unique, and we denote it by $D\Phi(x)$.
\medskip
{\sl\noindent{\bf Remark:\ }This definition only differs from the finite-dimensional version by the requirement that the derivative be bounded, which is unnecessary in the finite-dimensional case. Importantly, since $\opLin(E,F)$ is also a normed vector space, this allows us to iterate the concept of differentiability and thereby consider derivatives of arbitrary order.}
\medskip
We say that a function $\Phi:U\rightarrow F$ is $C^1$ whenever $D\Phi$ exists at every point of $U$ and defines a continuous function from $U$ into $\opLin(E,F)$. We define inductively the notion of higher order differentiability, and we say that $\Phi$ is $C^k$ whenever $D\Phi$ exists at every point of $U$ and defines a $C^{k-1}$ function from $U$ into $\opLin(E,F)$. We say that $\Phi$ is \gloss{smooth} whenever it is $C^k$ for all $k\in\Bbb{N}$.
\par
In order to make use of this concept, we require elementary rules for the construction of smooth functions over normed spaces. Indeed, Theorem \procref{SmoothFunctions}, below, will provide an important tool for the construction of a large family of smooth functions, which, in particular, includes almost every function that arises in geometry. At this stage, however, we recall the following three elementary rules, which are derived in exactly the same manner as in the finite-dimensional case.
\medskip
{\bf\noindent Chain Rule:\ }Let $E_1$, $E_2$ and $E_3$ be normed vector spaces. Let $U_1$ and $U_2$ be open subsets of $E_1$ and $E_2$ respectively, and let $\Phi:U_1\rightarrow U_2$ and $\Psi:U_2\rightarrow E_3$ be smooth mappings. The composition $\Psi\circ\Phi$ is smooth, and its first derivative is given by
$$
D(\Psi\circ\Phi)(x) = D\Psi(\Phi(x))D\Phi(x).
$$
{\bf\noindent Direct sums:\ }Let $E$, $F_1$,...,$F_n$ be normed vector spaces. Let $U$ be an open subset of $E$, and for $1\leq i\leq n$, let $\Phi_i:U\rightarrow F_i$ be a smooth mapping. The function $\Phi:=(\Phi_1,...,\Phi_n)$ defines a smooth mapping from $E$ into $F_1\oplus...\oplus F_n$, and its first derivative is given by
$$
D\Phi(x) = (D\Phi_1(x),...,D\Phi_n(x)).
$$
{\bf\noindent Multilinear forms:\ }Let $E_1,...,E_n$ and $F$ be normed vector spaces. Let $\Phi:E_1\oplus...\oplus E_n\rightarrow F$ be a bounded, multilinear map. $\Phi$ is smooth, and its first derivative is given by
$$
D\Phi(x_1,...,x_n)(V_1,...,V_n) = \Phi(V_1,x_2,...,x_n) + ... + \Phi(x_1,...,x_{n-1},V_n).
$$
{\sl\noindent{\bf Remark:\ }In particular, the product rule constitutes a special case of the above results.}
\newsubhead{Banach spaces}
Let $E$ be a normed vector space. We say that $E$ is a \gloss{Banach space} whenever it is complete. This extra hypothesis yields the inverse function theorem (c.f. \cite{Rudin}).
\subheadlabel{BanachSpaces}
\proclaim{Theorem \nextprocno}
\noindent Let $E$ and $F$ be Banach spaces. Let $U$ be an open subset of $E$ and let $\Phi$ be a smooth mapping from $U$ to $F$. If $D\Phi(x)$ is invertible at some point $x\in U$, then there exist neighbourhoods $V$ of $x$ in $U$, $W$ of $\Phi(x)$ in $F$ and a smooth mapping $\Psi:W\rightarrow V$ such that $W=\Phi(V)$ and
$$
\Psi\circ\Phi = \opId,\qquad \Phi\circ\Psi = \opId.
$$
\endproclaim
\proclabel{InverseFunctionTheorem}
Let $E$ be a Banach space. Let $X$ be a subset of $E$. For $n\in\Bbb{N}$, we say that $X$ is an $n$-dimensional \gloss{submanifold} of $E$ whenever there exists a Banach space $F$ with the property that for all $x\in X$, there exist neighbourhoods $U$ of $x$ in $E$ and $V$ of $(0,0)$ in $\Bbb{R}^n\times F$ and a smooth mapping $\Phi:U\rightarrow V$ with smooth inverse such that $\Phi(X\minter U) = (\Bbb{R}^n\times\left\{0\right\})\minter V$. We refer to the triplet $(\Phi,U,V)$ as a \gloss{trivialising chart} of $X$ about $x$. Recall that an abstract manifold is a separable metrisable space furnished with an atlas of charts all of whose transition maps are smooth.
\proclaim{Lemma \nextprocno}
\noindent Let $E$ be a Banach space. Let $X$ be a finite-dimensional submanifold of $E$ and let $e:X\rightarrow E$ be the canonical embedding. If $X$ is separable, then $X$ is a smooth, finite-dimensional manifold, and $e:X\rightarrow E$ is a smooth mapping.
\endproclaim
\proclabel{SubmanifoldIsManifold}
\proof By hypothesis, $X$ is separable. Furthermore, since $X$ is a subset of a normed space, it is itself a metrisable. It thus suffices to construct a smooth atlas of charts for $X$. Choose $x\in X$ and let $(\Phi,\tilde{U},\tilde{V})$ be a trivialising chart of $X$ about $x$. Denote $\phi:=\Phi|_{X\minter\tilde{U}}$, $U:=X\minter\tilde{U}$ and $V:=(\Bbb{R}^n\times\left\{0\right\})\minter\tilde{V}$. We see that $(\phi,U,V)$ defines a homeomorphism from an open subset of $X$ to an open subset of $\Bbb{R}^n$. We claim that the family of all such charts constitutes a smooth atlas for $X$. Indeed, fix $x'\in X$. Let $(\Phi',\tilde{U}',\tilde{V}')$ be another trivialising chart of $X$ about $x'$ and denote $\phi':=\Phi'|_{X\minter\tilde{U}'}$, $U':=X\minter\tilde{U}'$ and $V':=(\Bbb{R}^n\times\left\{0\right\})\minter\tilde{V}'$. Observe that $U\minter U'=\tilde{U}\minter\tilde{U}'\minter X$ and
$$
\phi'\circ(\phi^{-1})|_{\phi(U\minter U')} = \Phi'\circ(\Phi^{-1})|_{\Cal{\Phi}(\tilde{U}\minter\tilde{U}'\minter X)}.
$$
In particular, the transition map is smooth. Since $x,x'\in X$ are arbitrary, we conclude that the set of all such charts constitutes a smooth atlas, as desired. Finally, in the chart $(\phi,U,V)$, the canonical immersion coincides with $\phi^{-1}$. However, since
$$
\phi^{-1} = \Phi^{-1},
$$
it follows that this map is smooth, and this completes the proof.\qed
\medskip
Let $E$ and $F$ be two Banach spaces. Recall that a bounded linear mapping $A\in\opLin(E,F)$ is said to be \gloss{Fredholm} whenever it has closed image and both its kernel and cokernel are finite-dimensional. We define the \gloss{index} of a Fredholm mapping by
$$
\opInd(A) := \opDim(\opKer(A)) - \opDim(\opCoker(A)).
$$
Let $U$ be an open subset of $E$, and let $\Phi$ be a smooth mapping from $U$ into $F$. We say that $\Phi$ is \gloss{Fredholm} whenever $D\Phi(x)$ is a Fredholm mapping for all $x\in U$. Recall that the space of linear Fredholm mappings constitutes an open subset of $\opLin(E,F)$ and that two linear Fredholm mappings in the same connected component have the same index. It follows that if $U$ is connected, then $\opInd(D\Phi(x))$ is independent of $x\in U$, and we therefore refer to it as the index of the mapping $\Phi$. In addition, recall that the set of surjective, linear Fredholm mappings also constitutes an open subset of $\opLin(E,F)$. This is relevant to situations where we apply the following submersion theorem.
\proclaim{Theorem \nextprocno}
\noindent Let $E$ and $F$ be two Banach spaces. Let $U$ be an open subset of $E$ and let $\Phi:U\rightarrow F$ be a smooth, Fredholm map. If $D\Phi$ is surjective for all $x\in\Phi^{-1}(\left\{0\right\})$, then $\Phi^{-1}(\left\{0\right\})$ is a smooth $\opInd(\Phi)$-dimensional submanifold of $E$.
\endproclaim
\proclabel{SubmersionTheorem}
\proof Indeed, choose $x_0\in\Phi^{-1}(\left\{0\right\})$. Let $\opKer(D\Phi(x_0))$ be the kernel of $D\Phi(x_0)$. Since $D\Phi(x_0)$ is Fredholm and surjective, $\opDim(\opKer(D\Phi(x_0)))$ is equal to $\opInd(\Phi)$. By the Hahn-Banach theorem (Theorem $5.16$ of \cite{Rudin}), the identity map $\opId:\opKer(D\Phi(x_0))\rightarrow\opKer(D\Phi(x_0))$ extends to a bounded, linear projection $\pi$ from $E$ onto $\opKer(D\Phi(x_0))$. Define the mapping $\hat{\Phi}:U\rightarrow\opKer(D\Phi)(x_0)\times F$ by $\hat{\Phi}(x) := (\pi(x-x_0),\Phi(x))$. This function is smooth and, for all vectors $y\in E$, $D\hat{\Phi}(x_0)(y)=(\pi(y),D\Phi(x_0)(y))$. In particular, $D\hat{\Phi}(x_0)$ is bijective. Thus, by the closed graph theorem (c.f. Theorems $5.9$ and $5.10$ of \cite{Rudin}), $D\hat{\Phi}(x_0)$ is invertible with bounded, linear inverse. It follows from Theorem \procref{InverseFunctionTheorem} that there exist neighbourhoods $U$ of $x_0$ in $E$ and $V$ of $(0,0)$ in $\opKer(D\Phi(x_0))\times F$, and a smooth mapping $\Psi:V\rightarrow U$ such that $\hat{\Phi}(U)=V$, $\hat{\Phi}\circ\Psi=\opId$ and $\Psi\circ\hat{\Phi}=\opId$. We readily verify that $\hat{\Phi}(X\minter U)$ coincides with $(\opKer(D\Phi(x_0))\times\left\{0\right\})\minter V$, and this completes the proof.\qed
\newsubhead{Degree theory}
Let $M$ be a finite-dimensional manifold and let $E$ and $F$ be Banach spaces. Since the results of this section are local, we may suppose that $M$ is an open subset of some finite-dimensional vector space. Let $U$ be an open subset of $E$ and let $\Phi:M\times U\rightarrow F$ be a smooth Fredholm mapping of index equal to the dimension of $M$. We define the \gloss{solution space} of $\Phi$ by
\subheadlabel{DegreeTheory}
$$
\Cal{Z} := \left\{(p,x)\in M\times E\ |\ \Phi(p,x)=0\right\}.
$$
Observe that if $D\Phi(p,x)$ is surjective for all $(p,x)\in\Cal{Z}$ then, by Theorem \procref{SubmersionTheorem}, $\Cal{Z}$ is a smooth, finite-dimensional submanifold of $M\times E$ of dimension equal to $\opInd(\Phi)=\opDim(M)$. Let $\Pi:M\times E\rightarrow M$ be the projection onto the first factor. Denote by $\Pi_\Cal{Z}$ its restriction to $\Cal{Z}$. We consider the topological degree of this mapping.
\proclaim{Theorem \nextprocno}
\noindent If $D\Phi(p,x)$ is surjective for all $(p,x)\in\Cal{Z}$, and if $\Pi_\Cal{Z}$ is proper, then there exists an open, dense subset $M'\subseteq M$ with the property that, for all $p\in M'$, $\Pi_\Cal{Z}^{-1}(\left\{p\right\})$ is finite and for all $p,q\in M'$,
$$
\left|\Pi_\Cal{Z}^{-1}(\left\{p\right\})\right| = \left|\Pi_\Cal{Z}^{-1}(\left\{q\right\})\right|\text{\ Mod\ }2.
$$
\endproclaim
\proclabel{DifferentialTopologicalDegree}
\proof Denote $n=\opDim(M)=\opInd(\Phi)$. Since $\Phi$ is smooth and Fredholm, and since $D\Phi$ is surjective at every point of $\Cal{Z}=\Phi^{-1}(\left\{0\right\})$, by Theorem \procref{SubmersionTheorem}, $\Cal{Z}=\Phi^{-1}(\left\{0\right\})$ is a smooth, $n$-dimensional submanifold of $M\times E$. Since $M$ is a finite-dimensional manifold, it particular, it is separable. Thus, since $\Pi_{\Cal{Z}}:\Cal{Z}\rightarrow M$ is proper, $\Cal{Z}$ is also separable. By Lemma \procref{SubmanifoldIsManifold}, $\Cal{Z}$ is therefore a smooth, $n$-dimensional manifold and the canonical embedding $e:\Cal{Z}\rightarrow M\times E$ is a smooth mapping. In particular, $\Pi_\Cal{Z}=\Pi\circ e$ is also smooth.
\par
We now apply standard differential topological techniques to the mapping $\Pi_\Cal{Z}$. We use the terminology of \cite{GuillemanPollack}. Define $M'$ to be the set of regular values of $\Pi_{\Cal{Z}}$. If $p\in M'$, then $\Pi_{\Cal{Z}}^{-1}(\left\{p\right\})$ is discrete, and since $\Pi_{\Cal{Z}}$ is proper, this set is compact and therefore finite. Moreover, since $\Pi_{\Cal{Z}}$ is a smooth, proper mapping, by Sard's Theorem (c.f. \cite{GuillemanPollack}), $M'$ is open and dense, and the first assertion follows.
\par
Now choose $p,q\in M'$ and let $\gamma:[0,1]\rightarrow M$ be any smooth, embedded curve such that $\gamma(0)=x$ and $\gamma(1)=y$. By genericity (c.f. \cite{GuillemanPollack}), there exists another smooth, embedded curve $\tilde{\gamma}:[0,1]\rightarrow M$ which we may choose as close to $\gamma$ as we wish in the $C^\infty$ sense with the property that $\tilde{\gamma}(0)=\gamma(0)=p$, $\tilde{\gamma}(1)=\gamma(1)=q$ and $\tilde{\gamma}$ is transverse to $\Pi_{\Cal{Z}}$. If we denote by $\tilde{\Gamma}\subseteq M$ the image of $\tilde{\gamma}$, then, by transversality, $\Pi_{\Cal{Z}}^{-1}(\tilde{\Gamma})$ is a smooth, $1$-dimensional, embedded, submanifold of $\Cal{Z}$ with boundary given by
$$
\partial\Pi_{\Cal{Z}}^{-1}(\Gamma') = \Pi_{\Cal{Z}}^{-1}(\left\{p\right\})\munion\Pi_{\Cal{Z}}^{-1}(\left\{q\right\}).
$$
Since $\Pi_{\Cal{Z}}^{-1}$ is proper, $\Pi_{\Cal{Z}}^{-1}(\tilde{\Gamma})$ is compact, and therefore has an even number of boundary points. Thus
$$
\left|\Pi_{\Cal{Z}}^{-1}(\left\{p\right\})\right| + \left|\Pi_{\Cal{Z}}^{-1}(\left\{q\right\})\right| = \left|\partial\Pi_{\Cal{Z}}^{-1}(\Gamma')\right| = 0\text{\ Mod\ }2,
$$
\noindent as desired.\qed
\newsubhead{H\"older spaces and H\"older norms}
Let $E$ be a finite-dimensional normed vector space. For $\alpha\in]0,1]$ we denote by $[\cdot]_\alpha$ the \gloss{H\"older semi-norm} over $C^0(\overline{\Omega},E)$ of order $\alpha$. That is, for all $f\in C^0(\overline{\Omega},E)$,
\subheadlabel{HolderSpaces}
$$
[f]_\alpha := \msup_{x\neq y\in\overline{\Omega}}\frac{\|f(x)-f(y)\|}{\|x-y\|^\alpha}.
$$
We readily obtain
\proclaim{Lemma \nextprocno}
\noindent If $\overline{\Omega}$ is convex, then, for all continuously differentiable $f\in C^0(\overline{\Omega},E)$,
$$
[f]_1 = \|Df\|_0.
$$
\endproclaim
{\sl\noindent{\bf Remark:\ }In general, for a compact set $\Omega$ with rectifiable boundary, $[f]_1 < C(\Omega)\|Df\|_0$, where $C\geq 1$ depends on the geometry of $\Omega$. In fact, convex sets are characterised amongst all compact sets with rectifiable boundary by the property that $C(\overline{\Omega})=1$.}
\medskip
\noindent For all $\lambda = k + \alpha \in ]0,\infty[$, where $k\in\Bbb{N}$ and $\alpha\in]0,1]$, we denote by $\|\cdot\|_\lambda$ the \gloss{H\"older norm} over $C^k(\overline{\Omega})$ of order $\lambda$. That is, for all $f\in C^k(\overline{\Omega})$,
$$
\|f\|_\lambda := \sum_{i=0}^k\|D^if\|_0 + [D^kf]_\alpha.
$$
For all such $\lambda$, we denote by $C^\lambda(\overline{\Omega})$ the space of all functions $f\in C^k(\overline{\Omega})$ such that $\|f\|_\lambda<\infty$. We refer to $C^\lambda(\overline{\Omega})$ as the space of $\lambda$-times H\"older differentiable functions over $\overline{\Omega}$.
\par
We restate the classical Arzela-Ascoli theorem in the following form (c.f. \cite{Rudin}).
\proclaim{Theorem \nextprocno}
\noindent Choose $\lambda\in]0,\infty[$ and let $(f_m)_\minn$ be a sequence of functions in $C^\lambda(\overline{\Omega})$. If there exists $B>0$ such that $\|f_m\|_\lambda\leq B$ for all $m$, then there exists $f_\infty\in C^\lambda(\overline{\Omega})$ such that $\|f_\infty\|_\lambda<B$ and $(f_m)_\ninn$ subconverges to $f_\infty$ in the $C^\mu$ norm for all $\mu<\lambda$.
\endproclaim
\proclabel{ArzelaAscoli}
\noindent In particular, this yields
\proclaim{Lemma \nextprocno}
\noindent For all $\lambda\in]0,\infty[$, $(C^\lambda(\overline{\Omega}),\|\cdot\|_\lambda)$ is a Banach space.
\endproclaim
\proof Choose $\lambda>0$ and let $(f_m)_\minn$ be a Cauchy sequence of functions in $C^\lambda(\overline{\Omega})$. We need to show that $(f_m)_\ninn$ converges in $C^\lambda(\overline{\Omega})$. For all $i\in\Bbb{N}$, define the subset $F_i$ of $C^\lambda(\overline{\Omega})$ by
$$
F_i := \left\{ f_m\ |\ m\geq i\right\},
$$
and define $d_i$ to be its diameter. Since $(f_m)_\minn$ is a Cauchy sequence, the sequence $(d_i)_{i\in\Bbb{N}}$ converges to $0$.
\par
By Theorem \procref{ArzelaAscoli}, there exists a subsequence $(m_i)_{i\in\Bbb{N}}$ and a function $f_\infty$ in $C^\lambda(\overline{\Omega})$ such that $(f_{m_i})_{i\in\Bbb{N}}$ converges to $f_\infty$ in the $C^\mu$ norm for all $\mu<\lambda$. Moreover, we may assume that $m_i\geq i$ for all $i$. Consequently, for all $i$ and for all $j>i$, $f_{m_j}\in F_i$ and so
$$
\|f_{m_j} - f_i\| \leq d_i.
$$
Choose $i\in\Bbb{N}$. By Theorem \procref{ArzelaAscoli}, there exists $g\in C^\lambda(\overline{\Omega})$ such that $\|g\|_{\lambda}\leq d_i$ and $(f_{m_j}-f_i)_{j\in\Bbb{N}}$ subconverges to $g$ in the $C^\mu$ norm for all $\mu<\lambda$. However, since $(f_{m_j}-f_i)_{j\in\Bbb{N}}$ also converges to $(f_\infty - f_i)$, it follows that $g = f_\infty - f_i$. In particular, $\|f_\infty - f_i\|_\lambda = \|g\|_\lambda \leq d_i$ and we conclude that $(f_i)_{i\in\Bbb{N}}$ converges to $f_\infty$ as desired.\qed
\medskip
For all $\lambda\in]0,\infty[$, denote by $C^\lambda_0(\overline{\Omega})$ the linear subspace of $C^\lambda(\overline{\Omega})$ consisting of those functions which vanish along the boundary. Observe that, for all $\lambda$, $C^\lambda_0(\overline{\Omega})$ is a closed subspace of $C^\lambda(\overline{\Omega})$ and, in particular, is a Banach space in its own right.
\par
We leave the reader to verify that for all $\lambda\leq\mu$, $C^\mu(\overline{\Omega})$ (resp. $C^\mu_0(\overline{\Omega})$) canonically embeds as a subspace of $C^\lambda(\overline{\Omega})$ (resp. $C^\lambda_0(\overline{\Omega})$). Moreover, this embedding is continuous and if $\lambda<\mu$, it is also a compact mapping. In addition,
$$
C^\infty(\overline{\Omega}) = \minter_{\lambda>0}C^\lambda(\overline{\Omega}),\qquad C^\infty_0(\overline{\Omega}) = \minter_{\lambda>0}C^\lambda_0(\overline{\Omega}),
$$
and, moreover, a sequence $(f_m)_{\minn}$ in $C^\infty(\overline{\Omega})$ (resp. $C^\infty_0(\overline{\Omega})$) converges to a limit $f_\infty$ in $C^\infty(\overline{\Omega})$ (resp. $C^\infty_0(\overline{\Omega})$) if and only if it converges to $f_\infty$ in the $C^\lambda$-norm for all $\lambda>0$.
\newsubhead{Smooth mappings of H\"older spaces}
Let $E$ be a finite-dimensional vector space and let $U$ be an open subset of $E$. For all $\lambda>0$, we define $C^\lambda(\overline{\Omega},U)$ to be the open subset of $C^\lambda(\overline{\Omega},E)$ consisting of all functions $g$ such that $g(x)\in U$ for all $x$. Let $F$ be another finite dimensional vector space and let $\phi:\overline{\Omega}\times U\rightarrow F$ be a smooth mapping. We define the mapping $\Phi:C^0(\overline{\Omega},U)\rightarrow C^0(\overline{\Omega})$ by $\Phi(f)(x):=\phi(x,f(x))$. Together with the three construction rules already presented in Section \subheadref{SmoothFunctionals}, the following result allows us apply the techniques of the preceeding sections to almost every function that we will encounter.
\proclaim{Theorem \nextprocno}
\noindent For all $\lambda>0$ and for all $g\in C^\lambda(\overline{\Omega},U)$, $\Phi(g)\in C^\lambda(\overline{\Omega})$. Moreover, $\Phi$ defines a smooth mapping from $C^\lambda(\overline{\Omega},U)$ into $C^\lambda(\overline{\Omega})$, and for all $h\in C^\lambda(\overline{\Omega},E)$,
$$
(D\Phi(g)h)(x) = D_2\phi(g(x))h(x),
$$
where $D_2\phi$ is the partial derivative of $\phi$ with respect to the second component.
\endproclaim
\proclabel{SmoothFunctions}
Let $M$ be a finite dimensional manifold. Let $U$ be an open subset of $\oplus_{i=0}^2\opSymm(i,\Bbb{R}^n)$. Let $F:M\times\overline{\Omega}\times U\rightarrow\Bbb{R}$ be a smooth function. For all $\lambda\geq 2$, we define $\Cal{U}_0^\lambda(\overline{\Omega})$ to be the set of all functions $g$ in $C_0^\lambda(\overline{\Omega})$ such that $J^2g(x)\in U$ for all $x$. We define the mapping $\Cal{F}:M\times\Cal{U}^2_0(\overline{\Omega})\rightarrow C^0(\overline{\Omega})$ by $\Cal{F}(p,g)(x) = F(p,x,J^2g(x))$.
\proclaim{Lemma \nextprocno}
\noindent For all $\lambda\in]0,\infty[$ and for all $(p,g)\in M\times \Cal{U}_0^{\lambda+2}(\overline{\Omega})$, $\Cal{F}(p,g)\in C^\lambda(\overline{\Omega})$. Moreover, $\Cal{F}$ defines a smooth mapping from $M\times\Cal{U}_0^{\lambda+2}(\overline{\Omega})$ into $C^\lambda(\overline{\Omega})$ and its partial derivative with respect to the second component is given by
$$
(D_2\Cal{F}(p,g)h)(x) = D_3F(p,x,J^2g(x))J^2h(x),
$$
where $D_3F$ is the partial derivative of $F$ with respect to the third component.
\endproclaim
\proclabel{SecondOrderFunctionals}
\proof Choose $(p,g)\in M\times\Cal{U}^{\lambda + 2}_0(\overline{\Omega})$. Then $(p,J^2g)\in M\times C^\lambda(\overline{\Omega},\oplus_{i=0}^2\opSymm(i,\Bbb{R}^n))$ and, by Theorem \procref{SmoothFunctions}, $\Cal{F}(p,g)=F(p,x,J^2g)$ is an element of $C^\lambda(\overline{\Omega})$. Furthermore, since the mapping $g\mapsto J^2g$ defines a bounded linear map from $C^{\lambda+2}(\overline{\Omega})$ into $C^\lambda(\overline{\Omega})$, in particular it is smooth, and so, by Theorem \procref{SmoothFunctions} and the chain rule, $\Cal{F}$ defines a smooth mapping from $M\times\Cal{U}_0^{\lambda + 2}(\overline{\Omega})$ into $C^\lambda(\overline{\Omega})$ and, for all $h\in C^{\lambda+2}_0(\overline{\Omega})$,
$$
(D_2\Cal{F}(p,g)h)(x) = D_3F(p,x,J^2g(x))J^2h(x),
$$
as desired.\qed
\medskip
\noindent For all $\lambda>0$, we think of $\Cal{F}$ as a smooth family of mappings sending $\Cal{U}^{\lambda + 2}_0(\overline{\Omega})$ into $C^\lambda(\overline{\Omega})$ which is parametrised by $M$. Observe that for all $(p,g)\in M\times\Cal{U}_0^{\lambda + 2}(\overline{\Omega})$, $D_2\Cal{F}(p,g)$ is a second-order, linear, partial differential operator from $C_0^{\lambda + 2}(\overline{\Omega})$ into $C^\lambda(\overline{\Omega})$. As in Section \subheadref{HigherOrderBounds}, for all $\xi\in\oplus_{i=0}^1\opSymm(i,\Bbb{R}^n)$, we define $U_\xi\subseteq\opSymm(2,\Bbb{R}^n)$ by
$$
U_\xi := \left\{A\in\opSymm(2,\Bbb{R}^n)\ |\ (\xi,A)\in U\right\},
$$
and for all $(p,x,\xi)\in M\times\overline{\Omega}\times\oplus_{i=0}^1\opSymm(i,\Bbb{R}^n)$, we define $F_{p,x,\xi}:U_\xi\rightarrow\Bbb{R}$ by
$$
F_{p,x,\xi}(A) := F(p,x,\xi,A).
$$
We say that $F$ is an \gloss{elliptic function} whenever the derivative $DF_{p,x,\xi}(A)$ is a positive-definite matrix for all $(p,x,\xi,A)\in M\times\overline{\Omega}\times U$.
\proclaim{Lemma \nextprocno}
\noindent If $F$ is an elliptic function, then $D_2\Cal{F}(p,g)$ is an elliptic operator for all $(p,g)\in M\times\Cal{U}^2_0(\overline{\Omega})$.
\endproclaim
\proclabel{EllipticSecondOrderFunctionals}
\proof Denote by $D_3F$ the partial derivative of $F$ with respect to the third factor. By Lemma \procref{SecondOrderFunctionals}, for $(p,g)\in M\times\Cal{U}_0^{\lambda+2}(\overline{\Omega})$, for $h\in C_0^{\lambda + 2}(\overline{\Omega})$ and for $x\in\overline{\Omega}$,
$$
(D_2\Cal{F}(p,g)h)(x) = D_3F(p,x,J^2g(x))J^2h(x).
$$
If $\sigma_2(D_2\Cal{F}(p,g))(x)$ denotes the principal symbol of this operator at the point $x$ (c.f. \cite{FriedlanderJoshi}), then, for all $\xi\in\Bbb{R}^n$,
$$
\sigma_2(D_2\Cal{F}(p,g))(x)(\xi,\xi) = DF_{p,x,J^1g(x)}(D^2g(x))^{ij}\xi_i\xi_j.
$$
Since $F$ is elliptic, this is positive for all $\xi$, and since $x\in\overline{\Omega}$ is arbitrary, we conclude that $D_2\Cal{F}(p,g)$ is an elliptic operator, as desired.\qed
\medskip
Observe that the compactness result of the previous chapter only applies to smooth functions. In particular, it does not necessarily apply to functions which are only known to be H\"older differentiable of some finite order. However, the following regularity result, derived by inductively applying the classical Schauder estimates to difference quotients (c.f. Section $6.4$ of \cite{GilbTrud}) makes this distinction irrelevant.
\proclaim{Theorem \nextprocno}
\noindent Choose $\lambda>0$ and $(p,g)\in M\times\Cal{U}^{2+\lambda}(\overline{\Omega})$. If $F$ is an elliptic function and if $\Cal{F}(p,g)=0$, then $g$ is a smooth function.
\endproclaim
\proclabel{Regularity}
\noindent In order to apply the degree theory described in Sections \subheadref{SmoothFunctionals}, \subheadref{BanachSpaces} and \subheadref{DegreeTheory}, we require that $\Cal{F}$ be a Fredholm mapping. However, this readily follows from classical elliptic theory (c.f. \cite{GilbTrud}).
\proclaim{Lemma \nextprocno}
\noindent If $F$ is an elliptic function, then for all $\lambda\notin\Bbb{N}$, $\Cal{F}$ defines a Fredholm mapping from $M\times\Cal{U}_0^{\lambda+2}(\overline{\Omega})$ into $C^\lambda(\overline{\Omega})$. Moreover, $\opInd(\Cal{F})=\opDim(M)$.
\endproclaim
\proclabel{SecondOrderFredholmFunctionals}
{\sl\noindent{\bf Remark:\ }First, observe that $\Cal{F}$ is only Fredholm for non-integer values of $\lambda$. This is a significant limitation of elliptic theory in H\"older spaces. Second, we draw the reader's attention to the fact that the calculation of the index follows from general considerations. Indeed, if an elliptic operator sends sections of a bundle $E_1$ into sections of a bundle $E_2$, then the index of the operator only depends on the topology of the bundle $\opLin(E_1,E_2)$. In particular, in the case at hand, we may show that any elliptic operator sending $C_0^\infty(\overline{\Omega})$ into $C^\infty(\overline{\Omega})$ has index $0$. We refer the interested reader to \cite{Brezis}, \cite{GilbTrud} and \cite{RudinII} for more details.}
\medskip
\proof Let $D_1\Cal{F}$ and $D_2\Cal{F}$ be the partial derivatives of $\Cal{F}$ with respect to the first and second components respectively. By Lemma \procref{EllipticSecondOrderFunctionals}, for all $(p,g)\in M\times\Cal{U}_0^{\lambda + 2}(\overline{\Omega})$, $D_2\Cal{F}(p,g)$ is an elliptic operator. By classical elliptic theory, for all such $(p,g)$, $D_2\Cal{F}(p,g)$ is a Fredholm operator from $C^{\lambda + 2}_0(\overline{\Omega})$ into $C^\lambda(\overline{\Omega})$. Since $D_2\Cal{F}(p,g)$ acts on real valued functions, $\opInd(D_2\Cal{F}(p,g))=0$. Let $\pi_1$ and $\pi_2$ be the canonical projections of $M\times \Cal{U}_0^{\lambda+2}(\overline{\Omega})$ onto the first and second factors respectively. Trivially, $D\pi_2$ is Fredholm of index $\opDim(M)$. Since the composition of two Fredholm operators is Fredholm of index equal to the sum of the indices of each component, $D_2\Cal{F}(p,g)\circ D\pi_2$ is also Fredholm of index equal to $\opDim(M)$. Since $M$ has finite dimension, in particular, $D\pi_1$ has finite rank, and therefore so too does $D_1\Cal{F}(p,g)\circ D\pi_1$. Since the sum of a Fredholm operator and a finite rank operator is also a Fredholm operator of the same index, it follows that $D\Cal{F}(p,g) = D_1\Cal{F}(p,g)\circ D\pi_1 + D_2\Cal{F}(p,g)\circ D\pi_2$ is also Fredholm of index equal to $\opDim(M)$, as desired.\qed
\newsubhead{Existence}
We now recall the construction of Section \headref{TheCNSMethod}. Let $\Gamma\subseteq\opSymm(2,\Bbb{R}^n)$ be the open cone of positive-definite, symmetric matrices, and denote $U=(\oplus_{i=0}^1\opSymm(i,\Bbb{R}^n))\times\Gamma$. Let $G:\Bbb{R}^n\rightarrow\Bbb{R}$ be a smooth, convex function bounded below by $1$. Let $\phi\in C^\infty([0,1]\times\overline{\Omega},]0,\infty[)$ be a smooth family of smooth, positive functions over $\overline{\Omega}$. Let $X$ be a finite dimensional subspace of $C^\infty(\overline{\Omega})$ and let $r>0$ be a positive number, both of which we will chose presently (c.f. Lemma \procref{CorrectChoiceOfE}). We denote by $B_r$ the ball of radius $r$ about $0$ in $X$ and we define $F:[0,1]\times B_r\times\overline{\Omega}\times U\rightarrow\Bbb{R}$ by
$$
F(s,g,x,(t,\xi,A)) := \opDet(A) - (\phi_s(x) + g(x))(sG(\xi) + (1-s)).
$$\par
Observe that there exists $r>0$ such that for all $(s,g)\in[0,1]\times B_r$, $\phi_s + g>0$. Furthermore, by Lemma \procref{FirstDerivativeOfF}, $F$ is an elliptic function. For all $\lambda$, we now define the mapping $\Cal{F}:[0,1]\times B_r\times\Cal{U}^{\lambda + 2}_0(\overline{\Omega})\rightarrow C^\lambda(\overline{\Omega})$ by
$$
\Cal{F}(s,g,h)(x) := F(s,g,x,J^2h(x)).
$$
\proclaim{Lemma \nextprocno}
\noindent If $\lambda\notin\Bbb{N}$, then $\Cal{F}$ defines a smooth Fredholm mapping from $[0,1]\times B_r\times\Cal{U}_0^{\lambda+2}(\overline{\Omega})$ into $C^\lambda(\overline{\Omega})$. Moreover, $\opInd(\Cal{F})=\opDim(X)+1$.
\endproclaim
\proclabel{FunctionalIsFredholm}
\proof By Lemma \procref{SecondOrderFunctionals}, $\Cal{F}$ defines a smooth mapping from $[0,1]\times B_r\times\Cal{U}^{\lambda + 2}_0(\overline{\Omega})$ into $C^\lambda(\overline{\Omega})$, and, by Lemma \procref{SecondOrderFredholmFunctionals}, this mapping is Fredholm of Fredholm index equal to $\opDim(X)+1$. This completes the proof.\qed
\medskip
We now apply the degree theory of Section \subheadref{DegreeTheory} to $\Cal{F}$. Define the \gloss{solution space} $\Cal{Z}\subseteq[0,1]\times B_r\times\Cal{U}^{\lambda +2}_0(\overline{\Omega})$ by
$$
\Cal{Z} := \left\{ (s,g,f)\ |\ \Cal{F}(s,g,f) = 0\right\}.
$$
Let $\Pi:[0,1]\times B_r\times\Cal{U}^{\lambda+2}_0(\overline{\Omega})\rightarrow [0,1]\times B_r$ be the projection into the first two factors, and let $\Pi_\Cal{Z}$ be the restriction of $\Pi$ to $\Cal{Z}$. Let $\hat{f}\in C^\infty([0,1]\times\overline{\Omega})$ be a smooth family of strictly convex functions such that, for all $s$ and for all $x\in\partial\Omega$,
$$
\hat{f}_s(x)=0,
$$
and for all $s$ and for all $x\in\overline{\Omega}$,
$$
\Cal{F}(s,0,\hat{f}_s)(x) > 0.
$$
By compactness, upon reducing $r$ if necessary, we may suppose that for all $(s,g)\in [0,1]\times B_r$, and for all $x\in\overline{\Omega}$,
$$
\Cal{F}(s,g,\hat{f}_s)(x)>0.
$$
\proclaim{Lemma \nextprocno}
\noindent Suppose that for all $(s,g)\in[0,1]\times B_r$ and for all $x\in\overline{\Omega}$, $\phi_s(x) - g(x)>0$ and $\Cal{F}(s,g,\hat{f}_s)(x)>0$. If $\lambda\notin\Bbb{N}$, then $\Pi_\Cal{Z}$ is a proper mapping.
\endproclaim
\proclabel{PropernessIfPiInContextOfCurvature}
\proof Let $(s_m,g_m)_\minn$ be sequence in $[0,1]\times B_r$ converging to the limit, $(s_\infty,g_\infty)$, say in $[0,1]\times B_r$. Let $(f_m)_\minn$ be a sequence in $\Cal{U}^{\lambda+2}_0(\overline{\Omega})$ such that for all $m$, $(s_m,g_m,f_m)\in\Cal{Z}$. That is $\Cal{F}(s_m,g_m,f_m)=0$. Since $F$ is an elliptic function, by Theorem \procref{Regularity}, for all $m$, $f_m$ is smooth. By Lemma \procref{StrongMaximumPrinciple}, for all $m$, $f_m\geq\hat{f}_m$. Thus, by Theorem \procref{Compactness}, there exists $f_\infty\in C^\infty(\overline{\Omega})\subseteq C^{\lambda+2}_0(\overline{\Omega})$ towards which $(f_m)_\minn$ subconverges in the $C^\infty$ sense, and, in particular, $\Cal{F}(s_\infty,g_\infty,f_\infty)=0$. It remains to show that $f_\infty\in\Cal{U}_0^{\lambda+2}(\overline{\Omega})$, that is, that $f_\infty$ is strictly convex. Indeed, suppose the contrary. Since $f_\infty$ is a limit of a sequence of convex functions, it is convex. Since it is not strictly convex, there exists a point $x\in\overline{\Omega}$ at which $D^2f_\infty$ is degenerate. However, at this point, $\opDet(D^2f_\infty(x))=0$, and so $F(s_\infty,g_\infty,x,J^2f_\infty)(x)<0$. This is absurd, and we conclude that $f_\infty\in\Cal{U}_0^{\lambda+2}(\overline{\Omega})$ as asserted. In particular, $(s_\infty,g_\infty,f_\infty)\in\Cal{Z}$ and compactness follows.\qed
\medskip
\noindent Both $X$ and $r$ are now chosen to ensure surjectivity.
\proclaim{Theorem \nextprocno}
\noindent If $\lambda\notin\Bbb{N}$, then there exists a finite dimensional subspace $X\subseteq C^\infty(\overline{\Omega})$ and $r>0$ such that for all $(s,g,f)\in\Cal{Z}$, $D\Cal{F}(s,g,f)$ is surjective.
\endproclaim
\proclabel{CorrectChoiceOfE}
\proof Choose $\lambda\notin\Bbb{N}$. Define $\Cal{F}_0:[0,1]\times\Cal{U}^{\lambda+2}_0(\overline{\Omega})\rightarrow C^\lambda(\overline{\Omega})$ by
$$
\Cal{F}_0(s,f)(x) := F(s,0,x,J^2h(x)) = \opDet(D^2f(x)) - \phi_s(x)(sG(Df(x)) + (1-s)).
$$
Choose $(s_1,f_1)\in [0,1]\times\Cal{U}^{\lambda+2}_0(\overline{\Omega})$ such that $\Cal{F}_0(s_1,f_1)=0$. By Lemma \procref{FunctionalIsFredholm}, $D\Cal{F}_0(s_1,f_1)$ is a Fredholm operator. In particular, its cokernel is finite- dimensional. Let $X_1$ be the dual space to $\opIm(D\Cal{F}_0(s_1,f_1))$ in $C^\lambda(\overline{\Omega})$ with respect to the $L^2$ inner product. Observe that $E_1$ is finite-dimensional and
$$
C^\lambda(\overline{\Omega}) = \opIm(D\Cal{F}_0(s_1,f_1))\oplus X_1.
$$
Since $C^\infty(\overline{\Omega})$ is dense as a subset of $C^\lambda(\overline{\Omega})$ with respect to the $L^2$ norm, we may perturb $X_1$ to a subspace $X_1'$ of $C^\infty(\overline{\Omega})$ such that $C^\lambda(\overline{\Omega})=\opIm(D\Cal{F}_0(s_1,f_1))\oplus X_1'$. Since surjectivity of Fredholm mappings is an open property, there exists a neighbourhood $U_1$ of $(s_1,f_1)$ in $[0,1]\times\Cal{U}^{\lambda+2}_0(\overline{\Omega})$ such that, for all $(s,f)\in U_1$,
$$
C^\lambda(\overline{\Omega}) = \opIm(D\Cal{F}_0(s_1,f_1))\oplus X_1'.
$$
By Lemma \procref{PropernessIfPiInContextOfCurvature}, there exist finitely many points $(s_i,f_i)_{1\leq i\leq n}$ in $[0,1]\times\Cal{U}^{\lambda+2}_0(\overline{\Omega})$ such that $\Cal{F}_0^{-1}(\left\{0\right\})$ is contained in the union of the collection $(U_i)_{1\leq i\leq n}$. We therefore choose $X=X_1'+...+X_n'$, and for all $(s,f)\in\Cal{F}_0^{-1}(\left\{0\right\})$, we obtain,
$$
C^\lambda(\overline{\Omega}) = \opIm(D\Cal{F}_0(s_1,f_1))+ X \subseteq \opIm(D\Cal{F}(s_1,0,f_1)).
$$
Since surjectivity of Fredholm mappings is an open property, by Lemma \procref{PropernessIfPiInContextOfCurvature} again, there exists $r>0$ such that if $g\in B_r$ and if $(s,g,f)\in\Cal{Z}_r$, then $D\Cal{F}(s,g,f)$ is surjective, as desired.\qed
\proclaim{Theorem \nextprocno}
\noindent Let $\overline{\Omega}$ be a compact, convex subset of $\Bbb{R}^n$ with smooth boundary and non-trivial interior. Let $\phi\in C^\infty(\overline{\Omega})$ be a smooth, positive function. If there exists a strictly convex function $\hat{f}\in C^\infty_0(\overline{\Omega})$ such that
$$
F(D^2\hat{f}) > \phi G(D\hat{f}),\qquad \hat{f}|_{\partial\Omega} = 0,
$$
then there exists a unique strictly convex function $f\in C^\infty_0(\overline{\Omega})$ such that
$$
F(D^2f) = \phi G(Df),\qquad f|_{\partial\Omega} = 0.
$$
\endproclaim
\proclabel{FirstExistenceTheorem}
\proof First, if $f,f'\in C^\infty_0(\overline{\Omega})$ are both solutions, then, by Lemma \procref{StrongMaximumPrinciple}, both $f-f'$ and $f'-f$ attain their minimum values along the boundary. In particular, $f=f'$, and uniqueness follows.
\par
For all $t\in[0,1]$, define $G_t := tG + (1-t)$. Now fix $\alpha\in]0,1[$. Since $\hat{f}$ is strictly convex, so too is $\alpha\hat{f}$. Denote $\phi_0 := F(D^2(\alpha\hat{f}))$. Observe that $D^2\hat{f} > D^2(\alpha\hat{f})$, and so $F(D^2\hat{f})>F(D^2(\alpha\hat{f}))=\phi_0$. For $\delta>0$ and $t\in[0,1]$, define $\phi_t$ by
$$
\phi_t := \opMax((1-t/\delta)\phi_0,(1-(1-t)/\delta)\phi,\delta) > 0.
$$
For sufficiently small $\delta$, $\Cal{F}(t,0,\hat{f})>0$ for all $t\in[0,1]$. In addition, upon perturbing $\phi_t$ slightly, we may suppose that it is smooth. By Lemma \procref{CorrectChoiceOfE}, there exists a finite dimensional subspace $X\subseteq C^\infty(\overline{\Omega})$ and $r>0$ such that $D\Cal{F}$ is surjective at every point of $\Cal{Z}$. Finally, upon reducing $r$ further if necessary, we may suppose in addition that $\Cal{F}(t,g,\hat{f})>0$ for all $(t,g)\in[0,1]\times B_r(0)$.
\par
We define $f_0 = \alpha\hat{f}$. By construction, $f_0|_{\partial\Omega}=0$, $f_0$ is strictly convex, and $\Cal{F}(0,0,f_0)=0$. By uniqueness, it is the only function with these properties, so that $\Pi_\Cal{Z}^{-1}(\left\{(0,0)\right\})=\left\{(0,0,f_0)\right\}$. We now claim that $(0,0)$ is a regular value of $\Pi_\Cal{Z}$. Since $(0,0,f_0)$ is the only element of $\Pi_\Cal{Z}^{-1}(\left\{(0,0)\right\})$, it suffices to show that $D\Pi_{\Cal{Z}}$ is surjective at this point. However, let $L:=D_3\Cal{F}(0,0,f_0)$ be the partial derivative of $\Cal{F}$ with respect to the third component at $(0,0,f_0)$. We claim that $L$ is invertible. Indeed, by Lemmas \procref{FirstDerivativeOfF} and \procref{SecondOrderFunctionals},
$$
Lg = \frac{1}{n}F(D^2f_0)(D^2f_0^{-1})^{ij}g_{ij}.
$$
By classical elliptic theory, $L$ is Fredholm of index $0$. Furthermore, if $g\in\opKer(L)$ then, by the maximum principal, $g$ attains its maximum and minimum values along $\partial\Omega$. Since $g$ is also an element of $C^{\lambda+2}_0(\overline{\Omega})$, it follows that $g=0$. The kernel of $D_2\Cal{F}(0,0,f_0)$ is therefore trivial, and we conclude that $L$ is invertible, as asserted.
\par
Now let $(t,g)$ be any vector in $\Bbb{R}\times X$. By invertibility, there exists $h\in C^{\lambda+2}_0(\overline{\Omega})$ such that $D_2\Cal{F}(0,0,f_0)h=Lh=-D\Cal{F}(0,0,f_0)(t,g,0)$. In particular, $D\Cal{F}(0,0,f_0)(t,g,h) = 0$, so that $(t,g,h)$ is a tangent vector to $\Cal{Z}$ at $(0,0,f_0)$. However,
$$
D\Pi_{\Cal{Z}}(0,0,f_0)(t,g,h) = (t,g),
$$
and since $(t,g)\in\Bbb{R}\times E$ is arbitrary, we conclude that $D\Pi_{\Cal{Z}}$ is surjective at this point, as desired. It follows that $(0,0)$ is a regular value of $\Pi_{\Cal{Z}}$ and, in particular, the degree of $\Pi_\Cal{Z}$ is equal to $1$ modulo $2$.
\par
By Theorem \procref{DifferentialTopologicalDegree}, there exists a sequence $(t_m,g_m)_\minn$ of regular values of $\Pi_\Cal{Z}$ in $[0,1]\times B_r$ which converges to $(1,0)$. Since the degree of $\Pi_\Cal{Z}$ is nonzero modulo $2$, for all $m$, there exists a function $f_m\in\Cal{U}_0^{\lambda+2}(\overline{\Omega})$ such that $\Cal{F}(t_m,g_m,f_m)=0$. By Lemma \procref{PropernessIfPiInContextOfCurvature}, $(f_m)_{m\in\Bbb{N}}$ converges to a limit, $f_\infty$, say, in $\Cal{U}_0^{\lambda+2}(\overline{\Omega})$ such that $\Cal{F}(1,0,f_\infty)=0$. In other words, $f_\infty|_{\partial\Omega}=0$, $f_\infty$ is strictly convex, and
$$
F(D^2f_\infty) = \phi G(Df_\infty).
$$
Finally, by Theorem \procref{Regularity}, $f_\infty$ is smooth, and this completes the proof.\qed
\newhead{Singularities}
We study the singularities that arise in Hausdorff limits of smooth hypersurfaces of constant gaussian curvature. We show that the singular set of any such limit comprises precisely those points which posses the local geodesic property (defined in Section \subheadref{LocalGeodesicProperty}). This result, which forms the content of Theorem \procref{StructureOfSingularities}, immediately yields a global geometric characterisation of the singular set. Indeed, by Theorem \procref{GBPBdyPointsAreConvexHull} it is contained in the convex hull of some subset of the boundary.
\par
The analytic content of Theorem \procref{StructureOfSingularities} follows from interior a-priori estimates obtained using a technique which dates back to Pogorelov (c.f. \cite{Pog}, but see also \cite{Calabi} and \cite{ShengUrbasWang}), and which bears some similarities to that already used in Section \subheadref{BoundaryToGlobal} to derive global second-order bounds from second-order bounds along the boundary. However, in order to derive the full geometric consequences of these estimates, we require a thorough understanding of the elementary geometry of convex subsets of euclidean space, and this forms the content of Sections \subheadref{TheHausdorffTopology} to \subheadref{LocalGeodesicProperty} inclusive. Although some readers may find these sections elementary, we have attempted to include detailed, and hopefully clear, proofs of certain fundamental results which, to our knowledge, are not readily available elsewhere in the literature. Indeed, in Theorems \procref{GBPBdyPointsAreConvexHull} and \procref{ConvexHullHasLGP}, for example, we prove that the local geodesic property characterises convex hulls in euclidean space. Likewise, in Theorem \procref{ExistenceOfGraphFunction}, we provide a straightforward proof of the well known fact that every convex set with non-trivial interior is locally the graph of a convex, Lipschitz continuous function.
\newsubhead{The Hausdorff topology}
Let $X$ and $Y$ be two non-empty compact subsets of $\Bbb{R}^{n+1}$, we recall that $d_H(X,Y)$, the \gloss{Hausdorff distance} between $X$ and $Y$, is defined by
\headlabel{TheStructureOfSingularities}
\subheadlabel{TheHausdorffTopology}
\goodbreak
$$
d_H(X,Y) := \msup_{x\in X}\minf_{y\in Y}\|x-y\| + \msup_{y\in Y}\minf_{x\in X} \|x-y\|.
$$
We readily verify that $d_H$ defines a metric on the set of non-empty compact subsets of $\Bbb{R}^{n+1}$.
\proclaim{Lemma \nextprocno}
\noindent Let $K_1\supseteq K_2\supseteq...$ be a nested sequence of non-empty, compact subsets of $\Bbb{R}^{n+1}$ and denote
$$
K_\infty:=\minter_{m\in\Bbb{N}}K_m.
$$
Then $K_\infty$ is non-empty and $(K_m)_{\minn}$ converges to $K_\infty$ in the Hausdorff sense.
\endproclaim
\proclabel{NestedIntersectionIsHausdorffLimit}
\proof Since $K_\infty$ is the intersection of a countable, nested family of non-empty, compact sets, it is also non-empty and compact. Now suppose that $(K_m)_\minn$ does not converge to $K_\infty$ in the Hausdorff sense. Upon extracting a subsequence, we may suppose that there exists $\epsilon>0$ and a sequence $(x_m)_\minn$ such that for all $m$, $x_m\in K_m$ and $\|x_m-y\|\geq\epsilon$ for all $y\in K_\infty$. Since $x_m\in K_1$ for all $m$, by compactness we may assume that there exists $x_\infty$ towards which $(x_m)_\minn$ converges. However, for all $m$, and for all $n\geq m$, $x_n\in K_m$ and so, taking limits, $x_\infty\in K_m$. It follows that $x_\infty\in K_\infty$ and so $\|x_m-x_\infty\|\geq\epsilon$ for all $m$. This is absurd and the result follows.\qed
\proclaim{Theorem \nextprocno}
\noindent For all $R>0$, the set of non-empty, compact subsets of $\Bbb{R}^{n+1}$ contained in $\overline{B}_R(0)$ is compact in the Hausdorff topology.
\endproclaim
\proclabel{CompactnessOfCompactSets}
{\sl\noindent{\bf Remark:\ }We leave the reader to verify that the result generalises to the set of non-empty, compact subsets of any given compact metric space.}
\medskip
\proof Choose $R>0$ and let $(X_m)_\minn\subseteq\overline{B}_R(0)$ be a sequence of non-empty, compact sets. For all $k\in\Bbb{N}$, denote by $Q_k\subseteq\Bbb{R}^{n+1}$ the closed cube of side length $1/2^k$ based on the origin. That is,
$$
Q_k := [0,2^{-k}]^{n+1}.
$$
For every vector $\alpha\in\Bbb{Z}^{n+1}$, define
$$
Q_{k,\alpha} := Q_k + 2^{-k}\alpha.
$$
For all $m$ and for all $k$, define
$$
X_{m,k} := \munion_{Q_{k,\alpha}\minter X_m\neq\emptyset} Q_{k,\alpha}.
$$
For all $k$, the sequence $(X_{m,k})_\minn$ contains a subsequence converging in the Hausdorff sense to a compact limit, $X_{\infty,k}$, say, in $\Bbb{R}^{n+1}$. By a diagonal argument, we may suppose that, for all $k$, the whole sequence $(X_{m,k})_\minn$ converges to this limit. Define
$$
X_\infty := \minter_{k\in\Bbb{N}} X_{\infty,k}.
$$
For all $k$ and for all $m$, $X_{m,k+1}\subseteq X_{m,k}$. Thus, upon taking limits, $X_{\infty,k+1}\subseteq X_{\infty,k}$, and it follows by Lemma \procref{NestedIntersectionIsHausdorffLimit} that $X_\infty$ is non-empty and compact and $(X_{\infty,k})_{k\in\Bbb{N}}$ converges to $X_\infty$ in the Hausdorff sense.
\par
It remains to show that $(X_m)_\minn$ converges to $X_\infty$ in the Hausdorff sense. However, for all $k$ and for all $m$,
$$
d_H(X_m,X_{m,k})\leq 2^{-k}\sqrt{n+1}.
$$
Likewise, for all $k\leq l$ and for all $m$,
$$
d_H(X_{m,k},X_{m,l})\leq 2^{-k}\sqrt{n+1}.
$$
Taking limits yields, for all $k\leq l$,
$$
d_H(X_{\infty,k},X_{\infty,l}) \leq 2^{-k}\sqrt{n+1}.
$$
Letting $l$ tend to infinity yields, for all $k$,
$$
d_H(X_{\infty,k},X_{\infty}) \leq 2^{-k}\sqrt{n+1}.
$$
Now choose $\delta>0$ and $k>0$ such that $2^{-k}\sqrt{n+1}<\delta/3$ and let $M>0$ be such that for $m\geq M$, $d_H(X_{m,k},X_{\infty,k})<\delta/3$. Then, for $m\geq M$,
$$
d_H(X_m,X_\infty) \leq d_H(X_m,X_{m,k}) + d_H(X_{m,k},X_{\infty,k}) + d_H(X_{\infty,k},X_\infty) < \delta.
$$
Since $\delta$ may be chosen arbitrarily small, we conclude that $(X_m)_\minn$ converges to $X_\infty$ in the Hausdorff sense, as desired.\qed
\proclaim{Lemma \nextprocno}
\noindent For all $R>0$, the set of non-empty, compact, convex subsets of $\overline{B}_R(0)$ is a closed subset of the set of compact subsets of $\overline{B}_R(0)$ with respect to the Hausdorff topology. In particular, this set is also compact in the Hausdorff topology.
\endproclaim
\proclabel{ConvexSetsIsAClosedSet}
\proof Let $(K_m)_\minn$ be a sequence of non-empty, compact, convex subsets of $\overline{B}_R(0)$ converging to a compact limit $K_\infty\subseteq\overline{B}_R(0)$ in the Hausdorff sense. Choose two points $x_\infty,y_\infty\in K_\infty$. There exist sequences $(x_m)_\minn$ and $(y_m)_\minn$ converging to $x_\infty$ and $y_\infty$ respectively such that for all $m$, the points $x_m$ and $y_m$ are elements of $K_m$. Choose $t\in[0,1]$. By convexity, for all $m$, $(1-t)x_m + ty_m\in K_m$ and taking limits yields $(1-t)x_\infty + ty_\infty\in K_\infty$. Since $x_\infty, y_\infty\in K_\infty$ and $t\in [0,1]$ are arbitary, we conclude that $K_\infty$ is convex, as desired.\qed
\newsubhead{Supporting normals}
Let $\Sigma^n\subseteq\Bbb{R}^{n+1}$ be the unit sphere. We recall that if $K$ is a convex subset of $\Bbb{R}^{n+1}$, if $x$ is a point of $K$, and if $\msf{N}$ is a vector in $\Sigma^n$, then $\msf{N}$ is said to be a \gloss{supporting normal} to $K$ at $x$ whenever every other $y\in K$ satisfies
\subheadlabel{SupportingNormals}
$$
\langle y-x,\msf{N}\rangle \leq 0.
$$
Observe that supporting normals only exist at boundary points of $K$. However, when they do exist, they need not be unique. Hence, for any boundary point $x$ of $K$, we denote the set of supporting normals to $K$ at $x$ by $\Cal{N}(x;K)$, and when there is no ambiguity, we denote this set merely by $\Cal{N}(x)$. We now show that $\Cal{N}(x)$ is non-empty and compact for any boundary point $x$ of $K$, and, moreover, that this set varies semi-continuously with $x$ in a sense that will be made clear presently. We first prove a straightforward compactness result.
\proclaim{Lemma \nextprocno}
\noindent Let $(K_m)_\minn$ and $K_\infty$ be compact, convex subsets of $\Bbb{R}^{n+1}$ such that $(K_m)_\minn$ converges to $K_\infty$ in the Hausdorff sense. For all finite $m$, let $x_m$ be a boundary point of $K_m$ and let $\msf{N}_m$ be a supporting normal to $K_m$ at $x_m$. If $(x_m)_\minn$ and $(\msf{N}_m)_\minn$ converge to $x_\infty$ and $\msf{N}_\infty$ respectively, then $x_\infty$ is a boundary point of $K_\infty$ and $\msf{N}_\infty$ is a supporting normal to $K_\infty$ at $x_\infty$.
\endproclaim
\proclabel{ConvergenceOfSupportingNormals}
\proof Indeed, choose $y\in K_\infty$. There exists a sequence $(y_m)_\minn$ in $\Bbb{R}^{n+1}$ converging to $y$ such that $y_m$ is an element of $K_m$ for all $m$. For all $m$, since $\msf{N}_m$ is a supporting normal to $K_m$ at $x_m$,
$$
\langle y_m - x_m,\msf{N}_m\rangle \leq 0.
$$
Taking limits therefore yields
$$
\langle y - x_\infty,\msf{N}_\infty\rangle \leq 0.
$$
Since $y\in K_\infty$ is arbitrary, we conclude that $x_\infty\in\partial K_\infty$ and that $\msf{N}_\infty$ is a supporting normal to $K_\infty$ at $x_\infty$, as desired.\qed
\proclaim{Corollary \nextprocno}
\noindent Let $K$ be a compact, convex subset of $\Bbb{R}^{n+1}$. If $x$ is a boundary point of $K$ then $\Cal{N}(x)$ is compact.
\endproclaim
The supporting normals characterise the closest point to any exterior point of a given convex set in the following sense.
\proclaim{Lemma \nextprocno}
\noindent Let $K$ be a closed, convex subset of $\Bbb{R}^{n+1}$. Let $x$ be a point in the complement of $K$. Then $y\in K$ minimises distance to $x$ if and only if $y$ is a boundary point of $K$ and $(x-y)/\|x-y\|$ is a supporting normal to $K$ at $y$.
\endproclaim
\proclabel{CharacterisationOfClosestPoints}
\proof Suppose that $y$ is a boundary point and that $(x-y)/\|x-y\|$ is a supporting normal to $K$ at $y$. Then, for all other $z\in K$,
$$\eqalign{
\|z - x\|^2 &= \|(z-y) - (x-y)\|^2\cr
&=\|z-y\|^2 - 2\langle z - y,x-y\rangle + \|y-x\|^2\cr
&\geq \|y-x\|^2.\cr}$$
Since $z\in K$ is arbitrary, we conclude that $y$ minimises distance to $x$ in $K$ as desired.
\par
Conversely, suppose that $y\in K$ minimises distance to $x$. We claim that $\langle z - y,x - y\rangle\leq 0$ for all other $z\in K$. Indeed, suppose the contrary, so that there exists $z\in K$ such that $\langle z - y,x-y\rangle>0$. For all $t\in[0,1]$, denote $z_t=(1-t)y + tz$. By convexity, $z_t\in K$ for all $t$. Moreover,
$$
\partial_t\|z_t - x\|^2|_{t=0} = 2\langle z - y, y - x\rangle < 0.
$$
Thus, for sufficiently small $t$, $\|z_t - x\|^2<\|y-x\|^2$, which is absurd, and the assertion follows. We conclude that $y$ is a boundary point of $K$ and that $(x-y)/\|x-y\|$ is a supporting normal to $K$ at $y$, as desired.\qed
\proclaim{Theorem \nextprocno}
\noindent Let $K$ be a compact, convex subset of $\Bbb{R}^{n+1}$. For every boundary point $x$ of $K$ there exists a supporting normal $\msf{N}$ to $K$ at $x$. That is, for all $x\in\partial K$, $\Cal{N}(x)\neq\emptyset$.
\endproclaim
\proclabel{ExistenceOfSupportingNormals}
\proof Indeed, choose $x\in\partial K$. Let $(x_m)_\minn$ be a sequence of points in the complement of $K$ converging to $x$. By compactness, for all $m$, there exists a point $y_m$ in $K$ minimising distance to $x_m$. In particular, for all $m$, $d(x_m,y_m)\leq d(x_m,x)$ and so
$$
d(x,y_m) \leq d(x,x_m) + d(x_m,y_m) \leq 2d(x_m,x),
$$
so that $(y_m)_\minn$ therefore also converges to $x$. For all $m$, we define $\msf{N}_m\in\Sigma^n$ by
$$
\msf{N}_m := \frac{x_m - y_m}{\|x_m - y_m\|}.
$$
By Lemma \procref{CharacterisationOfClosestPoints}, for all $m$, $y_m$ is a boundary point of $K$ and $\msf{N}_m$ is a supporting normal to $K$ at $y_m$. After extracting a subsequence we may suppose that $(\msf{N}_m)_\minn$ converges to a unit vector $\msf{N}$, say. By Lemma \procref{ConvergenceOfSupportingNormals}, $\msf{N}$ is a supporting normal to $K$ at $x$, as desired.\qed
\medskip
We conclude by studying various properties of $\Cal{N}(x)$. We first examine the semi-continuous dependence of these sets on $x$. Thus, for two non-empty subsets $X$ and $Y$ of $\Sigma^n$, we define
$$
\delta(X,Y) := \msup_{y\in Y}\minf_{x\in X}d_\Sigma(x,y) = \msup_{y\in Y} d_\Sigma(X,y).
$$
Observe, that $\delta$ is not symmetric. However, by definition,
$$
d_{H,\Sigma}(X,Y)=\delta(X,Y)+\delta(Y,X) \geq\delta(X,Y).
$$
In particular, $\delta$ is continuous with respect to Hausdorff distance in the sphere.
\proclaim{Lemma \nextprocno}
\noindent Let $(K_m)_\minn$ and $K_\infty$ be compact, convex subsets of $\Bbb{R}^{n+1}$ such that $(K_m)_\minn$ converges to $K_\infty$ in the Hausdorff sense. For all $m$, let $x_m$ be a boundary point of $K_m$, and suppose that the sequence $(x_m)_\minn$ converges to $x_\infty$, say, in $K_\infty$. For all $\epsilon>0$ there exists $M\in\Bbb{N}$ such that for $m\geq M$, $\delta(\Cal{N}(x_\infty),\Cal{N}(x_m))<\epsilon$.
\endproclaim
\proclabel{ContinuityOfSupportingNormals}
\proof Suppose the contrary. Upon extracting a subsequence, we may suppose that there exists $\epsilon>0$ such that $\delta(\Cal{N}(x_\infty),\Cal{N}(x_m))\geq\epsilon$ for all $m$. For all $m$, there therefore exists $\msf{N}_m\in\Cal{N}(x_m)$ such that $d(\Cal{N}(x_\infty),\msf{N}_m)\geq\epsilon$. By compactness of the sphere, we may suppose that there exists $\msf{N}_\infty$ towards which $(\msf{N}_m)_\minn$ converges, and taking limits yields $d(\Cal{N}(x_\infty),\msf{N}_\infty)\geq\epsilon$. However, by Lemma \procref{ConvergenceOfSupportingNormals}, $\msf{N}_\infty\in\Cal{N}(x_\infty)$.  This is absurd, and the result follows.\qed
\medskip
We now show that $\Cal{N}(x)$ is in fact defined locally. Indeed, recall that if $K$ and $L$ are compact, convex sets, then so too is their intersection.
\proclaim{Lemma \nextprocno}
\noindent Let $K$ and $L$ be compact, convex subsets of $\Bbb{R}^{n+1}$. Let $x$ be a boundary point of $K\minter L$ and let $\msf{N}\in\Sigma^n$ be a supporting normal to $K\minter L$ at $x$. If $x\in L^o$, then $x$ is also a boundary point of $K$ and $\msf{N}$ is also a supporting normal to $K$ at $x$. That is, for all such $x$, $\Cal{N}(x;K)=\Cal{N}(x;K\minter L)$.
\endproclaim
\proclabel{SupportingNormalIsLocalProperty}
\proof It suffices to show that for all $y\in K$,
$$
\langle y - x,\msf{N}\rangle \leq 0.
$$
However, suppose the contrary. There exists $y\in K$ such that $\langle y - x,\msf{N}\rangle >0$. For all $t\in[0,1]$, denote $y_t := (1-t)x + ty$. Then, for all $t$,
$$
\langle y_t - x,\msf{N}\rangle = t\langle y - x,\msf{N}\rangle > 0.
$$
However, by convexity, $y_t$ is an element of $K$ for all $t$. Furthermore, since $x$ is an interior point of $L$, for sufficiently small $t$, $y_t$ is also an element of $L$. That is, for sufficiently small $t$, $y_t$ is an element of $K\minter L$, so that $\langle y_t - x,\msf{N}\rangle\leq 0$. This is absurd, and the result follows.\qed
\medskip
Finally, we include non-compact convex subsets of $\Bbb{R}^{n+1}$ into this framework.
\proclaim{Lemma \nextprocno}
\noindent Let $K$ be a closed, convex subset of $\Bbb{R}^{n+1}$. For every boundary point $x$ of $K$, there exists a supporting normal $\msf{N}$ to $K$ at $x$. That is, for all $x\in\partial K$, $\Cal{N}(x)\neq\emptyset$.
\endproclaim
\proclabel{SupportingNormalForClosedConvexSets}
\proof Indeed, let $x$ be a boundary point of $K$. For all $r>0$, denote $K_r:=K\minter\overline{B}_r(x)$, and observe that $K_r$ is compact and convex. Now fix $r>0$. Trivially, $\Cal{N}(x;K)\subseteq\Cal{N}(x;K_r)$. Conversely, by Lemma \procref{SupportingNormalIsLocalProperty}, for all $s>r$, $\Cal{N}(x;K_s)=\Cal{N}(x;K_r)$. In particular, if $\msf{N}\in\Cal{N}(x,K_r)$ and if $y\in K$, then since $y\in K_s$ for some $s>r$, $\langle y - x,\msf{N}\rangle \leq 0$. We conclude that $\Cal{N}(x;K_r)\subseteq\Cal{N}(x;K)$, and the two sets therefore coincide. In particular, by Theorem \procref{ExistenceOfSupportingNormals}, $\Cal{N}(x;K)=\Cal{N}(x;K_r)$ is non-empty, and the result follows.\qed
\newsubhead{Convex sets as graphs}
Let $K$ be a compact, convex subset of $\Bbb{R}^{n+1}$. Let $x$ be a boundary point of $K$ and let $\msf{N}$ be a supporting normal to $K$ at $x$. Upon applying an affine isometry, we may suppose that $x=0$ and that $\msf{N}$ is any given unit vector in the sphere so that the results which follow are completely general. We decompose $\Bbb{R}^{n+1}$ as $\Bbb{R}^n\times\Bbb{R}$ and we use the notation outlined in Appendix \headref{Terminology}. For $C,r>0$, we say that $\partial K$ is a $C$-Lipschitz \gloss{graph} over a radius $r$ near $0$ whenever there exists a $C$-Lipschitz function $f:B_r'(0)\rightarrow]-2Cr,2Cr[$ such that the intersection of $\partial K$ with $B_r'(0)\times]-2Cr,2Cr[$ coincides with the graph of $f$ over $B_r'(0)$.
\proclaim{Lemma \nextprocno}
\noindent Let $K$ be a compact, convex subset of $\Bbb{R}^{n+1}$ and suppose that $0$ is a boundary point of $K$.
Choose $\theta\in[0,\pi/2[$ and $r>0$ and suppose that for all $x\in\partial K\minter B_r(0)$, and for every supporting normal $\msf{N}$ to $K$ at $x$,
$$
\langle\msf{N},-e_{n+1}\rangle \geq \opCos(\theta).
$$
Then, for all $(x',s),(y',t)\in \partial K\minter B_r(0)$,
$$
\left|s-t\right| \leq \opTan(\theta)\|x'-y'\|.
$$
\endproclaim
\proclabel{LemmaGraphFunctionIsLipschitz}
\proof Indeed, let $\msf{N}$ be a supporting normal to $K$ at $x$ and let $\msf{N}'$ be its orthogonal projection onto $\Bbb{R}^n$. In particular,
$$
\|\msf{N}'\|^2= 1 - \langle\msf{N},e_{n+1}\rangle^2\leq 1  -\opCos^2(\theta) = \opSin^2(\theta).
$$
Using the Cauchy-Schwarz inequality, we obtain
$$\eqalign{
\langle y - x,\msf{N}\rangle &= \langle y - x, \langle\msf{N},e_{n+1}\rangle e_{n+1}\rangle + \langle y - x,\msf{N}'\rangle\cr
&=\langle\msf{N},-e_{n+1}\rangle(s-t) + \langle y' - x',\msf{N}'\rangle\cr
&\geq\langle\msf{N},-e_{n+1}\rangle(s-t) - \|x'-y'\|\opSin(\theta).\cr}$$
However, by definition of the supporting normal, $\langle y - x,\msf{N}\rangle\leq 0$, and so
$$
(s-t) \leq \|x'-y'\|\frac{\opSin(\theta)}{\langle\msf{N},-e_{n+1}\rangle} \leq \opTan(\theta)\|x'-y'\|.
$$
By symmetry, we conclude that
$$
\left|s-t\right| \leq \opTan(\theta)\|x' - y'\|,
$$
as desired.\qed
\proclaim{Theorem \nextprocno}
\noindent Let $K$ be a compact, convex subset of $\Bbb{R}^{n+1}$ and suppose that $0$ is a boundary point of $K$. Choose $\theta\in[0,\pi/2[$ and $r>0$ and suppose that for all $x\in\partial K\minter B_r(0)$, and for every supporting normal $\msf{N}$ to $K$ at $x$,
$$
\langle\msf{N},-e_{n+1}\rangle \geq \opCos(\theta).
$$
\noindent Then, denoting $C:=\opTan(\theta)$ and $\rho:=\frac{r}{\sqrt{1+4C^2}}$, there exists a unique function $f:B'_\rho(0)\rightarrow]-C\rho,C\rho[$ such that
\medskip
\myitem{(1)} $f(0)=0$;
\medskip
\myitem{(2)} $f$ is convex and $C$-Lipschitz; and
\medskip
\myitem{(3)} $(\partial K)\minter (B'_\rho(0)\times]-2C\rho,2C\rho[)$ coincides with the graph of $f$.
\endproclaim
\proclabel{ExistenceOfGraphFunction}
{\sl\noindent{\bf Remark:\ }In other words, $\partial K$ is a $C$-Lipschitz graph over a radius $\rho$ near $0$.}
\medskip
\proof Observe that $B'_\rho(0)\times]-2C\rho,2C\rho[\subseteq B_r(0)$. For all $x'\in B_\rho'(0)$, denote $L_{x'}:=\left\{x'\right\}\times]-2C\rho,2C\rho[$ and consider the set $L_{x'}\minter\partial K$. First, if $s,t\in]-2C\rho,2C\rho[$ are such that $(x',s),(x',t)\in\partial K$, then, by Lemma \procref{LemmaGraphFunctionIsLipschitz}, $\left|s-t\right|=0$, and so $s=t$. It follows that $L_{x'}\minter\partial K$ contains at most one point.
\par
We now prove existence. For all $t\in]-2C\rho,2C\rho[$, denote $B'_t := B_\rho'(0)\times\left\{t\right\}$. We first claim that $B'_{\pm C\rho}$ does not intersect $\partial K$. Indeed, otherwise, there exists $x'\in B'_\rho(0)$ such that $(x',\pm C\rho)\in\partial K$. However, since $(0,0)\in\partial K$, by Lemma \procref{LemmaGraphFunctionIsLipschitz},
$$
C\rho \leq C\|x'\| < C\rho.
$$
This is absurd, and the assertion follows. In particular, by connectedness, $B'_{\pm C\rho}$ is entirely contained either in the interior of $K$, or in $\Bbb{R}^{n+1}\setminus K$. However, for any supporting normal $\msf{N}$ to $K$ at $(0,0)$,
$$
\langle (0,-C\rho) - (0,0),\msf{N}\rangle = 2C\rho\langle -e_{n+1},\msf{N}\rangle \geq 2C\rho\opCos(\theta) > 0,
$$
so that the point $(0,-C\rho)$ does not lie in $K$. In particular, $B'_{-C\rho}$ intersects $\Bbb{R}^{n+1}\setminus K$ non-trivially, and is therefore entirely contained in $\Bbb{R}^{n+1}\setminus K$.
\par
We now show that $B'_{C\rho}$ is entirely contained in the interior of $K$. By hypothesis, $(0,0)$ lies in $L_0\minter\partial K$, and, by uniqueness, $L_0\minter\partial K$ contains no other point. However, by compactness and convexity, $L_0\minter K$ is a relatively closed, connected subset of $L_0$. Since, furthermore, the relative boundary of $L_0\minter K$ in $L_0$ is contained in $L_0\minter\partial K$, it follows that $L_0\minter K$ coincides with one of $\left\{(0,0)\right\}$, $\left\{0\right\}\times[0,2C\rho]$, $\left\{0\right\}\times[-2C\rho,0]$ or $\left\{0\right\}\times]-2C\rho,2C\rho[$. The last two are excluded since $(0,-C\rho)\notin K$ and it thus remains to show that $L_0\minter K\neq\left\{(0,0)\right\}$. However, suppose the contrary. For $0<\delta<\opMin(1,C)\rho$, the point $x:=(0,\delta)$ does not lie in $K$. Let $y:=(y',t)$ be the closest point in $K$ to $x$. By Lemma \procref{CharacterisationOfClosestPoints}, $\msf{N}:=(x-y)/\|x-y\|$ is a supporting normal to $K$ at $x$. However, $t\leq\delta$, since otherwise, by convexity, the point $(\delta y'/t,\delta)$ also lies in $K$, but is closer to $x$ than $y$. In particular, $\langle\msf{N},-e_{n+1}\rangle=t-\delta<0$. This is absurd, and we conclude that $L_0\minter K$ coincides with $\left\{0\right\}\times[0,2C\rho]$. In particular, the set $B'_{C\rho}$ intersects $K$ non-trivially, and is therefore entirely contained in the interior of $K$, as desired.
\par
Since $B'_{-C\rho}$ is contained in $\Bbb{R}^{n+1}\setminus K$ and since $B'_{C\rho}$ is contained in the interior of $K$, it follows that for all $x'\in B'_{\rho}(0)$, there exists a unique point $f(x')\in]-C\rho,C\rho[$ such that $(x',f(x'))\in\partial K$. In particular, by Lemma \procref{LemmaGraphFunctionIsLipschitz}, for all $x',y'\in B_\rho'(0)$,
$$
\left|f(x') - f(y')\right| \leq C\|x' - y'\|,
$$
so that the function $f$ is $C$-Lipschitz.
\par
Finally, denote $\opCyl:=B'_\rho(0)\times]-2C\rho,2C\rho[$ and $\hat{K}:=K\minter\opCyl$, and denote the graph of $f$ by $\opGr(f)$. Observe that $\opCyl\setminus\opGr(f)$ consists of two connected components. Furthermore, since $\hat{K}\minter\opCyl=\opGr(f)$, it follows that the set $\hat{K}\minter(\opCyl\setminus\opGr(f))$ is both open and closed in $\opCyl\setminus\opGr(f)$. Thus, since $B_{C\rho}\subseteq\hat{K}$ and $B_{-C\rho}\subseteq\opCyl\setminus\hat{K}$, it follows that $\hat{K}$ coincides with the closure of the connected component of $\opCyl\setminus\opGr(f)$ lying above $\opGr(f)$. That is,
$$
\hat{K} = \left\{(x',t)\ |\ t\geq f(x')\right\}.
$$
Now choose $x',y'\in B_r'(0)$ and, for all $s\in[0,1]$, denote $x'_s := (1-s)x' + sy'$ and $t_s := (1-s)f(x') + sf(y')$. By convexity, since $(x',f(x'))$ and $(y',f(y'))$ are both elements of $\hat{K}$, for all $s\in[0,1]$, so too is $(x'_s,t_s)$, so that
$$
f((1-s)x' + sy') \leq (1-s)f(x') + sf(y').
$$
Since $x',y'\in B_r'(0)$ are arbitrary, we conclude that $f$ is convex, and this completes the proof.\qed
\newsubhead{Convex hulls}
Let $X$ be a subset of $\Bbb{R}^{n+1}$. We define the \gloss{convex hull} of $X$ to be the intersection of all open, convex subsets of $\Bbb{R}^{n+1}$ containing $X$. We denote this set by $\opConv(X)$. Observe, in particular, that $\opConv(X)$ is also convex.
\proclaim{Lemma \nextprocno}
\noindent If $K$ is a convex set, then $\overline{K}$ and $K^o$ are also convex.
\endproclaim
\proclabel{ClosureOfConvexIsConvex}
\proof Choose $x,y\in\overline{K}$. Let $(x_m)_\minn$, $(y_m)_\minn$ be sequences of points in $K$ converging to $x$ and $y$ respectively. Choose $t\in[0,1]$. By convexity, $(1-t)x_m + ty_m\in K$ for all $m$. Taking limits, it follows that $(1-t)x + ty\in\overline{K}$. Since $x,y\in\overline{K}$ and $t\in[0,1]$ are arbitrary, we conclude that $\overline{K}$ is convex, as desired.
\par
Now choose $x,y\in K^o$, choose $\delta>0$ such that $B_\delta(x),B_\delta(y)\subseteq K$ and choose $t\in[0,1]$. By convexity, for all $z\in B_\delta(0)$,
$$
(1-t)x + ty + z = (1-t)(x+z) + t(y+z) \in K.
$$
It follows that $B_\delta((1-t)x + ty)\subseteq K$ and so $(1-t)x + ty\in K^o$. Since $x,y\in K^o$ and $t\in[0,1]$ are arbitrary, we conclude that $K^o$ is convex, as desired.\qed
\proclaim{Lemma \nextprocno}
\noindent If $X$ is compact, then $\opConv(X)$ is compact.
\endproclaim
\proclabel{ConvexHullOfCompactIsCompact}
\proof We first show that $\partial\opConv(X)\subseteq\opConv(X)$. Indeed, suppose the contrary, and choose $x\in\partial\opConv(X)\setminus\opConv(X)$. By definition, there exists an open, convex set $K$ such that $X\subseteq K$ but $x\notin K$. However, since $\opConv(X)\subseteq K$, in particular, $\partial\opConv(X)\subseteq\overline{K}$, so that $x\in\partial K$. By Lemma \procref{ClosureOfConvexIsConvex}, $\overline{K}$ is convex, and so, by Lemma \procref{SupportingNormalForClosedConvexSets}, there exists a supporting normal $\msf{N}$ to $\overline{K}$ at $x$. By definition of supporting normals, for all $y\in X\subseteq K\subseteq\overline{K}$, $\langle y-x,\msf{N}\rangle<0$. Thus, by compactness, there exists $\delta>0$ such that for all $y\in X$, $\langle y-x,\msf{N}\rangle<-\delta$. Now define
$$
K' := \left\{y\in K\ |\ \langle y-x,\msf{N}\rangle < -\delta\right\}.
$$
Since $K'$ is open and convex, and since $X\subseteq K'$, it follows that $\opConv(X)\subseteq K'$. However, $x\notin\overline{K}'\supseteq\partial\opConv(X)$. This is absurd, and it follows that $\partial\opConv(X)\subseteq\opConv(X)$, as desired. In particular, $\opConv(X)$ is closed. Finally, since $X$ is bounded, there exists $R>0$ such that $X\subseteq B_R(0)$. Since $B_R(0)$ is open and convex, it follows that $\opConv(X)\subseteq B_R(0)$, and is therefore also bounded. It follows by the Heine-Borel Theorem that $\opConv(X)$ is compact, as desired.\qed
\proclaim{Lemma \nextprocno}
\noindent If $X$ is open, then $\opConv(X)$ is open.
\endproclaim
\proclabel{ConvexHullOfOpenIsOpen}
\proof Choose $x\in\partial\opConv(X)$. By Lemma \procref{ClosureOfConvexIsConvex}, $\overline{\opConv}(X)$ is convex. By Lemma \procref{SupportingNormalForClosedConvexSets}, there exists a supporting normal $\msf{N}$ to $\overline{\opConv}(X)$ at $x$. Since $X$ is open, for all $y\in X\subseteq\opConv(X)\subseteq\overline{\opConv}(X)$, $\langle y-x,\msf{N}\rangle<0$. Thus, if we define
$$
K := \left\{y\in \Bbb{R}^{n+1}\ |\ \langle y-x,\msf{N}\rangle < 0\right\},
$$
then $X\subseteq K$ and since $K$ is open and convex, $\opConv(X)\subseteq K$. In particular, $x\in\Bbb{R}^{n+1}\setminus K\subseteq\Bbb{R}^{n+1}\setminus\opConv(X)$, and since $x\in\partial\opConv(X)$ is arbitrary, we conclude that $\partial\opConv(X)\minter\opConv(X)=\emptyset$, so that $\opConv(X)$ is open, as desired.\qed
\proclaim{Lemma \nextprocno}
\noindent If $K$ is compact and convex, then $\opConv(K)=K$.
\endproclaim
\proclabel{ConvexHullOfConvexIsTheSame}
{\sl\noindent{\bf Remark:\ }Observe that the analogous result for open convex sets follows immediately from our definition of the convex hull.}
\medskip
\proof Choose $x\in\Bbb{R}^{n+1}\setminus K$. Let $y\in K$ be a point minimising distance to $x$, and define $\msf{N} = (x - y)/\|x-y\|$. By Lemma \procref{CharacterisationOfClosestPoints}, $\msf{N}$ is a supporting normal to $K$ at $y$. By definition of supporting normals, for all $z\in K$, $\langle z - y,\msf{N}\rangle\leq 0$, and so
$$
\langle z - x,\msf{N}\rangle = \langle (z-y) + (y-x),\msf{N}\rangle \leq -\|x-y\| < 0.
$$
Define
$$
K' := \left\{z\in\Bbb{R}^{n+1}\ |\ \langle z - x,\msf{N}\rangle < 0\right\}.
$$
Since $K'$ is open and convex and since $K\subseteq K'$, it follows that $\opConv(K)\subseteq K'$. In particular, $x\notin\opConv(K)$. Since $x\in\Bbb{R}^{n+1}\setminus K$ is arbitrary, we conclude that $\Bbb{R}^{n+1}\setminus K\subseteq\Bbb{R}^{n+1}\setminus\opConv(K)$, so that $\opConv(K)\subseteq K$, and since $K$ is trivially contained in $\opConv(K)$, we conclude that the two sets coincide, as desired.\qed
\newsubhead{The local geodesic property}
We define an open straight-line segment in $\Bbb{R}^{n+1}$ to be any set $\Gamma$ of the form
\subheadlabel{LocalGeodesicProperty}
$$
\Gamma = \left\{ x + ty\ |\ a < t < b \right\},
$$
where $x$ and $y$ are points in $\Bbb{R}^{n+1}$ and $a<b$ are real numbers. Let $K$ be a compact, convex subset of $\Bbb{R}^{n+1}$ and let $x$ be any point of $K$. We say that $K$ satisfies the \gloss{local geodesic property} at $x$ whenever there exists an open straight-line segment $\Gamma$ such that $x\in\Gamma\subseteq K$. Every interior point of $K$ trivially satisfies the local geodesic property.
\proclaim{Lemma \nextprocno}
\noindent Let $K$ be a compact, convex subset of $\Bbb{R}^{n+1}$. Let $x$ be a boundary point of $K$ and let $\Gamma$ be an open straight-line segment such that $x\in\Gamma\subseteq K$. If $\msf{N}$ is a supporting normal to $K$ at $x$ then $\Gamma$ is contained in the hyperplane passing through $x$ normal to $\msf{N}$. In particular, $\Gamma$ is contained in the boundary of $K$.
\endproclaim
\proclabel{GBPLineIsContainedInBoundary}
\proof By definition, there exists $y\in\Bbb{R}^{n+1}$ and real numbers $a<0<b$ such that
$$
\Gamma = \left\{ x + ty\ |\ a<t<b\right\}.
$$
By definition of supporting normals, for all $t\in]a,b[$,
$$
t\langle y,\msf{N}\rangle = \langle (x+ty) - x,\msf{N}\rangle \leq 0.
$$
It follows that $\langle y,\msf{N}\rangle = 0$, and so $\langle (x+ty) - x,\msf{N}\rangle =0$ for all $t\in]a,b[$ as desired.\qed
\medskip
The local geodesic property characterises convex hulls in the following sense.
\proclaim{Theorem \nextprocno}
\noindent Let $K$ be a compact, convex set, let $X$ be a subset of $\partial K$, and let $Y$ be the set of all points of $\partial K$ satisfying the local geodesic property. If $X\munion Y$ is closed, then $Y\subseteq\opConv(X)$.
\endproclaim
\proclabel{GBPBdyPointsAreConvexHull}
\proof We prove this by induction on the dimension of the ambient space. First suppose that $n=1$. In particular, $K=:[a,b]$ is a compact interval. We claim that $\left\{a,b\right\}\subseteq X$. Indeed, observe that $]a,b[\subseteq Y$. Thus, since $X\munion Y$ is closed, it follows that $[a,b]\subseteq X\munion Y$. Since neither $a$ nor $b$ is an element of $Y$, it follows that $\left\{a,b\right\}\subseteq X$, as asserted. In particular, $K=[a,b]=\opConv(X)$, as desired.
\par
Now consider an ambient space of arbitrary dimension greater than $1$. Choose $y\in Y$ and let $H$ be a supporting hyperplane to $K$ at $y$. Denote $K':=K\minter H$, $Y':=Y\minter \partial K'$ and $X':=X\minter\partial K'$. Observe that $K'$ is a compact, convex subset of $H$ and that $X'\munion Y'$ is closed. We claim that $Y'$ coincides with the set of all boundary points of $K'$ satisfying the local geodesic property. Indeed, if $z\in Y'$, then there exists an open straight-line segment $\Gamma$ such that $z\in\Gamma\subseteq K$. Since $H$ is a supporting hyperplane to $K$, by Lemma \procref{GBPLineIsContainedInBoundary}, $\Gamma\subseteq H$. In particular, $\Gamma\subseteq K'$ and so $K'$ also satisfies the local geodesic property at $z$. Conversely, if $z$ is a boundary point of $K'$ and if $K'$ satisfies the local geodesic property at $z$, then $z$ is also a boundary point of $K$ and $K$ also satisfies the local geodesic property at $z$. The assertion follows and we conclude by the inductive hypothesis that $Y'\subseteq\opConv(X')\subseteq\opConv(X)$.
\par
If $y\in Y'$, then $y\in\opConv(X)$, and we are done. Otherwise, suppose that $y\in Y\setminus Y'$. That is, $y$ lies in the interior of $K'$. Let $V$ be any vector in $H$. Define $\gamma:\Bbb{R}\rightarrow H$ by $\gamma(t)=y + tV$. Denote $I=\gamma^{-1}(K')$. Since $K'$ is compact and convex, $I$ is a compact interval. Trivially, for all $t\in I^o$, $K$ satisfies the local geodesic property at $\gamma(t)$. That is, $I^o\subseteq\gamma^{-1}(Y)\subseteq\gamma^{-1}(Y\munion X)$, so that taking closures yields, $I\subseteq\gamma^{-1}(Y\munion X)$. Since, in addition, $\partial I\subseteq\gamma^{-1}(\partial K')$, it follows that $\partial I\subseteq \gamma^{-1}(Y'\munion X')\subseteq\gamma^{-1}(\opConv(X'))$. It follows that $y\in\opConv(X')$ and since $y\in Y$ is arbitrary, we conclude that $Y'\subseteq\opConv(X')$, and the result now follows by induction.\qed
\medskip
Conversely, we have
\proclaim{Theorem \nextprocno}
\noindent If $X$ is a compact subset of $\Bbb{R}^{n+1}$, then $\opConv(X)$ satisfies the local geodesic property at every point of $\opConv(X)\setminus X$.
\endproclaim
\proclabel{ConvexHullHasLGP}
\proof We prove this by induction on the dimension. The result trivially holds when the ambient space is $1$-dimensional. Now choose $x\in\opConv(X)\setminus X$. If $x$ is an interior point of $\opConv(X)$, then we are done. We therefore assume that $x$ is a boundary point of $\opConv(X)$. By Lemma \procref{ExistenceOfSupportingNormals}, there exists a supporting normal $\msf{N}$ to $\opConv(X)$ at $x$. Let $H$ be the hyperplane normal to $\msf{N}$ passing through $x$, and denote $X':=H\minter X$. We claim that $\opConv(X')=\opConv(X)\minter H$. Indeed, since $X'$ is compact, by Lemma \procref{ConvexHullOfCompactIsCompact}, so too is $\opConv(X')$. Choose $x'\in H\setminus\opConv(X')$ and let $y'\in H$ be a point in $\opConv(X')$ minimising distance to $x'$. Denote $\msf{N}':=(x'-y')/\|x'-y'\|$. By Lemma \procref{CharacterisationOfClosestPoints}, $\msf{N}'$ is a supporting normal to $\opConv(X')$ at $y'$. In particular, for all $z'\in X'\subseteq\opConv(X')$,
$$
\langle z' - y',\msf{N}'\rangle \leq 0,
$$
so that
$$
\langle z' - x',\msf{N}'\rangle = \langle (z' - y') + (y' - x'),\msf{N}'\rangle \leq -\|x'-y'\| < 0.
$$
However, by definition of $\msf{N}$ and $H$, for all $z\in X\setminus X' = X\setminus H\subseteq\opConv(X)\setminus H$,
$$
\langle z - x',\msf{N}\rangle = \langle z - x,\msf{N}\rangle + \langle x - x',\msf{N}\rangle < 0.
$$
Combining these relations yields, for all $z\in X$,
$$
\langle z - x',(\msf{N} + \msf{N}')\rangle < 0.
$$
Define $K\subseteq\Bbb{R}^{n+1}$ by,
$$
K := \left\{z\ |\ \langle z - x',(\msf{N}+\msf{N}')\rangle < 0 \right\}.
$$
Since $K$ is open and convex and since $X\subseteq K$, it follows that $\opConv(X)\subseteq K$. In particular, since $x'$ is not an element of $K$, it is not an element of $\opConv(X)$ either. Since $x'\in H\setminus\opConv(X')$ is arbitrary, we conclude that $\opConv(X)\minter H\subseteq\opConv(X')$. Conversely, let $K$ be an open, convex set containing $X$. Then $K\minter H$ is also convex and relatively open. Since $K\minter H$ contains $X'$, by definition, it also contains $\opConv(X')$. Upon taking the intersection over all such open sets, we conclude that $\opConv(X')\subseteq\opConv(X)\minter H$, and the two sets therefore coincide as asserted. It follows by the inductive hypothesis that $\opConv(X)\minter H = \opConv(X')$ satisfies the local geodesic property at $x$ and therefore so too does $\opConv(X)$. This completes the proof.\qed
\medskip
We introduce an alternative characterisation of the local geodesic property which will be of use in the sequel. Let $X\subseteq\Sigma^n$ be any closed subset. We say that $X$ is \gloss{strictly contained in a hemisphere} whenever there exists a unit vector $\msf{N}\in\Sigma^n$ such that for all $\msf{N}'\in X$,
$$
\langle\msf{N},\msf{N}'\rangle < 0.
$$
\proclaim{Lemma \nextprocno}
\noindent Let $K$ be a compact, convex subset of $\Bbb{R}^{n+1}$. Let $x$ be a point in $K$ and suppose that $K$ satisfies the local geodesic property at $x$. Then, for all sufficiently small $r>0$, $K\minter\partial B_r(x)$ is not strictly contained in a hemisphere.
\endproclaim
\proclabel{GBPMeansNotStrictlyContainedInAHemisphere}
\proof Let $\Gamma$ be the open straight-line segment passing through $x$ contained in $K$. By definition, there exists a point $y\in\Bbb{R}^{n+1}$ and real numbers $a<0<b$ such that
$$
\Gamma = \left\{ x + ty\ | a<t<b\right\}.
$$
Choose $r<\opMin(-a,b)$. Then $x\pm ry\in K\minter\partial B_r(x)$. In particular, if there exists $\msf{N}\in\Sigma^n$ such that $\langle (x\pm ry) - x,\msf{N}\rangle < 0$, then, $\pm\langle y,\msf{N}\rangle < 0$. This is absurd, and the assertion follows.\qed
\proclaim{Lemma \nextprocno}
\noindent If $X\subseteq\Sigma^n$ is a closed subset not strictly contained in a hemisphere, then $0$ is an element of $\opConv(X)$.
\endproclaim
\proclabel{ZeroIsInConvexHull}
\proof Suppose the contrary. By Lemma \procref{ConvexHullOfCompactIsCompact}, $\opConv(X)$ is compact. Let $x\in\opConv(X)$ be the point minimising distance in $\opConv(X)$ to $0$ and denote $\msf{N}:=-x/\|x\|$. By Lemma \procref{CharacterisationOfClosestPoints}, $\msf{N}$ is a supporting normal to $\opConv(X)$ at $x$. Thus, for all $y\in X\subseteq\opConv(X)$, $\langle y-x,\msf{N}\rangle\leq 0$, and so
$$
\langle y,\msf{N}\rangle = \langle (y-x) + x,\msf{N}\rangle \leq -\|x\| < 0.
$$
Since $y\in X$ is arbitrary, we conclude that $X$ is strictly contained in a hemisphere. This is absurd, and so $0$ is an element of $\opConv(X)$ as desired.\qed
\proclaim{Lemma \nextprocno}
\noindent Let $K$ be a compact, convex subset of $\Bbb{R}^{n+1}$. Let $x$ be a boundary point of $K$ and let $\msf{N}$ be a supporting normal to $K$ at $x$. Suppose that $K$ does not satisfy the local geodesic property at $x$. Then, for all $r>0$ there exists $\msf{N}'$, which we may choose as close to $\msf{N}$ as we wish, with the property that for all $y\in K\minter (B_r(x))^c$,
$$
\langle y - x,\msf{N}'\rangle < 0.
$$
\endproclaim
\proclabel{NoLGPMeansPointyBit}
{\sl{\bf\noindent Remark:\ }In particular, using the terminology of links which we introduce in Section \subheadref{Links}, below, for all $r>0$, the closure of $\Cal{L}_r(x;K)$ is strictly contained in a hemisphere.}
\medskip
\proof Upon applying an affine isometry, we may suppose that $x=0$. Choose $r>0$. Since $K$ does not satisfy the local geodesic property at $0$, by Lemma \procref{ConvexHullHasLGP}, $0$ does not lie in the convex hull of $K\minter\partial B_r(0)$. By Lemma \procref{ZeroIsInConvexHull}, $K\minter\partial B_r(0)$ is strictly contained in a hemisphere. There therefore exists a unit vector $\msf{N}'\in\Sigma^n$ such that for all $y\in K\minter\partial B_{r}(0)$,
$$
\langle y,\msf{N}'\rangle < 0.
$$
For all $\epsilon>0$, denote $\msf{N}_\epsilon := (\msf{N} + \epsilon\msf{N}')/\|\msf{N} + \epsilon\msf{N}'\|$. We claim that $\msf{N}_\epsilon$ has the desired properties for all $\epsilon>0$. Indeed, for all $\epsilon>0$, and for all $y\in K\minter\partial B_r(0)$, we obtain
$$
\|\msf{N} + \epsilon\msf{N}'\|\langle y,\msf{N}_\epsilon\rangle = \langle y ,\msf{N}\rangle + \epsilon\langle y ,\msf{N}'\rangle\leq \epsilon\langle y,\msf{N}'\rangle< 0.$$
However, if $z\in K\minter(B_r(x))^c$, then, by convexity, $z=sy$ for some $y\in K\minter\partial B_r(0)$ and some $s\geq 1$, so that
$$
\langle z,\msf{N}_\epsilon\rangle = s\langle y,\msf{N}_\epsilon\rangle < 0,
$$
as desired.\qed
\newsubhead{Interior a-priori bounds}
We now return to the framework of Section \headref{TheCNSMethod}, and consider smooth, convex functions $f:\overline{\Omega}\rightarrow]-\infty,0]$ which are solutions of \eqnref{BasicEquation}. We develop a-priori estimates that will allow us in the following section to describe the local geometric structure of singularities of uniform limits of sequences of such functions. We achieve this by once again using the maximum principal via an argument that dates back to Pogorelov (c.f. \cite{Pog}, but see also \cite{Calabi} and \cite{ShengUrbasWang}). We will use the notation of Section \subheadref{BoundaryToGlobal}, and we first complement Lemma \procref{ThirdBarrierFunction} with the following two estimates.
\proclaim{Lemma \nextprocno}
\noindent There exists $C>0$ which only depends on $\|\phi\|_0$ such that
$$
\Cal{L}_f\opLog(-f) \geq -C(-f)^{-1} - B^{ij}(\partial_i\opLog(-f))(\partial_j\opLog(-f)).
$$
\endproclaim
\proclabel{FourthBarrierFunction}
\proof By Corollary \procref{HomogeneityOfF},
$$
DF(D^2f)(D^2f) = F(D^2f).
$$
Furthermore, since $\phi$ is positive and since $G$ is convex,
$$
\phi DG(Df)(Df) \geq \phi G(Df) - \phi G(0).
$$
Thus, by definition of $\Cal{L}_f$,
$$
\Cal{L}_f f \leq F(D^2f) - \phi G(Df) + \phi G(0) = \phi G(0).
$$
There therefore exists $C>0$ which only depends on $\|\phi\|_0$ such that,
$$
\Cal{L}_f(-f) \geq -C.
$$
Thus, by Lemma \procref{DerivativeOfLogarithm},
$$
\Cal{L}_f\opLog(-f) \geq -C(-f)^{-1} - B^{ij}(\partial_i\opLog(-f))(\partial_j\opLog(-f)),
$$
as desired.\qed
\proclaim{Lemma \nextprocno}
\noindent There exists $C>0$ which only depends on $\|\phi\|_1$ and $\|f\|_1$ such that
$$
\Cal{L}_f\|Df\|^2 \geq -C + \frac{2\phi}{n}\lambda_n(f).
$$
\endproclaim
\proclabel{FifthBarrierFunction}
\proof By the product rule, for all $i$, using the summation convention,
$$\eqalign{
\partial_i\|Df\|^2 &= 2f_{ik}f_k\cr
\partial_i\partial_j\|Df\|^2 &= 2f_{ijk}f_k + 2f_{ik}f_{jk}.\cr}
$$
Thus, by definition of $B^{ij}$, and using the symmetry of the derivatives of $f$,
$$\eqalign{
\Cal{L}_f\|Df\|^2&= 2f_kDF(D^2f)(D^2f_k) + 2B^{ij}f_{ik}f_{jk} - 2f_k\phi DG(Df)(Df_k)\cr
&=2f_k\Cal{L}_ff_k + \frac{2}{n}F(D^2f)f_{kk}.\cr}
$$
Since $f$ satisfies \eqnref{BasicEquation}, since $\phi$ is positive, since $G\geq 1$ and since $D^2f$ is positive definite,
$$
\frac{2}{n}F(D^2f)f_{kk} = \frac{2\phi}{n} G(Df)\opTr(D^2f) \geq \frac{2\phi}{n}\lambda_n(f).
$$
On the other hand, differentiating \eqnref{BasicEquation} once in the $e_k$ direction yields
$$
\Cal{L}_f f_k = \phi_k G(Df).
$$
Combining these relations, we obtain
$$
\Cal{L}_f\|Df\|^2 \geq 2 f_k\phi_k G(Df) + \frac{2\phi}{n}\lambda_n(f).
$$
There therefore exists $C>0$ which only depends on $\|\phi\|_1$ and $\|f\|_1$ such that
$$
\Cal{L}_f\|Df\|^2 \geq -C + \frac{2\phi}{n}\lambda_n(f),
$$
as desired.\qed
\proclaim{Theorem \nextprocno}
\noindent There exists $C>0$ which only depends on $\|\phi\|_2$, $\|f\|_1$ and $\minf_{x\in\overline{\Omega}}\phi(x)$ such that
$$
\msup_{x\in\overline{\Omega}}\left|f(x)\right|^2\|D^2f(x)\|\leq C.
$$
\endproclaim
\proclabel{PogorelovsInteriorSecondOrderBound}
\proof For $\epsilon\in]0,1[$, define the function $\psi_\epsilon:\Omega\rightarrow\Bbb{R}$ by
$$
\psi_\epsilon := \mu_n(f) + 2\opLog(-f) + \epsilon\|Df\|^2.
$$
It suffices to obtain a-priori bounds for $\psi_\epsilon$ for some $\epsilon>0$. Fix $\epsilon$ and observe that $\psi_\epsilon$ is continuous and tends to $-\infty$ near the boundary of $\Omega$. It therefore attains its maximum at some interior point $x\in\Omega$. Upon applying an affine isometry, we may suppose that $x=0$ and that $e_1,...,e_n$ are the eigenvectors of $D^2f(0)$ corresponding to the eigenvalues $\lambda_1,...,\lambda_n$ respectively. By Lemmas \procref{ThirdBarrierFunction}, \procref{FourthBarrierFunction} and \procref{FifthBarrierFunction}, there exists $C_1>0$, which only depends on $\|\phi\|_2$, $\|f\|_1$ and $\minf_{x\in\overline{\Omega}}\phi(x)$ such that
$$\eqalign{
\Cal{L}_f\psi_\epsilon &\geq -\frac{C_1}{(-f)^2}+\frac{2\epsilon\phi}{n}\lambda_n(f) - \frac{2}{n\lambda_n}(\partial_n(\psi_\epsilon - 2\opLog(-f) - \epsilon\|Df\|^2))^2\cr
&\qquad +B^{ij}\partial_i(\psi_\epsilon - 2\opLog(-f) - \epsilon\|Df\|^2)\partial_j(\psi_\epsilon - 2\opLog(-f) - \epsilon\|Df\|^2)\cr
&\qquad -2B^{ij}\partial_i\opLog(-f)\partial_j\opLog(-f)\cr}
$$
in the weak sense. Let $\alpha:\Omega\rightarrow\Bbb{R}$ be such that $\alpha\leq\psi_\epsilon$ and $\alpha(0)=\psi_\epsilon(0)$. In particular, $\alpha$ attains its maximum at $0$, and so $(\partial_n\alpha)(0)=0$ so that, at the origin,
$$\eqalign{
\Cal{L}_f\alpha &\geq -\frac{C_1}{(-f)^2}+\frac{2\epsilon\phi}{n}\lambda_n(f) - \frac{2}{n\lambda_n}(2\partial_n\opLog(-f) + \epsilon\partial_n\|Df\|^2)^2\cr
&\qquad +B^{ij}\partial_i(2\opLog(-f) + \epsilon\|Df\|^2)\partial_j(2\opLog(-f) + \epsilon\|Df\|^2)\cr
&\qquad -2B^{ij}\partial_i\opLog(-f)\partial_j\opLog(-f).\cr}
$$
Let $A_1$ and $A_2$ denote respectively the third term and the sum of the last two terms on the right-hand side of the above equation. By definition of $\lambda_n$,
$$
A_1=-\frac{2}{n\lambda_n}\left(\frac{2}{(-f)}(-f_n) + 2\epsilon\lambda_n f_n\right)^2.
$$
However, by Lemma \procref{LowerBoundsOfTraceOfHessian}, by definition of $G$, and since $f$ solves \eqnref{BasicEquation},
$$
\lambda_n \geq \frac{1}{n}\opTr(D^2f) \geq F(D^2f) = \phi G(Df) \geq \phi.
$$
There therefore exists $C_2>0$, which only depends on $\|f\|_1$ and $\inf_{x\in\overline{\Omega}}\phi(x)$ such that
$$
A_1 \geq -\frac{C_2}{(-f)^2}-C_2\epsilon^2\lambda_n.
$$
Since $B^{ij}$ is symmetric and positive definite,
$$
A_2 \geq 4\epsilon B^{ij}(\partial_i\opLog(-f))(\partial_j\|Df\|^2).
$$
Thus, by definition of $B^{ij}$, and since $f$ solves \eqnref{BasicEquation},
$$
A_2 \geq \frac{-8\epsilon}{n(-f)}F(D^2f)\|Df\|^2 = \frac{-8\epsilon\phi}{n(-f)}G(Df)\|Df\|^2 \geq \frac{-8\epsilon\phi}{n(-f)}\|Df\|^2.
$$
In particular, there exists $C_3>0$, which only depends on $\|\phi\|_0$ and $\|f\|_1$ such that
$$
A_2 \geq -\frac{C_3}{(-f)^2}.
$$
Combining these relations, we conclude that there exists $C_4>0$, which only depends on $\|\phi\|_2$, $\|f\|_1$ and $\minf_{x\in\overline{\Omega}}\phi(x)$ such that, at the origin,
$$
\Cal{L}_f\alpha \geq -\frac{C_4}{(-f)^2} + \epsilon\left(\frac{2\phi}{n} - C_4\epsilon\right)\lambda_n.
$$
In particular, for $\epsilon>0$ sufficiently small, there exists $C_5>0$ which only depends on $\|\phi\|_2$, $\|f\|_1$, $\minf_{x\in\overline{\Omega}}\phi(x)$ and $\epsilon$ such that, at the origin,
$$
\Cal{L}_f\alpha \geq \frac{1}{(-f)^2}\left(-C_5 + \frac{1}{C_5} e^{\alpha}\right).
$$
However, since $\alpha$ attains its maximum at the origin, $\Cal{L}_f\alpha\leq 0$, so that
$$
\psi_\epsilon(0) = \alpha(0) \leq 2\opLog(C_5),
$$
as desired.\qed
\newsubhead{The structure of singularities}
We now describe the local structure of singularities that arise upon taking limits. We first require some preliminary results.
\proclaim{Lemma \nextprocno}
\noindent Let $K$ be a compact, convex subset of $\Bbb{R}^{n+1}$. If $K$ has non-trivial interior, then $\Cal{N}(x)$ is strictly contained in a hemisphere for every boundary point $x$ of $K$.
\endproclaim
\proclabel{NonTrivialInteriorMeansSupportingNormalsContainedInHypersphere}
\proof Let $y$ be an interior point of $K$. Let $x$ be a boundary point of $K$. Denote $\msf{N}=(y-x)/\|y-x\|$. Choose $\msf{M}\in\Cal{N}(x)$. Since $y$ is an interior point of $K$, $\langle\msf{N},\msf{M}\rangle<0$, and since $\msf{M}\in\Cal{N}(x)$ is arbitrary, conclude that $\Cal{N}(x)$ is strictly contained in a hemisphere as desired.\qed
\proclaim{Lemma \nextprocno}
\noindent Let $K$ be a compact, convex subset of $\Bbb{R}^{n+1}$. Let $x$ be a boundary point of $K$. If $\Cal{N}(x)$ is strictly contained in a hemisphere, then there exists a supporting normal $\msf{N}$ to $K$ at $x$ such that $\langle\msf{N},\msf{M}\rangle>0$ for all $\msf{M}\in\Cal{N}(x)$.
\endproclaim
\proclabel{OptimalSupportingNormal}
\proof Upon applying a linear isometry, we may suppose that $\langle\msf{N},e_{n+1}\rangle<0$ for all $\msf{N}\in\Cal{N}(x)$. Observe that $\Cal{N}(x)$ is closed. If $-e_{n+1}\in\Cal{N}(x)$, then we are done. Otherwise, suppose that $-e_{n+1}\notin\Cal{N}(x)$. Define
$$
\opCone(\Cal{N}(x)) := \left\{t\msf{N}\ |\ t\in[0,\infty[,\ \msf{N}\in\Cal{N}(x)\right\}.
$$
Observe that $\opCone(\Cal{N}(x))$ is closed and convex. Moreover, $-e_{n+1}\notin\opCone(\Cal{N}(x))$. Let $y\in\opCone(\Cal{N}(x))$ minimise distance in $\opCone(\Cal{N}(x))$ to $-e_{n+1}$. We claim that $\msf{N}:=y/\|y\|$ has the desired property. Indeed, denote $\hat{N}:=-(e_{n+1}+y)/\|e_{n+1}+y\|$. By Lemma \procref{CharacterisationOfClosestPoints}, $\hat{N}$ is a supporting normal to $\opCone(\Cal{N}(x))$ at $y$. Thus, for all $t\in[0,\infty[$, since $ty\in\opCone(\Cal{N}(x))$,
$$
\langle ty - y,\hat{\msf{N}}\rangle \leq 0,
$$
and differentiating this relation at $t=1$ yields $\langle y,\hat{\msf{N}}\rangle=0$. Now choose $\msf{N}\in\Cal{N}(x)\subseteq\opCone(\Cal{N}(x))$. Using the fact that $\langle\msf{M},-e_{n+1}\rangle>0$, we have
$$\eqalign{
0&\geq \langle \msf{M} - y,\hat{\msf{N}}\rangle\cr
&= \langle \msf{M},\hat{\msf{N}}\rangle\cr
&= \|e_n + y\|^{-1}\langle \msf{M}, -e_{n+1} - y\rangle\cr
&> \|e_n + y\|^{-1}\|y\|\langle \msf{M}, -y/\|y\|\rangle,\cr
}$$
so that $\langle\msf{M},\msf{N}\rangle>0$, and since $\msf{M}\in\Cal{N}(x)$ is arbitrary, the result follows.\qed
\proclaim{Theorem \nextprocno}
\noindent Let $(K_m)_\minn$ and $K_\infty$ be compact, convex subsets of $\Bbb{R}^{n+1}$. Suppose that $K_\infty$ has non-trivial interior and that $(K_m)_\minn$ converges to $K_\infty$ in the Hausdorff sense. Let $k>0$ be a real number, let $U$ be an open subset of $\Bbb{R}^{n+1}$ and suppose that for all $m$, $(\partial K_m)\minter U$ is smooth with constant gaussian curvature equal to $k$. If $y\in(\partial K_\infty)\minter U$ then either
\medskip
\myitem{(1)} there exists $r>0$ such that $(\partial K_\infty)\minter B_r(y)$ is smooth with constant gaussian curvature equal to $k$; or
\medskip
\myitem{(2)} $K_\infty$ satisfies the local geodesic property at $y$.
\endproclaim
\proclabel{StructureOfSingularities}
\proof Choose $y\in(\partial K_\infty)\minter U$ and suppose that $K_\infty$ does not satisfy the local geodesic property at $y$. Let $(y_m)_\minn$ be a sequence converging to $y$ such that $y_m\in\partial K_m$ for all $m$. Upon applying a convergent sequence of affine isometries, we may suppose that $y_m=0$ for all $m$. Since $K_\infty$ has non-trivial interior, by Lemma \procref{NonTrivialInteriorMeansSupportingNormalsContainedInHypersphere}, $\Cal{N}(0;K_\infty)$ is strictly contained in a hemisphere. In particular, by Lemma \procref{OptimalSupportingNormal}, we may suppose that $-e_{n+1}\in\Cal{N}(0;K_\infty)$ and that $\langle\msf{N},-e_{n+1}\rangle>3\opCos(\theta)$ for all other $\msf{N}\in\Cal{N}(0;K_\infty)$ and for some $\theta\in[0,\pi/2[$. Denote $C:=\opTan(\theta)$.
\par
By Lemma \procref{ConvergenceOfSupportingNormals}, upon extracting a subsequence, we may suppose that there exists $r>0$ such that $B_r(0)\subseteq U$ and for all $m$, for all $x\in(\partial K_m)\minter B_r(0)$, and, for all $\msf{N}\in\Cal{N}(x;K_m)$, $\langle\msf{N},-e_{n+1}\rangle>2\opCos(\theta)$. Denote $\rho:=r/\sqrt{1+4C^2}$. By Lemma \procref{NoLGPMeansPointyBit}, upon applying a small rotation, we may suppose that for all $x\in K_\infty\setminus B_{\rho/2}(0)$, $\langle x,-e_{n+1}\rangle<0$. Choosing this rotation sufficiently small, we may continue to assume that for all $m$, for all $x\in(\partial K_m)\minter B_r(0)$ and for all $\msf{N}\in\Cal{N}(x;K_m)$, $\langle \msf{N},-e_{n+1}\rangle > \opCos(\theta)$.
\par
By Theorem \procref{ExistenceOfGraphFunction}, for all $m$, there exists a convex, $C$-Lipschitz function, $f_m:B_\rho'(0)\rightarrow]-C\rho,C\rho[$ such that $f_m(0)=0$ and $(\partial K_m)\minter (B_\rho'(0)\times]-2C\rho,2C\rho[)$ coincides with the graph of $f_m$ over $B_\rho'(0)$. By the Arzela-Ascoli theorem, every subsequence of $(f_m)_\minn$ has a subsubsequence converging in the local uniform sense over $B_\rho'(0)$ to some limit $f_\infty'$ say. Since $(K_m)_\minn$ converges to $K_\infty$ in the Hausdorff sense, it follows that $f_\infty'=f_\infty$, and we conclude that $(f_m)_\minn$ converges in the local uniform sense over $B_\rho'(0)$ to $f_\infty$.
\par
Observe that $f_m$ is smooth for all $m<\infty$. Furthermore, there exists $\delta>0$ such that, $f_\infty(x')>4\delta$ for all $x'\in\partial B_{\rho/2}'(0)$. Since $(f_m)_\minn$ converges locally uniformly to $f_\infty$ over $B_\rho'(0)$, we may suppose that, for all $m$ and for all $x'\in B_{\rho/2}'(0)$, $f_m(x')>2\delta$. For all $m<\infty$, since $f_m$ is $C$-Lipschitz,
$$
\|f_m|_{\overline{B}_{\rho/2}'(0)}\|_1 \leq C(1+\rho/2).
$$
Thus, by Theorem \procref{PogorelovsInteriorSecondOrderBound}, there exists $C_2>0$ such that for all $m<\infty$, and for all $x\in B'_{\rho/2}(0)$,
$$
\left|2\delta - f_m(x)\right|^2\|D^2f_m(x)\| \leq C_2.
$$
Since $f_\infty(0)=0$, by continuity, there exists $s\in]0,\rho/2[$ such that $f_\infty(x')\leq \delta/2$ for all $x'\in\overline{B}_s'(0)$. Since $(f_m)_\minn$ converges locally uniformly to $f_\infty$, we may suppose that for all $m<\infty$ and for all $x\in\overline{B}_s'(0)$, $f_m(x')\leq\delta$, so that
$$
\|D^2f_m(x)\| \leq C_2/\delta^2.
$$
By the Krylov estimates (c.f Theorem \procref{Krylov}) and the Schauder estimates (c.f. Theorem \procref{Schauder}), for all $k\in\Bbb{N}$, there exists $C_k>0$ such that for all $m<\infty$,
$$
\|f_m|_{\overline{B}'_{s/2}(0)}\|_k \leq C_k.
$$
By the Arzela-Ascoli theorem, every subsequence of $(f_m)_\minn$ has a subsubsequence which converges in the $C^\infty$ sense over $\overline{B}'_{s/4}(0)$ to some limit $f_\infty'$ say. Since $(f_m)_\minn$ converges uniformly to $f_\infty$, it follows that $f_\infty'=f_\infty$. We conclude that $(f_m)_\minn$ converges to $f_\infty$ in the $C^\infty$ sense over $\overline{B}'_{s/4}(0)$. In particular, $f_\infty$ is smooth over $B'_{s/4}(0)$ and its graph has constant gaussian curvature equal to $k$. In other words, $(\partial K_\infty)\minter(B'_{s/4}(0)\times]-2C\rho,2C\rho[)$ is smooth and has constant gaussian curvature equal to $k$, as desired.\qed 
\newhead{Duality of Convex Sets}
Before analysing the general Plateau problem, it will be useful to continue our study of the elementary geometry of convex sets. In particular, we review the concept of duality for subsets of the sphere, showing how it is closely related to the concept of the convex hull, which we introduced in the preceeding chapter. We introduce the infinitesimal link of a given boundary point of a given compact, convex subset of $\Bbb{R}^{n+1}$ with non-trivial interior. This is defined to be an open subset of the sphere, and we show that it coincides with the dual of the set of supporting normal vectors to the convex set at that point. This allows us to prove the most important result of this chapter, namely Theorem \procref{SupportingNormalsOfIntersection}, which determines the supporting normal set  of the intersection of two given convex sets at any point on the boundary of this intersection.
\par
Although the results of this chapter are of use in the sequel, they are only tangential to the main flow of this text. In particular, Theorem \procref{SupportingNormalsOfIntersection}, although interesting, may be substituted by ad-hoc arguments in the relatively straightforward cases where it will be applied.
\newsubhead{Open half spaces and convex hulls}
Let $\msf{N}$ be a unit vector and let $t>0$ be a positive real number (possibly $+\infty$). We define the subset $H(\msf{N},t)$ of $\Bbb{R}^{n+1}$ by
\headlabel{WeakBarriers}
\subheadlabel{OpenHalfSpacesAndConvexHulls}
$$
H(\msf{N},t) := \left\{ x\ |\ \langle x, \msf{N}\rangle < t \right\},
$$
and we refer to this set as the \gloss{open half-space} normal to $\msf{N}$ of height $t$. Observe that this definition incorporates the degenerate case $\Bbb{R}^{n+1}=H(\msf{N},\infty)$.
\proclaim{Lemma \nextprocno}
\noindent If $K$ is an open, convex subset of $\Bbb{R}^{n+1}$, then $K=(\overline{K})^o$.
\endproclaim
\proclabel{InteriorOfClosureIsSame}
\proof Since $K$ is open, $K\subseteq (\overline{K})^o$. We now show that $(\overline{K})^o\subseteq K$. By Lemma \procref{ClosureOfConvexIsConvex}, $\overline{K}$ is convex. Choose $x\in(\overline{K})^o$. Without loss of generality, we may suppose that $x=0$. Choose $\delta>0$ such that $B_\delta(0)\subseteq\overline{K}$. Then $K\minter B_\delta(0)$ is a dense subset of $B_\delta(0)$. Upon applying a homothety, we may suppose that $\delta=2$. Choose $x_1,...,x_k\in B_1(0)$ such that
$$
\partial B_1(0) \subseteq \munion_{i=1}^k B_1(x_i).
$$
Since $K\minter B_\delta(0)$ is dense, upon perturbing $x_1,...,x_k$ if necessary, we may suppose that $x_i\in K$ for all $i$. Let $L$ be the convex hull of $\left\{x_1,...,x_k\right\}$. In particular, $L$ is a compact, convex subset of $K$. We claim that $0$ is an element of $L$. Indeed, suppose the contrary. Let $y\in L$ be the point minimising distance to $0$ and denote $\msf{N}:=-y/\|y\|$. By Lemma \procref{CharacterisationOfClosestPoints}, $\msf{N}$ is a supporting normal to $L$ at $y$. Thus, for all $z\in L$, $\langle z - y,\msf{N}\rangle \leq 0$. In particular, for all $1\leq i\leq k$,
$$
\langle x_i,\msf{N}\rangle = \langle x_i - y,\msf{N}\rangle + \langle y,\msf{N}\rangle \leq -\|y\| < 0,
$$
so that, for all $i$,
$$
\|x_i - \msf{N}\|^2 = \|x_i\|^2 - 2\langle x_i,\msf{N}\rangle + \|\msf{N}\|^2 > 2.
$$
However, by definition of $x_1,...,x_k$, there exists $1\leq i\leq k$ such that $\|x_i-\msf{N}\|^2<1$. This is absurd, and we conclude that $0\in L\subseteq K$ as asserted. Since $x\in (\overline{K})^o$ is arbitrary, it follows that $(\overline{K})^o\subseteq K$, and the two sets therefore coincide, as desired.\qed
\proclaim{Theorem \nextprocno}
\noindent For any subset $X$ of $\Bbb{R}^{n+1}$, the convex hull of $X$ coincides with the intersection of all open half-spaces containing $X$.
\endproclaim
\proclabel{ConvexHullIsIntersectionOfOpenHalfSpaces}
\proof Denote by $\hat{X}$ the intersection of all open half-spaces containing $X$. Since every open half-space is also convex, by definition of the convex hull, $\opConv(X)\subseteq\hat{X}$. We now show that $\hat{X}\subseteq\opConv(X)$. Indeed, choose $x\in\Bbb{R}^{n+1}\setminus\opConv(X)$. Let $K$ be an open, convex set such that $X\subseteq K$ and $x\notin K$. By Lemma \procref{ClosureOfConvexIsConvex}, $\overline{K}$ is also convex. We now have two cases to consider. Suppose first that $x\in\Bbb{R}^{n+1}\setminus\overline{K}$. Let $y$ be a point in $\overline{K}$ minimising distance to $x$ and denote $\msf{N} := (x-y)/\|x-y\|$. By Lemma \procref{SupportingNormalForClosedConvexSets}, $\msf{N}$ is a supporting normal to $\overline{K}$ at $y$, and so, for all $z\in K\subseteq\overline{K}$,
$$
\langle z - x,\msf{N}\rangle = \langle (z-y) + (y-x),\msf{N}\rangle \leq -\|x-y\| < 0.
$$
Now suppose that $x\in \overline{K}\setminus K$. By Lemma \procref{InteriorOfClosureIsSame}, $K=(\overline{K})^o$, and so $x\in\overline{K}\setminus(\overline{K})^o=\partial\overline{K}$.  Let $\msf{N}$ be a supporting normal to $\overline{K}$ at $x$. For all $z\in K$, $\langle z - x,\msf{N}\rangle\leq 0$, and, since $K$ is open, $\langle z - x,\msf{N}\rangle < 0$. In both cases, we conclude that $K\subseteq H(\msf{N},\langle\msf{N},x\rangle)$ so that, by definition, $\hat{X}\subseteq H(\msf{N},\langle\msf{N},x\rangle)$. In particular, $x\notin\hat{X}$, and since $x\notin\opConv(X)$ is arbitrary, we conclude that $\Bbb{R}^{n+1}\setminus\opConv(X)\subseteq\Bbb{R}^{n+1}\setminus\hat{X}$. Taking complements yields $\hat{X}\subseteq\opConv(X)$, and the two sets therefore coincide, as desired.\qed
\medskip
\noindent We also have the following complement of Lemma \procref{InteriorOfClosureIsSame}.
\proclaim{Lemma \nextprocno}
\noindent If $K$ is a compact, convex subset of $\Bbb{R}^{n+1}$ with non-trivial interior, then $K=\overline{K^o}$.
\endproclaim
\proclabel{ClosureOfInteriorIsSame}
\proof Since $K^o\subseteq K$ and since $K$ is closed, $\overline{K^o}\subseteq K$. Conversely, choose $x\in K$. Since $K$ has non-trivial interior, there exists $y\in K^o$. Choose $\delta>0$ such that $B_\delta(y)\subseteq K$. Bearing in mind that $K$ is convex, for all $t\in]0,1]$ and for all $z\in B_{t\delta}(0)$,
$$
(1-t)x + ty + z = (1-t)x + t(y+t^{-1}z) \in K.
$$
In other words, for all $t\in]0,1]$, $B_{t\delta}((1-t)x + ty)\subseteq K$, so that $(1-t)x+ty\in K^o$. It follows that $x\in\overline{K^o}$, and since $x\in K$ is arbitrary, we conclude that $K\subseteq\overline{K^o}$, and the two sets therefore coincide, as desired.\qed
\newsubhead{Convex subsets of the sphere}
Let $\msf{N}_0,\msf{N}_1\in\Sigma^{n}$ be points in the sphere. We say that $\msf{N}_0$ and $\msf{N}_1$ are \gloss{non-antipodal} whenever $\msf{N}_0+\msf{N}_1\neq 0$. In this case, we define the curve $\msf{N}:[0,1]\rightarrow\Sigma^n$ by
\subheadlabel{ConvexSubsetsOfTheSphere}
$$
\msf{N}(s) := \frac{(1-s)\msf{N}_0 + s\msf{N}_1}{\|(1-s)\msf{N}_0 + s\msf{N}_1\|}.
$$
We refer to $\msf{N}$ as the \gloss{great-circular arc} joining $\msf{N}_0$ to $\msf{N}_1$. This terminology is justified by the following result.
\proclaim{Lemma \nextprocno}
\noindent If $\msf{N}_0,\msf{N}_1\in\Sigma^n$ are distinct, non-antipodal points of the sphere, then there exists a unique great circle $C$ passing through $\msf{N}_0$ and $\msf{N}_1$. Moreover, if $\msf{N}$ is the great-circular arc joining $\msf{N}_0$ to $\msf{N}_1$ then, for all $s\in[0,1]$, $\msf{N}(s)$ is an element of $C$.
\endproclaim
\proclabel{GreatCircularArcs}
\proof Observe that every great circle in $\Sigma^n$ coincides with the intersection of $\Sigma^n$ with a plane in $\Bbb{R}^{n+1}$ containing the origin. Conversely, the intersection of any such plane with $\Sigma^n$ is a great circle. Now let $\msf{N}_0$ and $\msf{N}_1$ be distinct, non-antipodal points. In particular, they are linearly independent. There therefore exists a unique plane, $E\subseteq\Bbb{R}^{n+1}$, which passes through $\msf{N}_0$, $\msf{N}_1$ and the origin, and the intersection $C:=E\minter\Sigma^n$ is therefore the unique great circle passing through these two points. Finally, for all $s$, $\msf{N}(s)\in E$ and $\msf{N}(s)\in\Sigma^n$, so that $\msf{N}(s)\in E\minter\Sigma^n=C$, as desired.\qed
\medskip
Let $X$ be a subset of $\Sigma^n$ which is strictly contained in a hemisphere. In particular, no two points of $X$ are antipodal and so there is a well defined great-circular arc joining any two of them. We say that $X$ is \gloss{convex} whenever it has in addition the property that for all $\msf{N}_0,\msf{N}_1\in X$, the great-circular arc joining $\msf{N}_0$ and $\msf{N}_1$ is also contained in $K$.
\par
For all $\msf{N}\in\Sigma^n$, we define the subset $\Sigma_-^n(\msf{N})$ of $\Sigma^n$ by
$$
\Sigma_-^n(\msf{N}) := \left\{ x\ |\ \langle x,\msf{N}\rangle<0\right\}.
$$
We refer to $\Sigma_-^n(\msf{N})$ is the \gloss{open hemisphere} defined by $\msf{N}$. In particular, when $\msf{N}=e_{n+1}$, we define the \gloss{southern hemisphere} of $\Sigma^n$ by $\Sigma_-^n := \Sigma_-(e_{n+1})$. We now identify $\Bbb{R}^{n+1}$ with the product $\Bbb{R}^n\times\Bbb{R}$. Observe that $\Sigma_-^n$ then coincides with the intersection of $\Sigma$ with $\Bbb{R}^n\times]-\infty,0[$. We define the mapping $P:\Sigma^n_-\rightarrow\Bbb{R}^n$ by
$$
P(x',t) := -x'/t,
$$
and we refer to $P$ as the \gloss{affine projection} of $\Sigma^n_-$ onto $\Bbb{R}^n$.
\proclaim{Lemma \nextprocno}
\noindent $P$ defines a smooth diffeomorphism from $\Sigma_-^n$ onto $\Bbb{R}^n$.
\endproclaim
\proclabel{PIsADiffeomorphism}
\proof Define $\hat{P}:\Bbb{R}^n\times]-\infty,0[\rightarrow\Bbb{R}^n$ by $\hat{P}(x',t) := -x'/t$. Since $\hat{P}$ is smooth, and since $P$ coincides with its restriction to $\Sigma^n_-$, $P$ is also smooth. Now define $Q:\Bbb{R}^n\rightarrow\Sigma^n_-$ by $Q(x') := (x',-1)/\sqrt{1 + \|x'\|^2}$. Observe that $Q$ is smooth. Moreover, for all $(x',t)\in\Sigma^n_-$, bearing in mind that $\|x'\|^2 + t^2 = 1$,
$$
(Q\circ P)(x',t) = Q(-x'/t) = (-x'/t, -1)/\sqrt{1 + x^2/t^2} = (x',t).
$$
Conversely, for all $x'\in\Bbb{R}^n$,
$$
(P\circ Q)(x') = P((x',-1)/\sqrt{1 + \|x'\|^2}) = x'.
$$
We conclude that $P$ is a smooth diffeomorphism with inverse $Q$ as desired.\qed
\proclaim{Lemma \nextprocno}
\noindent $P$ maps the set of great-circular arcs in $\Sigma^n_-$ bijectively onto the set of straight-line segments in $\Bbb{R}^n$.
\endproclaim
\proclabel{PMapsGreatCircularArcsToLineSegments}
\proof Since straight-line segments and great-circular arcs are uniquely defined by their end points, it suffices to show that for any two distinct points $\msf{N}_0,\msf{N}_1\in\Sigma^n_-$, if $\msf{N}:[0,1]\rightarrow\Sigma^n_-$ is the great-circular arc joining these two points, then $P\circ\msf{N}$ is, up to reparametrisation, the straight-line segment joining $P(\msf{N}_0)$ to $P(\msf{N}_1)$. However, for all $t$,
$$
(P\circ\msf{N})(t)=(1-\tau(s))P(\msf{N}_0) + \tau(s)P(\msf{N}_1),
$$
where
$$
\tau(s) := \frac{s\langle\msf{N}_1,e_{n+1}\rangle}{s\langle\msf{N}_1,e_{n+1}\rangle + (1-s)\langle\msf{N}_0,e_{n+1}\rangle}.
$$
Observe that $\tau(0)=0$, $\tau(1)=1$, and
$$
\tau'(s) = \frac{\langle\msf{N}_0,e_{n+1}\rangle\langle\msf{N}_1,e_{n+1}\rangle}{s\langle\msf{N}_1,e_{n+1}\rangle + (1-s)\langle\msf{N}_0,e_{n+1}\rangle} > 0.
$$
The function $\tau$ is therefore a reparametrisation of the unit interval, and so the image of $\msf{N}$ under $P$ is a reparametrised straight-line segment from $P(\msf{N}_0)$ to $P(\msf{N}_1)$, as desired.\qed
\medskip
In technical terms, Lemma \procref{PMapsGreatCircularArcsToLineSegments} means that $\Sigma^n_-$ is affine equivalent to $\Bbb{R}^n$. In particular, this immediately yields
\proclaim{Theorem \nextprocno}
\noindent If $X$ is a subset of $\Sigma^n_-$, then $X$ is convex if and only if $P(X)$ is convex.
\endproclaim
\proclabel{PMapsConvexToConvex}
\noindent Convex subsets of $\Sigma^n_-$ therefore possess all the properties of convex subsets of $\Bbb{R}^{n+1}$ studied in Chapter \headref{TheStructureOfSingularities}. In particular, if $X$ is a subset of $\Sigma^n_-$ which is strictly contained in a hemisphere, then we define the convex hull of $X$ to be the intersection of all open convex sets in $\Sigma^n$ containing $X$. We denote this set by $\opConv(X)$.
\proclaim{Lemma \nextprocno}
\noindent Let $X$ and $K$ be subsets of $\Sigma^n$ which are strictly contained in hemispheres. Then,
\medskip
\myitem{(1)} if $X$ is compact, then $\opConv(X)$ is compact;
\medskip
\myitem{(2)} if $X$ is open, then $\opConv(X)$ is open;
\medskip
\myitem{(3)} if $K$ is compact and convex, then $\opConv(K)=K$; and
\medskip
\myitem{(4)} if $K$ is open and convex, then $\opConv(K)=K$.
\endproclaim
\proclabel{LemmaBasicPropertiesOfConvexSubsetsOfSpheres}
\proof Upon applying rotations, we may suppose that $X,K\subseteq\Sigma^n_-$. Observe that if $K'$ is an open convex subset of $\Sigma^n_-$ containing $X$, then so too is $K'\minter\Sigma^n_-$. It follows that $\opConv(X)$ coincides with the intersection of all open, convex subsets of $\Sigma^n_-$ containing $X$. Thus, since $P$ is a diffeomorphism mapping convex sets to convex sets,
$$
P(\opConv(X)) = \opConv(P(X)).
$$
$(1)$ and $(2)$ now follow by Lemmas \procref{ConvexHullOfCompactIsCompact} and \procref{ConvexHullOfOpenIsOpen}. $(3)$ follows by Lemma \procref{ConvexHullOfConvexIsTheSame}. $(4)$ is trivial, and this completes the proof.\qed
\proclaim{Lemma \nextprocno}
\noindent $P$ maps the set of open hemispheres in $\Sigma^n$ bijectively onto the set of open half-spaces in $\Bbb{R}^n$.
\endproclaim
\proclabel{PMapsOpenHemispheresToOpenHalfSpaces}
{\sl{\bf\noindent Remark:\ }We include here the empty set as the trivial open half-space.}
\medskip
\proof The operation of intersection defines a bijection between the set of open linear half-spaces in $\Bbb{R}^{n+1}$ and the set of open hemispheres in $\Sigma^n$. Now identify $\Bbb{R}^n$ with the affine hyperplane $R:=\Bbb{R}^n\times\left\{-1\right\}$ in $\Bbb{R}^{n+1}$. The operation of intersection defines a bijection between the set of open linear half-spaces in $\Bbb{R}^{n+1}$ and the set of open half-spaces in $R$. However, for every open linear half-space $H$ in $\Bbb{R}^{n+1}$,
$$
P(H\minter\Sigma^n\minter\Sigma^n_-) = H\minter R,
$$
and the result follows.\qed
\medskip
\noindent This yields an alternative characterisation of the convex hull of a subset of the sphere.
\proclaim{Theorem \nextprocno}
\noindent Let $X$ be a subset of $\Sigma^n$. If $X$ is strictly contained in a hemisphere, then $\opConv(X)$ coincides with the intersection of all open hemispheres containing $X$.
\endproclaim
\proclabel{ConvIsIntersectionOfOpenHemispheres}
\proof This follows from Theorem \procref{ConvexHullIsIntersectionOfOpenHalfSpaces} and Lemma \procref{PMapsOpenHemispheresToOpenHalfSpaces}.\qed
\newsubhead{Duality}
Let $X$ be a subset of $\Sigma^n$. We define the \gloss{dual} subset $X^*$ to $X$ by
$$
X^* := \left\{\msf{M}\in\Sigma^n\ |\ \langle\msf{N},\msf{M}\rangle<0\ \forall\ \msf{N}\in X\right\}
= \minter_{\msf{N}\in X}\Sigma_-(\msf{N}).
$$
\proclaim{Lemma \nextprocno}
\noindent $X$ is non-empty if and only if $X^*$ is strictly contained in a hemisphere.
\endproclaim
\proclabel{NonEmptyToStrictlyContainedInAHemisphere}
\proof Suppose $X$ is non-empty. Choose $\msf{N}\in X$. For all $\msf{M}\in X^*$, $\langle\msf{M},\msf{N}\rangle<0$ and so $X^*$ is strictly contained in a hemisphere, as desired. Conversely, suppose that $X$ is strictly contained in a hemisphere. Let $\msf{M}\in\Sigma^n$ be such that $\langle\msf{N},\msf{M}\rangle<0$ for all $\msf{N}\in X$. By definition, $M\in X^*$ and so $X^*$ is non-empty, as desired.\qed
\proclaim{Lemma \nextprocno}
\noindent If $X$ is closed, then $X^*$ is open. If $X$ is open, then $X^*$ is closed.
\endproclaim
\proclabel{OpenToClosed}
\proof Suppose that $X$ is closed. Choose $\msf{M}\in X^*$. By compactness of $X$, there exists $\epsilon>0$ such that $\langle\msf{N},\msf{M}\rangle\leq-\epsilon$ for all $\msf{N}$ in $X$. If $\msf{M}'\in B_\epsilon(\msf{M})\minter\Sigma^n$, then, using the Cauchy-Schwarz inequality, we obtain, for all $\msf{N}\in X$,
$$\eqalign{
\langle\msf{N},\msf{M}'\rangle
&= \langle\msf{N},\msf{M}' - \msf{M}\rangle + \langle\msf{N},\msf{M}\rangle\cr
&\leq \|\msf{N}\|\|\msf{M}'-\msf{M}\| + \langle\msf{N},\msf{M}\rangle\cr
&<\epsilon - \epsilon\cr
&=0.\cr}
$$
Since $\msf{M}'\in B_\epsilon(\msf{M})\minter\Sigma^n$ is arbitrary, we conclude that $B_\epsilon(\msf{M})\minter\Sigma^n\subseteq X^*$, and since $\msf{M}\in X^*$ is arbitrary, we conclude that $X^*$ is open, as desired.
\par
Now suppose that $X$ is open. Choose $\msf{M}\in\Sigma^n\setminus X^*$. There exists $\msf{N}\in X$ such that $\langle\msf{N},\msf{M}\rangle\geq 0$. For all $s>0$, denote $\msf{N}_s := (\msf{N} + s\msf{M})/\|\msf{N} + s\msf{M}\|$. For all $s>0$, we have
$$\eqalign{
\langle\msf{N}_s,\msf{M}\rangle
&= \frac{1}{\|\msf{N}+ s\msf{N}\|}\langle\msf{N} + s\msf{M},\msf{M}\rangle\cr
&\geq \frac{s}{\|\msf{N}+ s\msf{N}\|}\cr
&> 0.\cr}
$$
Since $X$ is open, for sufficiently small $s$, $\msf{N}_s\in X$. Thus, upon replacing $\msf{N}$ with $\msf{N}_s$, we may suppose that $\langle\msf{N},\msf{M}\rangle=:\epsilon>0$. If $\msf{M}'\in B_\epsilon(\msf{M})\minter\Sigma^n$, then, using the Cauchy-Schwarz inequality, we obtain
$$\eqalign{
\langle\msf{N},\msf{M}'\rangle&=\langle\msf{N},\msf{M}'-\msf{M}\rangle + \langle\msf{N},\msf{M}\rangle\cr
&\geq -\|\msf{N}\|\|\msf{M}' - \msf{M}\| + \langle\msf{N},\msf{M}\rangle\cr
&\geq -\epsilon + \epsilon\cr
&= 0.\cr}
$$
Since $\msf{M}'\in B_\epsilon(\msf{N})\minter\Sigma^n$ is arbitrary, we conclude that $B_\epsilon(\msf{M})\minter\Sigma^n\subseteq\Sigma^n\setminus X^*$, and since $\msf{M}\in\Sigma^n\setminus X^*$ is arbtrary, we conclude that $\Sigma^n\setminus X^*$ is open, so that $X^*$ is closed, as desired.\qed
\proclaim{Lemma \nextprocno}
\noindent If $X$ is non-empty, then $X^*$ is convex.
\endproclaim
\proclabel{NonEmptyToConvex}
\proof By Lemma \procref{NonEmptyToStrictlyContainedInAHemisphere}, $X^*$ is strictly contained in a hemisphere. Choose $\msf{M}_0,\msf{M}_1\in X^*$. In particular, $\msf{M}_0$ and $\msf{M}_1$ are not antipodal. Let $\msf{M}$ be the great-circular arc joining $\msf{M}_0$ to $\msf{M}_1$. For all $\msf{N}\in X$ and for all $s\in [0,1]$,
$$\eqalign{
\langle\msf{N},\msf{M}(s)\rangle
&= \frac{1}{\|(1-s)\msf{M}_0 + s\msf{M}_1\|}\langle \msf{N},(1-s)\msf{M}_0 + s\msf{M}_1\rangle\cr
&= \frac{(1-s)}{\|(1-s)\msf{M}_0 + s\msf{M}_1\|}\langle \msf{N},\msf{M}_0\rangle + \frac{s}{\|(1-s)\msf{M}_0 + s\msf{M}_1\|}\langle \msf{N},\msf{M}_1\rangle\cr
&< 0,\cr}
$$
so that $\msf{M}(s)$ is an element of $X^*$ for all $s\in[0,1]$. Since $\msf{M}_0,\msf{M}_1\in X^*$ are arbitrary, we conclude that $X^*$ is convex, as desired.\qed
\proclaim{Lemma \nextprocno}
\noindent If $X$ is strictly contained in a hemisphere, then $X^{**}=\opConv(X)$.
\endproclaim
\proclabel{DoubleDualOfClosedIsConvexHull}
\proof By Lemma \procref{ConvIsIntersectionOfOpenHemispheres}, $\opConv(X)$ is the intersection of all open hemispheres containing $X$. However, by definition, $\msf{M}$ is an element of $X^*$ if and only if $X$ is contained in $\Sigma_-(\msf{M})$. That is, $X^*$ parametrises the set of open hemispheres containing $X$. Thus,
$$
X^{**} = \minter_{\msf{M}\in X^*}\Sigma_-(\msf{M}) = \opConv(X),
$$
as desired.\qed
\proclaim{Lemma \nextprocno}
\noindent If $X_1$ and $X_2$ are both strictly contained in the same hemisphere, then $(X_1\munion X_2)^* = X_1^*\minter X_2^*$.
\endproclaim
\proclabel{DualOfUnionIsIntersectionOfDual}
\proof Suppose $\msf{M}\in (X_1\munion X_2)^*$. Then $\langle\msf{N},\msf{M}\rangle<0$ for all $\msf{N}\in X_1\munion X_2$ and so $\msf{M}\in X_1^*\minter X_2^*$. Conversely, if $\msf{M}\in X_1^*\minter X_2^*$, then $\langle\msf{N},\msf{M}\rangle<0$ for all $\msf{N}\in X_1\munion X_2$ and so $\msf{M}\in (X_1\munion X_2)^*$. These two sets therefore coincide, as desired.\qed
\proclaim{Lemma \nextprocno}
\noindent Let $X_1$ and $X_2$ be convex subsets of $\Sigma^n$ which are both strictly contained in a hemisphere. If $X_1$ and $X_2$ are either both open or both closed, and if $X_1\minter X_2$ is non-empty, then $(X_1\minter X_2)^* = \opConv(X_1^*\munion X_2^*)$.
\endproclaim
\proclabel{DualOfIntersectionIsConvexHullOfUnionOfDual}
\proof Suppose that $X_1$ and $X_2$ are open (resp. closed). Since they are both convex, by Lemmas \procref{LemmaBasicPropertiesOfConvexSubsetsOfSpheres} and \procref{DoubleDualOfClosedIsConvexHull}, $X_1=\opConv(X_1)=X_1^{**}$ and $X_2=\opConv(X_2)=X_2^{**}$.  We denote $Y_1=X_1^*$ and $Y_2=X_2^*$. Observe that, if $\msf{N}\in X_1\minter X_2$, then $Y_1=X_1^*\in\Sigma_-(\msf{N})$ and $Y_2=X_2^*\subseteq\Sigma_-(\msf{N})$. That is, since $X_1\minter X_2$ is non-empty, $Y_1$ and $Y_2$ are both strictly contained in the same hemisphere. Thus, using Lemma \procref{DualOfUnionIsIntersectionOfDual}, we obtain
$$
(X_1^*\munion X_2^*)^* = (Y_1\munion Y_2)^* = Y_1^*\minter Y_2^* = X_1\minter X_2,
$$
so that, by Lemma \procref{DoubleDualOfClosedIsConvexHull},
$$
\opConv(X_1^*\munion X_2^*) = (X_1^*\munion X_2^*)^{**} = (X_1\minter X_2)^*,
$$
as desired.\qed
\newsubhead{Links}
Let $K$ be a compact, convex set with non-trivial interior. Let $x$ be a boundary point of $K$. For $r>0$, we define $\Cal{L}_r(x;K)\subseteq\Sigma^n$, the \gloss{link} of $K$ of radius $r$ about $x$ by,
\subheadlabel{Links}
$$
\Cal{L}_r(x;K) := \left\{\msf{N}\ |\ x + r\msf{N}\in K^o\right\}.
$$
When there is no ambiguity, we denote $\Cal{L}_r(x)=\Cal{L}_r(x;K)$.
\proclaim{Lemma \nextprocno}
\noindent For every boundary point $x$ of $K$ and for all $r<s$, $\Cal{L}_s(x)\subseteq\Cal{L}_r(x)$.
\endproclaim
\proof Indeed, choose $\msf{N}\in\Cal{L}_s(x)$. Then $x + s\msf{N}\in K^o$. Viewing $\msf{N}$ as an element of $\Bbb{R}^{n+1}$, there exists $\delta>0$ such that for all $V\in B_\delta(0)$, $x+s\msf{N} + V\in K$. Thus, by convexity, for all $V\in B_{r\delta/s}(0)$,
$$
x + r\msf{N} + V = (1-\frac{r}{s})x + \frac{r}{s}(x + s\msf{N} + \frac{s}{r}V) \in K.
$$
It follows that $x + r\msf{N}\in K^o$, and so $\msf{N}\in \Cal{L}_r(x)$. Since $\msf{N}\in\Cal{L}_s(x)$ is arbitrary, we conclude that $\Cal{L}_s(x)\subseteq\Cal{L}_r(x)$ as desired.\qed
\medskip
\noindent $(\Cal{L}_r(x))_{r>0}$ therefore constitutes an increasing, nested family of open sets. We define $\Cal{L}(x;K)\subseteq\Sigma^{n}$, the \gloss{link} of $K$ at $x$ by,
$$
\Cal{L}(x;K) := \munion_{r>0}\Cal{L}_r(x;K).
$$
When there is no ambiguity, we denote $\Cal{L}(x)=\Cal{L}(x;K)$. Since it is the union of a family of open sets, $\Cal{L}(x)$ is also open.
\proclaim{Lemma \nextprocno}
\noindent Let $K$ be a compact, convex set with non-trivial interior. Then for every boundary point $x$ of $K$, $\Cal{N}(x;K)=\Cal{L}(x;K)^*$.
\endproclaim
\proclabel{SupportingNormalsIsDualOfLink}
\proof Suppose that $\msf{N}\in\Cal{N}(x;K)$. For all $z\in K$, $\langle z - x,\msf{N}\rangle\leq 0$, and so, for all $z\in K^o$, $\langle z - x,\msf{N}\rangle<0$. Choose $\msf{M}\in\Cal{L}(x;K)$. Choose $r>0$ such that $\msf{M}\in\Cal{L}_r(x;K)$. Then $x + r\msf{M}\in K^o$, and so
$$
\langle\msf{M},\msf{N}\rangle = \frac{1}{r}\langle (x + r\msf{M}) - x,\msf{N}\rangle < 0.
$$
Since $\msf{M}\in\Cal{L}(x;K)$ is arbitrary, we conclude that $\msf{N}\in\Cal{L}(x;K)^*$, and since $\msf{N}\in\Cal{N}(x;K)$ is arbitrary, we conclude that $\Cal{N}(x;K)\subseteq\Cal{L}(x;K)^*$.
\par
Conversely, choose $\msf{N}\in\Cal{L}(x;K)^*$. Since $K$ has non-trivial interior, by Lemma \procref{ClosureOfInteriorIsSame}, $K=\overline{K^o}$. Choose $y\in K^o$. Denote $r:=\|y-x\|$. Then $(y-x)/r\in\Cal{L}_r(x;K)$, and so
$$
\langle\msf{N},y-x\rangle = r\langle\msf{N},(y-x)/r\rangle < 0,
$$
Thus, by continuity, for all $y\in\overline{K^o}=K$,
$$
\langle\msf{N},y - x\rangle \leq 0,
$$
so that $\msf{N}\in\Cal{N}(x;K)$. Since $\msf{N}\in\Cal{L}(x;K)^*$ is arbitrary, we conclude that $\Cal{L}(x;K)^*\subseteq\Cal{N}(x;K)$, and the two sets therefore coincide, as desired.\qed
\proclaim{Lemma \nextprocno}
\noindent For every compact, convex set $K$ with non-trivial interior, and for every boundary point $x$ of $K$, $\Cal{N}(x;K)$ is closed, convex and strictly contained in a hemisphere.
\endproclaim
\proof By Lemma \procref{SupportingNormalsIsDualOfLink}, $\Cal{N}(x;K) = \Cal{L}(x;K)^*$. Since $K$ has non-trivial interior, $\Cal{L}(x;K)$ is non-empty, and so, by Lemma \procref{NonEmptyToStrictlyContainedInAHemisphere}, $\Cal{N}(x;K)$ is strictly contained in a hemisphere. Since $\Cal{L}(x;K)$ is open, by Lemma \procref{LemmaBasicPropertiesOfConvexSubsetsOfSpheres}, $\Cal{N}(x;K)$ is closed. Finally, by Lemma \procref{NonEmptyToConvex}, $\Cal{N}(x;K)$ is convex, and this completes the proof.\qed
\proclaim{Theorem \nextprocno}
\noindent Let $K_1$ and $K_2$ be two compact, convex sets whose intersection has non-trivial interior. Choose $x\in\partial(K_1\minter K_2)$. Then,
\medskip
\myitem{(1)} if $x\in(\partial K_1)\minter K_2^o$, then $\Cal{N}(x;K_1\minter K_2) = \Cal{N}(x;K_1)$;
\medskip
\myitem{(2)} if $x\in K_1^o\minter (\partial K_2)$, then $\Cal{N}(x;K_1\minter K_2) = \Cal{N}(x;K_2)$; and
\medskip
\myitem{(3)} if $x\in (\partial K_1)\minter (\partial K_2)$, then $\Cal{N}(x;K_1\minter K_2) = \opConv(\Cal{N}(x;K_1)\munion\Cal{N}(x;K_2))$.
\endproclaim
\proclabel{SupportingNormalsOfIntersection}
\proof Cases $(1)$ and $(2)$ follow from Lemma \procref{SupportingNormalIsLocalProperty}. However, using Lemmas \procref{DualOfIntersectionIsConvexHullOfUnionOfDual} and \procref{SupportingNormalsIsDualOfLink}, for $x\in (\partial K_1)\minter(\partial K_2)$, we obtain
$$\eqalign{
\Cal{N}(x;K_1\minter K_2) &= \Cal{L}(x;K_1\minter K_2)^* \cr
&= (\Cal{L}(x;K_1)\minter\Cal{L}(x;K_2))^* \cr
&= \opConv(\Cal{L}(x;K_1)^*\munion\Cal{L}(x;K_2)^*)\cr
&= \opConv(\Cal{N}(x;K_1)\munion\Cal{N}(x;K_2)),\cr}$$
as desired.\qed

\newhead{Weak Barriers}
For $k$ a positive real number, the set of weak barriers of gaussian curvature at least $k$ is essentially the closure in the Hausdorff topology of the set of compact, convex sets with smooth boundary of gaussian curvature at least $k$. This concept allows us to solve the Plateau problem in euclidean space for very general data. Once solutions have been found, the theory developed in Section \headref{TheStructureOfSingularities} is then applied to identify their singular sets. In particular, under suitable conditions on the boundary, these are shown to be empty, so that the solutions are actually smooth.
\par
Existence is proven via the Perron method. The main requirement for the application of this technique is the closure of the family of weak barriers under the operation of intersection. That is, if $K_1$ and $K_2$ are weak barriers, then so too is $K_1\minter K_2$. The proof of this result, which is encapsulated in Theorem \procref{SmoothIntersectionHasHighCurvature} is rather technical, and forms the content of Sections \subheadref{DistanceFunctions} to \subheadref{SmoothingTheIntersection} inclusive. The techniques used are mostly elementary, although some experience of the theory of distributions will be required, and we refer the reader to \cite{FriedlanderJoshi} for a clear and straightforward introduction.
\par
The experienced reader will notice that weak barriers are always viscosity supersolutions of the Gauss curvature equation (c.f. \cite{CrandhallIshiiLyons}). Like the space of weak barriers, the space of viscosity supersolutions is closed with respect to the Hausdorff topology and is closed under finite intersections. Furthermore, in contrast to weak barriers, these properties for viscosity supersolutions are almost trivial. However, viscosity supersolutions, on the other hand, do not obviously possess the properties required for us to apply Theorem \procref{FirstExistenceTheorem} as the regularising operation in the application of the Perron method. It is precisely for this reason that the more technical notion of weak barriers is required.
\par
Finally, the results of Sections \subheadref{DistanceFunctions} to \subheadref{SmoothingTheIntersection} inclusive are very general, and are useful for constructing convex barriers in a wide range of settings. In particular, we leave the enthusiastic reader to verify that they remain valid in any riemannian manifold.
\newsubhead{Distance functions}
Let $K$ be a closed, convex subset of $\Bbb{R}^{n+1}$. Let $d_K:\Bbb{R}^{n+1}\rightarrow[0,\infty[$ be the distance in $\Bbb{R}^{n+1}$ to $K$. That is,
\subheadlabel{DistanceFunctions}
\headlabel{TheTheoryOfWeakSolutions}
$$
d_K(x) := \inf_{y\in K}\|x - y\|.
$$
Since it is the infimum of a family of convex functions, $d_K$ is also convex. We now consider the closest point projection from $\Bbb{R}^{n+1}\setminus K$ onto $K$. First, we prove
\proclaim{Lemma \nextprocno}
\noindent Let $K$ be a closed, convex subset of $\Bbb{R}^{n+1}$. Choose $x\in\Bbb{R}^{n+1}\setminus K$. There is at most one point $y$ in the boundary of $K$ with the property that $x = y + t\msf{N}$ for some $t>0$ and for some supporting normal $\msf{N}$ to $K$ at $y$.
\endproclaim
\proclabel{PointMappingToXIsUnique}
\proof Suppose the contrary. Let $y$ and $y'$ be two such boundary points. Let $\msf{N}$ and $\msf{N}'$ be supporting normals to $K$ at $y$ and $y'$ respectively and let $t,t'>0$ be such that $x=y+t\msf{N}=y'+t'\msf{N}'$. By definition of the supporting normal,
$$
\langle y' - y,\msf{N}\rangle,\ \langle y - y',\msf{N}'\rangle \leq 0.
$$
In particular,
$$
\langle y' - y,x - y\rangle,\ \langle y - y',x - y'\rangle\leq 0.
$$
Summing these two relations yields $\|y'-y\|^2\leq 0$, so that $\|y'-y\|=0$, and so $y'=y$, as desired.\qed
\proclaim{Lemma \nextprocno}
\noindent Let $K$ be a closed, convex subset of $\Bbb{R}^{n+1}$. For all $x\in\Bbb{R}^{n+1}$, the point $y\in K$ minimising distance to $x$ is unique.
\endproclaim
\proclabel{ClosestPointIsUnique}
\proof Choose $x\in\Bbb{R}^{n+1}$. Let $y\in K$ minimise distance to $x$. If $x\in K$, then $y=x$ is unique, as desired. Otherwise, denote $\msf{N}:=(x-y)/\|x-y\|$. By Lemma \procref{CharacterisationOfClosestPoints}, $y$ is a boundary point of $K$ and $\msf{N}$ is a supporting normal to $K$ at $y$. In particular $x = y + \|x-y\|\msf{N}$, and by Lemma \procref{PointMappingToXIsUnique}, $y$ is unique, as desired. This completes the proof.\qed
\medskip
\noindent We define $\Pi_K:\Bbb{R}^{n+1}\rightarrow K$ to be the closest point projection. We now relate $\Pi_K$ to the derivative of $d_K$.
\proclaim{Lemma \nextprocno}
\noindent If $K$ as a closed, convex subset of $\Bbb{R}^{n+1}$, then $d_K$ is differentiable at every point $x$ in $\Bbb{R}^{n+1}\setminus K$ and, for all such $x$,
$$
\Pi_K(x) = x - d_K(x)Dd_K(x).
$$
\endproclaim
\proclabel{FirstDerivativeOfDistance}
\proof Choose $x\in\Bbb{R}^{n+1}\setminus K$. Denote $\msf{N} = (x-\Pi_K(x))/\|x - \Pi_K(x)\|$. By Lemma \procref{CharacterisationOfClosestPoints}, $\msf{N}$ is a supporting normal to $K$ at $\Pi_K(x)$. Since the image of $\Pi_K$ is contained in $K$, it follows that, for all $y\in\Bbb{R}^{n+1}$, $\langle\Pi_K(y) - \Pi_K(x),\msf{N}\rangle\leq 0$. Using the Cauchy-Schwarz inequality and the fact that $\msf{N}$ has unit length, we therefore obtain, for all $y$,
$$\eqalign{
d_K(y) &= \|y - \Pi_K(y)\|\cr
&\geq \langle y - \Pi_K(y),\msf{N}\rangle\cr
&= \langle y - \Pi_K(x),\msf{N}\rangle + \langle \Pi_K(x) - \Pi_K(y),\msf{N}\rangle\cr
&\geq \langle y - \Pi_K(x),\msf{N}\rangle.\cr}
$$
On the other hand, $d_K(y)\leq d(y,\Pi_K(x))$, so that
$$
\langle y - \Pi_K(x),\msf{N}\rangle \leq d_K(y) \leq d(y,\Pi_K(x)).
$$
The first and the last functions in this inequality are smooth at $x$. Moreover, they coincide at $x$ and have derivative equal to $\msf{N}$ at this point. It follows that $d_K$ is differentiable at $x$ and $Dd_K(x)=\msf{N}$. In particular,
$$
\Pi_K(x) = x - \|x - \Pi_K(x)\|\msf{N} = x - d_K(x)Dd_K(x),
$$
as desired.\qed
\proclaim{Lemma \nextprocno}
\noindent Let $K$ be a closed, convex subset of $\Bbb{R}^{n+1}$. If $x$ is a point of $\Bbb{R}^{n+1}\setminus K$, then $d_K$ is twice differentiable at $x$ if and only if $\Pi_K$ is differentiable at $x$. Moreover, at any such point, for all vectors $V$ and $W$,
$$
\langle D\Pi_K(x)V,W\rangle = \langle\pi(V),\pi(W)\rangle - d_K(x)D^2d_K(x)(V,W),
$$
where $\pi$ is the orthogonal projection from $\Bbb{R}^{n+1}$ onto $\langle Dd_K(x)\rangle^\perp$.
\endproclaim
\proclabel{TwiceDifferentiableWheneverPiIsDifferentiable}
\proof By Lemma \procref{FirstDerivativeOfDistance}, for all $x\in\Bbb{R}^{n+1}\setminus K$, $d_K$ is differentiable at $x$ and $\Pi_K(x) = x - d_K(x)Dd_K(x)$. Since $d_K(x)>0$, it follows by the product and quotient rules that $Dd_K$ is differentiable at $x$ if and only if $\Pi_K$ is. Furthermore, at any such point
$$\eqalign{
\langle D\Pi_K(x)V,W\rangle &= \langle V,W\rangle - \langle V,Dd_K(x)\rangle\langle W,Dd_K(x)\rangle - d_K(x)D^2d_K(x)(V,W)\cr
&= \langle\pi(V),\pi(W)\rangle - d_K(x)D^2d_K(x)(V,W),\cr
}$$
as desired.\qed
\medskip
We now consider the regularity of $\Pi_K$.
\proclaim{Lemma \nextprocno}
\noindent If $K$ is a closed, convex subset of $\Bbb{R}^{n+1}$, then $\Pi_K$ is $1$-Lipschitz.
\endproclaim
\proclabel{ProjectionIsLipschitz}
\proof If $x,x'\in K$, then $\Pi_K(x)=x$ and $\Pi_K(x')=x'$. In particular, $\|\Pi_K(x) - \Pi_K(x')\|=\|x-x'\|$, as desired. If $x\in K$ and if $x'\in\Bbb{R}^{n+1}\setminus K$, then $\pi_K(x)=x$. Define $y':=\Pi_K(x')$. By Lemma \procref{CharacterisationOfClosestPoints}, $(x'-y')/\|x'-y'\|$ is a supporting normal to $K$ at $y'$. In particular, $\langle x - y',x'-y'\rangle\leq 0$, and so
$$\eqalign{
\| x -x'\|^2 &= \|(x-y') - (x'-y')\|^2\cr
&=\|x-y'\|^2 -2\langle x-y',x'-y'\rangle + \|x'-y'\|^2\cr
&\geq \|x-y'\|^2 + \|x'-y'\|^2\cr
&\geq \|x-y'\|^2,\cr
}$$
so that $\|\Pi_K(x) - \Pi_K(x')\|\leq\|x - x'\|$, as desired. Finally, choose $x,x'\in\Bbb{R}^{n+1}\setminus K$. Denote $y:=\Pi_K(x)$, $y':=\Pi_K(x')$. By Lemma \procref{CharacterisationOfClosestPoints}, $(x-y)/\|x-y\|$ and $(x'-y')/\|x'-y'\|$ are supporting normals to $K$ at $y$ and $y'$ respectively. In particular,
$$
\langle y' - y,x - y\rangle, \langle y - y',x' - y'\rangle \leq 0.
$$
Consequently,
$$
\langle x - x', y - y'\rangle = \langle x - y, y - y'\rangle + \langle y - y',y-y'\rangle + \langle y' - x',y - y'\rangle\geq \| y - y'\|^2.
$$
Using the Cauchy-Schwarz inequality, this yields
$$
\|y-y'\|^2 \hfill\leq \langle x - x',y-y'\rangle\leq \|x-x'\|\|y - y'\|,
$$
so that,
$$
\|y - y'\|(\|x-x'\| - \|y - y'\|) \geq 0,
$$
and we conclude that $\|y-y'\|\leq \|x-x'\|$, as desired.\qed
\proclaim{Lemma \nextprocno}
\noindent If $K$ is a closed, convex subset of $\Bbb{R}^{n+1}$, then $\Pi_K$ is differentiable almost everywhere. Moreover, the pointwise derivative of $\Pi_K$ coincides with its distributional derivative, and $\|D\Pi_K(x)\|_{L^\infty}\leq 1$.
\endproclaim
\proclabel{PiIsDifferentiable}
\proof Since $\Pi_K$ is Lipschitz, it follows from Rademacher's Theorem (c.f. Theorem $5.2$ of \cite{Simon}) that $\Pi_K$ is differentiable almost everywhere and, moreover, that its pointwise derivative coincides with its distributional derivative. Furthermore, since $\Pi_K$ is $1$-Lipschitz, it follows that $\|D\Pi_K\|_{L^\infty}\leq 1$, and this completes the proof.\qed
\proclaim{Lemma \nextprocno}
\noindent If $K$ is a closed, convex subset of $\Bbb{R}^{n+1}$, then $d_K$ is twice differentiable almost everywhere in $\Bbb{R}^{n+1}\setminus K$. Moreover, the pointwise second derivative of $d_K$ coincides with its second-order distributional derivative, and if $d_K$ is twice differentiable at $x\in\Bbb{R}^{n+1}\setminus K$, then $\|D^2d_K(x)\|\leq 2/d_K(x)$.
\endproclaim
\proclabel{DIsTwiceDifferentiable}
\proof By Lemma \procref{TwiceDifferentiableWheneverPiIsDifferentiable}, $d_K$ is twice differentiable wherever $\Pi_K$ is differentiable, and so, by Lemma \procref{PiIsDifferentiable}, $d_K$ is twice differentiable almost everywhere. By Lemma \procref{FirstDerivativeOfDistance}, for all $x\in\Bbb{R}^{n+1}\setminus K$,
$$
Dd_K(x) = (x - \Pi_K(x))/d_K(x).
$$
Since $d_K(x)$ never vanishes over this set, using the quotient rules for pointwise derivatives and for distributional derivatives, it follows from Lemma \procref{PiIsDifferentiable} again that the pointwise second-order derivative of $d_K$ coincides with its second-order distributional derivative. Furthermore, at any point $x$ where $d_K$ is twice differentiable, for all vectors $V$ and $W$,
$$
D^2d_K(V,W) = \frac{1}{d_K(x)}(\langle\pi(V),\pi(W)\rangle - \langle D\Pi_K(x)V,W\rangle),
$$
where $\pi$ is the orthogonal projection from $\Bbb{R}^{n+1}$ onto $\langle Dd_K(x)\rangle^\perp$. In particular, since both $\pi$ and $D\Pi_K(x)$ have norm $1$,
$$
\left|D^2d_K(x)(V,W)\right| \leq \frac{2}{d_K(x)}\|V\|\|W\|,
$$
so that $\|D^2d_K(x)\|\leq 2/d_K(x)$, as desired.\qed
\medskip
We also show that the second derivatives of $d_K$ are almost everywhere symmetric.
\proclaim{Lemma \nextprocno}
\noindent For almost all $x\in\Bbb{R}^{n+1}\setminus K$, $d_K$ is twice differentiable at $x$ and its second derivative is symmetric at that point.
\endproclaim
\proclabel{SecondDerivativeOfDIsSymmetric}
\proof By Lemma \procref{DIsTwiceDifferentiable}, $d_K$ has $L^\infty_\oploc$ second-order, distributional derivatives over $\Bbb{R}^{n+1}\setminus K$. Denote this second-order distributional derivative by $A$. Then, for any $\phi\in C^\infty_\oploc(\Bbb{R}^{n+1}\setminus K)$, and for all $1\leq i,j\leq n$,
$$\eqalign{
\int_{\Bbb{R}^{n+1}\setminus K} A(x)(\partial_i,\partial_j)\phi(x)\opdVol_x &= \int_{\Bbb{R}^{n+1}\setminus K} d_K(x)D^2\phi(x)(\partial_j,\partial_i)\opdVol_x\cr
&= \int_{\Bbb{R}^{n+1}\setminus K}d_K(x)D^2\phi(x)(\partial_i,\partial_j)\opdVol_x\cr
&= \int_{\Bbb{R}^{n+1}\setminus K}A(x)(\partial_j,\partial_i)\phi(x)\opdVol_x.\cr}
$$
Since $\phi\in C^\infty_\oploc(\Bbb{R}^{n+1}\setminus K)$ is arbitrary, we conclude that $A(x)(\partial_i,\partial_j)=A(x)(\partial_j,\partial_i)$ for almost all $x\in\Bbb{R}^{n+1}\setminus K$, and since $1\leq i,j\leq n$ are arbitrary, we conclude that $A(x)$ is symmetric for almost all $x\in\Bbb{R}^{n+1}\setminus K$. However, by Lemma \procref{DIsTwiceDifferentiable}, again, for almost all $x\in\Bbb{R}^{n+1}\setminus K$, $d_K$ is twice differentiable at $x$ in the classical sense and $D^2d_K(x)=A(x)$, so that $D^2d_K(x)$ is almost everywhere defined and symmetric, as desired.\qed
\newsubhead{Convex sets with smooth boundary}
Let $K$ be a closed, convex subset of $\Bbb{R}^{n+1}$. Let $U$ be an open subset of $\Bbb{R}^{n+1}$. We denote $U(K) = U\minter (\partial K)$, and we suppose that $U(K)$ is smooth. We now use the terminology of riemannian geometry (c.f. \cite{doCarmoII}). Let $\msf{N}:U(K)\rightarrow\Sigma^n$ be the outward-pointing, unit, \gloss{normal} vector field over $U(K)$. Let $A$ be the \gloss{shape operator} of $U(K)$ associated to this normal. That is, for all $x\in U(K)$ and for any vector $V$ tangent to $U(K)$ at $x$, $A(x)V = D\msf{N}(x)V$.
\par
If $M\in\opSymm(2,\Bbb{R}^{n+1})$ is a symmetric matrix over $\Bbb{R}^{n+1}$, and if $E$ is any subspace of $\Bbb{R}^{n+1}$, we denote by $\opDet(M;E)$ the determinant of the restriction of $M$ to $E$. We are interested in estimating $\opDet(D^2d_K;\langle Dd_K\rangle^\perp)$ near $U(K)$. This quantity will be used in the sequel to estimate the gaussian curvature of smooth hypersurfaces approximating $K$.
\par
In this section, we study the functions $d_K$ and $\Pi_K$ over the set $\Pi_K^{-1}(U(K))$. We define $\Phi:U(K)\times[0,\infty[\rightarrow\Bbb{R}^{n+1}$ by $\Phi(x,t) = x + t\msf{N}(x)$.
\proclaim{Lemma \nextprocno}
\noindent $\Phi$ defines a smooth diffeomorphism from $U(K)\times[0,\infty[$ onto $\Pi_K^{-1}(U(K))$.
\endproclaim
\proof We first show that $\opIm(\Phi)=\Pi_K^{-1}(U(K))$. Indeed, choose $(x,t)\in U(K)\times[0,\infty[$. Then $\Phi(x,t) = x + t\msf{N}(x)$. By Lemma \procref{CharacterisationOfClosestPoints}, $x$ minimises distance to $\Phi(x,t)$ in $K$ so that $x=(\Pi_K\circ\Phi)(x,t)$, and, in particular, $\Phi(x,t)\in\Pi_K^{-1}(U(K))$. Since $(x,t)\in U(K)\times[0,\infty[$ is arbitrary, we conclude that $\opIm(\Phi)\subseteq\Pi_K^{-1}(U(K))$. Conversely, choose $y\in\Pi_K^{-1}(U(K))$. Denote $x=\Pi_K(y)\in U(K)$. By definition, $x$ minimises distance in $K$ to $y$. There are two cases to consider. First, if $y\in K$, then $y=x=\Phi(x,0)$, so that $y\in\opIm(\Phi)$. Second, if $y\in\Bbb{R}^{n+1}\setminus K$, then, by Lemma \procref{CharacterisationOfClosestPoints}, there exists $t>0$ such that $y = x+ t\msf{N}(x) = \Phi(x,t)$, so that $y\in\opIm(\Phi)$ in this case also. Since $y\in\Pi_K^{-1}(U(K))$ is arbitrary, we conclude that $\Pi_K^{-1}(U(K))\subseteq\opIm(\Phi)$, and the two sets therefore coincide, as desired.
\par
If $x,x'\in U(K)$ and $t,t'\in[0,\infty[$ are such that $x + t\msf{N}(x) = x' + t'\msf{N}(x')$, then, by Lemma \procref{PointMappingToXIsUnique}, $x=x'$ and $t=t'$, and it follows that $\Phi$ is injective. It remains to show that $\Phi$ is smooth with smooth inverse. Choose $(x,t)\in U(K)\times[0,\infty[$. Denote by $\partial_t$ the unit vector in the $t$ direction. Observe that
$$\triplealign{
&D\Phi(x,t)(0,\partial_t)&= \msf{N}(x)\cr
\Rightarrow&\|D\Phi(x,t)(0,\partial_t)\|^2&=1.\cr}
$$
Let $X$ be a tangent vector to $\partial K$ at $x$. Then,
$$\triplealign{
&D\Phi(x)(X,0) &= X + t D\msf{N}(x)X\cr
& &= X + tA(x)X\cr
\Rightarrow&\|D\Phi(x)(X,0)\|^2&=\|(\opId + tA(x))(X)\|^2.\cr}
$$
However, by convexity, $A(x)$ is non-negative definite, and so,
$$
\|D\Phi(x)(X,0)\|^2 = \|(\opId + tA(x))(X)\|^2 \geq \|X\|^2.
$$
Finally, bearing in mind that $\langle A(x)X,\msf{N}(x)\rangle = 0$,
$$
\langle D\Phi(x,t)(X,0),D\Phi(x,t)(0,\partial_t)\rangle = \langle X + tA(x)X,\msf{N}(x)\rangle = 0.
$$
It follows that $\|D\Phi(x,t)(V)\|^2>0$ for all non-zero $V$ and so $D\Phi(x,t)$ is invertible. Since $(x,t)\in U\times[0,\infty[$ is arbitrary, we conclude from the inverse function theorem that $\Phi$ is everywhere a smooth local diffeomorphism. By injectivity, it is a smooth global diffeomorphism, and this completes the proof.\qed
\proclaim{Lemma \nextprocno}
\noindent $\Pi_K$ and $d_K$ define smooth functions over $\Pi_K^{-1}(U(K))\setminus K$. Moreover, for all vectors $V$ and $W$,
$$\eqalign{
Dd_K(x)(V) &= \langle\msf{N}(\Pi_K(x)),V\rangle,\cr
D^2d_K(x)(V,W) &= \langle A(\Pi_K(x))D\Pi_K(x)V,W\rangle.\cr}
$$
\endproclaim
\proclabel{DistanceAndProjectionAreSmooth}
\proof Choose $(x,t)\in U(K)\times]0,\infty[$. Since $\Phi(x,t) = x + t\msf{N}(x)$, by Lemma \procref{CharacterisationOfClosestPoints}, $x$ minimises distance in $K$ to $\Phi(x,t)$. It follows that $(d_K\circ\Phi)(x,t)=t$ and $(\Pi_K\circ\Phi)(x,t)=x$. In particular, $\Pi_K\circ\Phi$ and $d_K\circ\Phi$ are both smooth, and, composing with $\Phi^{-1}$, we conclude that $d_K$ and $\Pi_K$ are also both smooth, as desired. Now choose $x\in\Pi^{-1}_K(U(K))\setminus K$. Observe that $\msf{N}(\Pi_K(x))$ is the unique supporting normal to $K$ at $\Pi_K(x)$. Thus, by Lemma \procref{FirstDerivativeOfDistance},
$$
Dd_K(x) = \frac{1}{d_K(x)}\left(x - \Pi_K(x)\right) = \msf{N}(\Pi_K(x)).
$$
The formula for the second derivative of $d_K$ follows by differentiating this relation, and this completes the proof.\qed
\proclaim{Lemma \nextprocno}
\noindent For every compact subset $X$ of $U$, there exists $C>0$ such that for all $x\in\Pi^{-1}_K(X\minter U(K))\setminus K$,
$$
\|D\Pi_K(x) - \pi\| \leq Cd_K(x),
$$
where $\pi$ is the orthogonal projection onto $\langle Dd_K(x)\rangle^\perp$.
\endproclaim
\proclabel{PiIsCloseToOrthogonalProjection}
\proof Let $C$ be such that $\|A(y)\|\leq C$ for all $y\in X\minter U(K)$. By Lemma \procref{PiIsDifferentiable}, $\|D\Pi_K(x)\|\leq 1$. Thus, by Lemma \procref{DistanceAndProjectionAreSmooth}, for all vectors $V$ and $W$,
$$\eqalign{
\left|D^2d_K(x)(V,W)\right| &= \left|\langle A(\Pi_K(x))D\Pi_K(x)V,W\rangle\right| \cr
&\leq \|A(\Pi_K(x))\|\|V\|\|W\| \cr
&\leq C\|V\|\|W\|.\cr}
$$
However, by Lemma \procref{TwiceDifferentiableWheneverPiIsDifferentiable}, for all vectors $V$ and $W$,
$$
\langle D\Pi_K(x)V,W\rangle  = \langle\pi(V),\pi(W)\rangle - d_K(x)D^2d_K(x)(V,W),
$$
where $\pi$ is the orthogonal projection from $\Bbb{R}^{n+1}$ onto $\langle Dd_K(x)\rangle^\perp$. Since $\langle\pi(V),\pi(W)\rangle=\langle\pi(V),W\rangle$, it follows that,
$$
\left|\langle D\Pi_K(x)V - \pi(V),W\rangle\right| \leq d_K(x)\left|D^2d_K(x)(V,W)\right| \leq Cd_K(x)\|V\|\|W\|,
$$
so that $\|D\Pi_K(x) - \pi\|\leq Cd_K(x)$, as desired.\qed
\proclaim{Lemma \nextprocno}
\noindent Choose $k>0$ and suppose that $U(K)$ has gaussian curvature everywhere at least $k$. For every compact subset $X$ of $U$ and for all $\epsilon>0$, there exists $r>0$ such that for all $x\in\Pi_K^{-1}(X\minter U(K))\setminus K$, if $d_K(x)<r$, then $\opDet(D^2d_K(x);\langle Dd_K(x)\rangle^\perp)\geq(k-\epsilon)^n$.
\endproclaim
\proclabel{LowerBoundsOnOrthogonalComponentOfDeterminantNearHypersurface}
\proof By compactness, there exists $\delta>0$ such that for all $x\in X\minter U(K)$ and for all $M\in B_\delta(A(x))$, $\opDet(M;\langle Dd_K(x)\rangle^\perp)\geq (k-\epsilon)^n$. Let $C_1$ be such that for all $y\in X\minter U(K)$, $\|A(y)\|\leq C_1$. Let $C_2$ be as in Lemma \procref{PiIsCloseToOrthogonalProjection}. If $x\in \Pi_K^{-1}(X\minter U(K))$ is such that $d_K(x)<\delta/C_1C_2$, then, for all vectors $V$ and $W$ in $\langle Dd_K(x)\rangle^\perp$,
$$\eqalign{
\left|D^2d_K(x)(V,W) - \langle A(\Pi_K(x))V,W\rangle\right|
&=\left|\langle A(\Pi_K(x))(D\Pi_K(x)(V) - V),W\rangle\right|\cr
&<\delta\|V\|\|W\|,\cr}
$$
so that $\opDet(D^2d_K(x);\langle Dd_K(x)\rangle^\perp)\geq (k-\epsilon)^n$, as desired.\qed
\newsubhead{Intersecting convex sets}
Let $K_1$ and $K_2$ be compact, convex subsets of $\Bbb{R}^{n+1}$ whose intersection has non-trivial interior. Let $U$ be an open subset of $\Bbb{R}^{n+1}$ and suppose that $U(K_1)$ and $U(K_2)$ are both smooth of gaussian curvature at least $k$. We denote $K:=K_1\minter K_2$, we denote by $\msf{N}_1$ and $\msf{N}_2$ the outward-pointing, unit, normal vector fields over $K_1(U)$ and $K_2(U)$ respectively and we denote by $A_1$ and $A_2$ their respective shape operators. Moreover, we denote $d:=d_{K_1\minter K_2}$, $d_1 := d_{K_1}$ and $d_2 := d_{K_2}$, and $\Pi := \Pi_{K_1\minter K_2}$, $\Pi_1:=\Pi_{K_1}$ and $\Pi_2:=\Pi_{K_2}$. We recall by Lemma \procref{SecondDerivativeOfDIsSymmetric} that $d$ is almost everywhere twice differentiable with symmetric second derivative. We are now interested in estimating lower bounds for $\opDet(D^2d;\langle Dd\rangle^\perp)$. There are four different cases to consider.
\subheadlabel{IntersectingConvexSets}
\proclaim{Lemma \nextprocno, {\bf Case 1}}
\noindent If $x\in\Pi^{-1}(U(K)\minter(\partial K_1)\minter K_2^o)\setminus K$ then $d=d_1$ and $\Pi=\Pi_1$.
\endproclaim
\proclabel{SecondDerivativeWhenProjectionIsInFirstConvex}
{\sl\noindent{\bf Remark:\ }Observe that this set is open, and so $Dd=Dd_1$ and $D^2d=D^2d_1$ over this set.}
\medskip
\proof Denote $y:=\Pi(x)$. Denote $\msf{N}:=(x - y)/\|x-y\|$. By Lemma \procref{CharacterisationOfClosestPoints}, $\msf{N}$ is a supporting normal to $K$ at $y$. Since $y\in\partial K_1\minter K_2^o$, By Lemma \procref{SupportingNormalIsLocalProperty}, $\msf{N}$ is also a supporting normal to $K_1$ at $y$. By Lemma \procref{CharacterisationOfClosestPoints}, $y$ minimises distance in $K_1$ to $x$. In particular, $d_1(x) = \|x-y\| = d(x)$, and $\Pi_1(x)=y=\Pi(x)$, as desired.\qed
\proclaim{Lemma \nextprocno, {\bf Case 2}}
\noindent If $x\in\Pi^{-1}(U(K)\minter(\partial K_2)\minter K_1^o)\setminus K$ then $d=d_2$ and $\Pi=\Pi_2$.
\endproclaim
\proclabel{SecondDerivativeWhenProjectionIsInSecondConvex}
\proof Denote $y:=\Pi(x)$. Denote $\msf{N}:=(x - y)/\|x-y\|$. By Lemma \procref{CharacterisationOfClosestPoints}, $\msf{N}$ is a supporting normal to $K$ at $y$. Since $y\in\partial K_2\minter K_1^o$, by Lemma \procref{SupportingNormalIsLocalProperty}, $\msf{N}$ is also a supporting normal to $K_2$ at $y$. By Lemma \procref{CharacterisationOfClosestPoints}, $y$ minimises distance in $K_2$ to $x$. In particular, $d_2(x) = \|x-y\| = d(x)$, and $\Pi_2(x)=y=\Pi(x)$, as desired.\qed
\proclaim{Lemma \nextprocno, {\bf Case 3}}
\noindent If $x\in\Pi^{-1}(U(K)\minter(\partial K_1)\minter(\partial K_2))\setminus K$, if $(\msf{N}_1\circ\Pi)(x) = (\msf{N}_2\circ\Pi)(x)$, and if $d$ is twice differentiable at $x$, then, for every vector $V$,
$$
D^2d(x)(V,V) \geq \opMax(D^2d_1(x)(V,V),D^2d_2(x)(V,V)).
$$
\endproclaim
\proclabel{SecondDerivativeLowerBoundWhenNormalsMatch}
\proof Denote $y:=\Pi(x)$. By Lemma \procref{CharacterisationOfClosestPoints}, $Dd(x)$ is a supporting normal to $K$ at $y$. By Theorem \procref{SupportingNormalsOfIntersection}, the set of supporting normals to $K$ at $y$ is the convex hull of $\left\{\msf{N}_1(y),\msf{N}_2(y)\right\}$. Since these two points coincide, this convex hull consists of a single point, and so $Dd(x)=\msf{N}_1(y)=\msf{N}_2(y)$. In particular, $Dd(x)$ is also a supporting normal to both $K_1$ and $K_2$ at $y$. It follows from Lemma \procref{CharacterisationOfClosestPoints} that $y$ minimises distance in $K_1$ to $x$. In particular, $d(x)=d_1(x)$. However, since $K\subseteq K_1$, for all $y$,
$$
d(y) = \minf_{z\in K}\|y-z\| \geq \minf_{z\in K_1}\|y-z\| = d_1(y).
$$
Thus, by differentiating, for all vectors $V$,
$$
D^2d(x)(V,V) \geq D^2d_1(x)(V,V).
$$
In like manner, we show that $D^2d(x)(V,V)\geq D^2d_2(x)(V,V)$, and so
$$
D^2d(x)(V,V) \geq \opMax(D^2d_1(V,V),D^2d_2(V,V)),
$$
as desired.\qed
\medskip
\noindent Before treating the fourth case, we require the following preliminary result.
\proclaim{Lemma \nextprocno}
\noindent If $x\in\Pi^{-1}(U(K)\minter(\partial K_1)\minter (\partial K_2))\setminus K$ and if $(\msf{N}_1\circ\Pi)(x)\neq(\msf{N}_2\circ\Pi)(x)$, then there exists a unique $s\in[0,1]$ such that
$$
Dd(x) = \frac{1-s}{l}(\msf{N}_1\circ\Pi)(x) + \frac{s}{l}(\msf{N}_2\circ\Pi)(x),
$$
where
$$
l=\|(1-s)(\msf{N}_1\circ\Pi)(x) + s(\msf{N}_2\circ\Pi)(x)\|.
$$
In particular,
$$
l \geq \frac{1}{2}\|(\msf{N}_1\circ\Pi)(x) + (\msf{N}_2\circ\Pi)(x)\|.
$$
\endproclaim
\proclabel{NormalLiesBetweenTheNormalsOfEachHypersurface}
\proof Denote $y:=\Pi(x)$. By Lemma \procref{CharacterisationOfClosestPoints}, $Dd(x)$ is a supporting normal to $K$ at $\Pi(x)$. By Theorem \procref{SupportingNormalsOfIntersection}, the set of supporting normals to $K$ at $\Pi(x)$ is the convex hull of $\left\{\msf{N}_1(y),\msf{N}_2(y)\right\}$. This coincides with the great-circular arc joining $\msf{N}_1(y)$ to $\msf{N}_2(y)$ (c.f. Section \subheadref{ConvexSubsetsOfTheSphere}), and the first assertion follows. Now observe that the vectors $\msf{N}_1(y)+\msf{N}_2(y)$ and $\msf{N}_1(y)-\msf{N}_2(y)$ are orthogonal. Thus
$$\eqalign{
l^2 &= \|(1-s)\msf{N}_1(y) + s\msf{N}_2(y)\|^2\cr
&= \|\frac{1}{2}(\msf{N}_1(y) + \msf{N}_2(y)) + \frac{(1-2s)}{2}(\msf{N}_1(y)-\msf{N}_2(y))\|^2\cr
&= \frac{1}{4}\|\msf{N}_1(y) + \msf{N}_2(y)\|^2 + \frac{(1-2s)^2}{4}\|\msf{N}_1(y)-\msf{N}_2(y)\|^2\cr
&\geq \frac{1}{4}\|\msf{N}_1(y) + \msf{N}_2(y)\|^2,\cr}
$$
and the second assertion follows. This completes the proof.\qed
\proclaim{Lemma \nextprocno, {\bf Case 4a}}
\noindent If $x\in\Pi^{-1}(U(K)\minter(\partial K_1)\minter (\partial K_2))\setminus K$, if $(\msf{N}_1\circ\Pi)(x)\neq(\msf{N}_2\circ\Pi)(x)$, and if $d$ is twice differentiable at $x$, then for every vector $V$ and for every vector $W$ which is orthogonal to both $(\msf{N}_1\circ\Pi)(x)$ and $(\msf{N}_2\circ\Pi)(x)$,
$$
D^2d(x)(V,W) = \frac{1-s}{l}\langle A_1(\Pi(x))D\Pi(x)V,W\rangle + \frac{s}{l}\langle A_2(\Pi(x))D\Pi(x)V,W\rangle,
$$
where $s$ and $l$ are as in Lemma \procref{NormalLiesBetweenTheNormalsOfEachHypersurface}.
\endproclaim
\proclabel{SecondDerivativesCaseIVA}
{\sl\noindent{\bf Remark:\ }Upon applying an isometry, we may suppose that the linear span of $\left\{e_n,e_{n+1}\right\}$ coincides with that of $\left\{(\msf{N}_1\circ\Pi)(x),(\msf{N}_2\circ\Pi)(x)\right\}$. Consequently, when $D^2d(x)$ is symmetric, this result determines every component of $D^2d(x)$ except $D^2d(x)(e_i,e_j)$ for $(i,j)\in\left\{n,n+1\right\}^2$.}
\medskip
\proof Denote $y:=\Pi(x)$. Let $V$ be a vector in $\Bbb{R}^{n+1}$. Define $\gamma:\Bbb{R}\rightarrow\Bbb{R}^{n+1}$ by $\gamma(t) = x + tV$. Let $(t_m)_\minn$ be a sequence of points in $\Bbb{R}$ converging to $0$. For all $m$, denote $x_m:=\gamma(t_m)$ and $y_m:=(\Pi\circ\gamma)(t_m)$. Upon extracting a subsequence, we may suppose that one of the following holds.
\medskip
{\bf\noindent 1:\ }$x_m\in\Pi^{-1}(U(K)\minter(\partial K_1)\minter K_2^o)\setminus K$ for all $m$. By Lemma \procref{SecondDerivativeWhenProjectionIsInFirstConvex}, for all $m$, $Dd(x_m)=Dd_1(x_m)$. Taking limits and bearing in mind Lemma \procref{DistanceAndProjectionAreSmooth}, it follows that
$$
Dd(x) = Dd_1(x) = (\msf{N}_1\circ\Pi)(x) = \msf{N}_1(y).
$$
In particular, $s=0$ and $l=1$. Moreover, for all vectors $W$ and for all $m$,
$$
\frac{1}{t_m}\langle Dd(x_m) - Dd(x),W\rangle= \frac{1}{t_m}\langle\msf{N}_1(y_m) - \msf{N}_1(y),W\rangle,
$$
so that, by the chain rule, upon taking limits, we obtain
$$
D^2d(x)(V,W) = \langle A_1(y)D\Pi(x)V,W\rangle,
$$
as desired.
\medskip
{\bf\noindent 2:\ }$x_m\in\Pi^{-1}(U(K)\minter(\partial K_2)\minter K_1^o)\setminus K$ for all $m$. As in Step $(1)$, we show that $s=1$, $l=1$ and
$$
D^2d(x)(V,W) = \langle A_2(y)D\Pi(x)V,W\rangle,
$$
as desired.
\medskip
{\bf\noindent 3:\ }$x_m\in\Pi^{-1}(U(K)\minter(\partial K_1)\minter (\partial K_2))\setminus K$ for all $m$. For all $m$, denote $\msf{N}_{1,m}:=\msf{N}_1(y_m)$ and $\msf{N}_{2,m}:=\msf{N}_2(y_m)$. Observe that, for sufficiently large $m$, $\msf{N}_{1,m}\neq\msf{N}_{2,m}$. Thus, by Lemma \procref{NormalLiesBetweenTheNormalsOfEachHypersurface}, for all $m\in\Bbb{N}$, there exists a unique $s_m\in[0,1]$ such that
$$
Dd(x_m) = \frac{1 - s_m}{l_m}\msf{N}_{1,m} + \frac{s_m}{l_m}\msf{N}_{2,m},
$$
where $l_m := \|(1-s_m)\msf{N}_{1,m} + s_m\msf{N}_{2,m}\|$. Since $\msf{N}_1$, $\msf{N}_2$, $Dd$ and $\Pi$ are continuous, $(s_m)_\minn$ and $(l_m)_\minn$ converge to the limits $s_\infty$ and $l_\infty$ respectively. Let $W$ be a vector normal to both $\msf{N}_1(y)$ and $\msf{N}_2(y)$. In particular $W$ is normal to $Dd(x)$. For all $m$,
$$\eqalign{
\frac{1}{t_m}\langle Dd(x_m) - Dd(x),W\rangle &= \frac{1}{t_m}\langle Dd(x_m),W\rangle\cr
&=\frac{1-s_m}{l_mt_m}\langle\msf{N}_{1,m},W\rangle + \frac{1-s_m}{l_mt_m}\langle\msf{N}_{2,m},W\rangle\cr
&=\frac{1-s_m}{l_mt_m}\langle\msf{N}_{1,m} - \msf{N}_1(y),W\rangle + \frac{1-s_m}{l_mt_m}\langle\msf{N}_{2,m} - \msf{N}_2(y),W\rangle.\cr
}$$
By the chain rule, upon taking limits, we obtain
$$
D^2d(x)(V,W) = \frac{1-s}{l}\langle A_1(y)D\Pi(x)V,W\rangle + \frac{s}{l}\langle A_2(y)D\Pi(x)V,W\rangle,
$$
as desired.\qed
\proclaim{Lemma \nextprocno, {\bf Case 4b}}
\noindent If $x\in\Pi^{-1}(U(K)\minter(\partial K_1)\minter (\partial K_2))\setminus K$, if $(\msf{N}_1\circ\Pi)(x)\neq(\msf{N}_2\circ\Pi)(x)$, and if $d$ is twice differentiable at $x$, then for every vector $V$ and for every vector $W$ which is tangent to the linear span of $\left\{(\msf{N}_1\circ\Pi)(x),(\msf{N}_2\circ\Pi)(x)\right\}$ and normal to $Dd(x)$,
$$
D^2d(x)(V,W) = \frac{1}{d(x)}\langle V,W\rangle.
$$
\endproclaim
\proclabel{SecondDerivativesCaseIVB}
{\sl\noindent{\bf Remark:\ }Upon applying an isometry, we may suppose that $e_{n+1}=Dd(x)$ and that the linear span of $\left\{e_n,e_{n+1}\right\}$ coincides with that of $\left\{(\msf{N}_1\circ\Pi)(x),(\msf{N}_2\circ\Pi)(x)\right\}$. Consequently when $D^2d(x)$ is symmetric, this result along with Lemma \procref{SecondDerivativesCaseIVA} determines every component of $D^2d(x)$ except for $D^2d(x)(e_{n+1},e_{n+1})$. In fact, we readily show that $D^2d(x)(e_{n+1},e_{n+1})=0$, but since this is not necessary for our purposes, we leave it as an easy exercise for the interested reader.}
\medskip
\proof Denote $y:=\Pi(x)$. By Theorem \procref{SupportingNormalsOfIntersection}, the set of supporting normals to $K$ at $y$ coincides with the convex hull of $\left\{\msf{N}_1(y),\msf{N}_2(y)\right\}$, which in turn coincides with the great-circular arc joining $\msf{N}_1(y)$ to $\msf{N}_2(y)$ (c.f. Section \subheadref{ConvexSubsetsOfTheSphere}). Denote this great-circular arc by $\msf{N}:[0,1]\rightarrow\Sigma$. In particular, for all $r$, $\msf{N}(r)$ is a supporting normal to $K$ at $y$. We define $\gamma:[0,1]\rightarrow\Bbb{R}^{n+1}$ by $\gamma(r) = y + d(x)\msf{N}(r)$. By Lemma \procref{CharacterisationOfClosestPoints}, for all $t$, $y$ minimises distance in $K$ to $\gamma(r)$. In particular, by Lemma \procref{FirstDerivativeOfDistance}, for all $r$, $(Dd\circ\gamma)(r)=\msf{N}(r)$. Let $s$ be as in Lemma \procref{NormalLiesBetweenTheNormalsOfEachHypersurface}. Since $W$ lies in the plane spanned by $\msf{N}_1(y)$ and $\msf{N}_2(y)$ but is normal to $Dd(x)$, upon multiplying by a scalar factor, we may suppose that $W=(\partial_r\gamma)(s)$. Thus
$$\eqalign{
D^2d(x)(W,V) &= \langle \partial_r(Dd\circ\gamma)(s),V\rangle\cr
&=\langle (\partial_r\msf{N})(s),V\rangle\cr
&=\frac{1}{d(x)}\langle(\partial_r\gamma)(s),V\rangle\cr
&=\frac{1}{d(x)}\langle W,V\rangle,\cr}
$$
as desired.\qed
\proclaim{Lemma \nextprocno}
\noindent For every compact subset $X$ of $U$ there exists $C>0$ with the property that for every $x\in(\partial K_1)\minter(\partial K_2)\minter X$,
$$
\|\msf{N}_1(x) + \msf{N}_2(x)\| \geq \frac{1}{C}.
$$
\endproclaim
\proclabel{LowerBoundOfL}
\proof Suppose the contrary. By compactness, there exists $x\in(\partial K_1)\minter(\partial K_2)\minter X$ such that $\msf{N}_1(x)+\msf{N}_2(x)=0$. By definition of supporting normals, for all $y\in K_1\minter K_2$, $\langle y-x,\msf{N}_1(x)\rangle\leq 0$ and $\langle y-x,\msf{N}_2(x)\rangle\leq 0$. Since $\msf{N}_1(x)=-\msf{N}_2(x)$, it follows that for all $y\in K_1\minter K_2$, $\langle y-x,\msf{N}_1\rangle=\langle y-x,\msf{N}_2\rangle=0$. In other words $K_1\minter K_2$ is contained in the hyperplane normal to $\msf{N}_1=-\msf{N}_2$ passing through $x$, and therefore has trivial interior. This is absurd, and the result follows.\qed
\proclaim{Lemma \nextprocno}
\noindent For every compact subset $X$ of $U$ there exists $C>0$ with the property that if $x\in\Pi^{-1}(X\minter U(K)\minter(\partial K_1)\minter(\partial K_2))\setminus K$, if $(\msf{N}_1\circ\Pi)(x)\neq(\msf{N}_2\circ\Pi)(x)$, if $d$ is twice differentiable at $x$, and if $D^2d(x)$ is symmetric, then
$$
\|D\Pi(x) - \pi^{1,2}\| \leq Cd(x),
$$
where $\pi^{1,2}$ is the orthonogonal projection from $\Bbb{R}^{n+1}$ onto $\langle(\msf{N}_1\circ\Pi)(x),(\msf{N}_2\circ\Pi)(x)\rangle^\perp$.
\endproclaim
\proclabel{DerivativeOfClosestPointProjectionIsLikeNormalProjection}
\proof Denote $y:=\Pi(x)$. Let $C_1>0$ be such that $\|A_1(z)\|\leq C_1$ and $\|A_2(z)\|\leq C_1$ for all $z$ in $X\minter U(K_1)$ and $X\minter U(K_2)$ respectively. Let $C_2>0$ be as in Lemma \procref{LowerBoundOfL}. By Lemma \procref{TwiceDifferentiableWheneverPiIsDifferentiable}, $\Pi$ is differentiable at $x$ and, for all vectors $V$ and $W$,
$$
\langle D\Pi(x)V,W\rangle = \langle \pi(V),\pi(W)\rangle - d(x)D^2d(x)(V,W),
$$
where $\pi$ is the orthogonal projection from $\Bbb{R}^{n+1}$ onto $\langle Dd(x)\rangle^\perp$. Observe, in particular, that since $D^2d$ is symmetric, so too is $D\Pi$. Let $V$ be any vector in $\Bbb{R}^{n+1}$. Define $\gamma:\Bbb{R}\rightarrow\Bbb{R}^{n+1}$ by $\gamma(t) := x + tV$. Since $(\Pi\circ\gamma)(t)\in K$ for all $t$, it follows that for each $\msf{N}\in\left\{\msf{N}_1(y),\msf{N}_2(y)\right\}$ and for all $t$,
$$
\langle (\Pi\circ\gamma)(t) - y,\msf{N}\rangle \leq 0.
$$
By the chain rule, differentiating this relation yields
$$
\langle D\Pi(x)V,\msf{N}\rangle = 0.
$$
Thus, by linearity and symmetry, for any vector $W$ in the linear span of $\left\{\msf{N}_1(y),\msf{N}_2(y)\right\}$,
$$
\langle D\Pi(x)W,V\rangle = \langle D\Pi(x)(V),W\rangle = 0.
$$
Now let $V$ and $W$ both be orthogonal to $\langle\msf{N}_1(y),\msf{N}_2(y)\rangle$. In particular, $V$ and $W$ are both orthogonal to $Dd(x)$. Thus, by Lemma \procref{TwiceDifferentiableWheneverPiIsDifferentiable},
$$
\langle D\Pi(x)V,W\rangle = \langle V,W\rangle - d(x)D^2d(x)(V,W).
$$
Let $s$ and $l$ be as in Lemma \procref{NormalLiesBetweenTheNormalsOfEachHypersurface}. Then, by Lemma \procref{SecondDerivativesCaseIVA} and bearing in mind Lemma \procref{PiIsDifferentiable},
$$\eqalign{
\left|D^2d(x)(V,W)\right| &= \left|\frac{1-s}{l}\langle A_1(y)D\Pi(x)V,W\rangle + \frac{s}{l}\langle A_2(y)D\Pi(x)V,W)\rangle\right|\cr
&\leq \frac{1-s}{l}C_1\|V\|\|W\| + \frac{s}{l}C_1\|V\|\|W\|\cr
&\leq 2C_1C_2\|V\|\|W\|.\cr}
$$
Thus
$$
\left|\langle D\Pi(x)V,W\rangle - \langle V,W\rangle\right| \leq 2C_1C_2d(x)\|V\|\|W\|.
$$
Combining these relations, we conclude that $\|D\Pi(x) - \pi^{1,2}\|\leq 2C_1C_2d(x)$, as desired.\qed
\proclaim{Lemma \nextprocno}
\noindent For every compact subset $X$ of $U$ and for all $\epsilon>0$, there exists $r>0$ with the property that if $x\in\Pi^{-1}(X\minter U(K))\setminus K$, if $d(x)<r$ and if $D^2d(x)$ is defined and is symmetric, then
$$
\opDet(D^2d(x);\langle Dd(x)\rangle^\perp) \geq (k-\epsilon)^n.
$$
\endproclaim
\proclabel{SecondDerivativeIsInTheCorrectSet}
\proof We consider the following cases.
\medskip
{\bf\noindent 1:\ }Suppose that $x\in \Pi^{-1}(X\minter U(K)\minter (\partial K_1)\minter K_2^o)\setminus K$. By Lemma \procref{SecondDerivativeWhenProjectionIsInFirstConvex}, $D^2d(x)=D^2d_1(x)$, and the result follows by Lemma \procref{LowerBoundsOnOrthogonalComponentOfDeterminantNearHypersurface}.
\medskip
{\bf\noindent 2:\ }Suppose that $x\in \Pi^{-1}(X\minter U(K)\minter (\partial K_2)\minter K_1^o)\setminus K$. By Lemma \procref{SecondDerivativeWhenProjectionIsInSecondConvex}, $D^2d(x)=D^2d_2(x)$, and the result follows by Lemma \procref{LowerBoundsOnOrthogonalComponentOfDeterminantNearHypersurface}.
\medskip
{\bf\noindent 3:\ }Suppose that $x\in \Pi^{-1}(X\minter U(K)\minter (\partial K_1)\minter (\partial K_2))\setminus K$ and $\msf{N}_1(\Pi(x))=\msf{N}_2(\Pi(x))$. By Lemma \procref{SecondDerivativeLowerBoundWhenNormalsMatch}, for all vectors $V\in\Bbb{R}^{n+1}$,
$$
D^2d(x)(V,V)\geq\opMax(D^2d_1(x)(V,V),D^2d_2(x)(V,V)).
$$
In particular, bearing in mind that $Dd(x)=Dd_1(x)=Dd_2(x)$,
$$
\opDet(D^2d(x);\langle Dd(x)\rangle^\perp)\geq\opDet(D^2d_1(x);\langle Dd_1(x)\rangle^\perp),\opDet(D^2d_2(x);\langle Dd_2(x)\rangle^\perp),
$$
and the result now follows by Lemma \procref{LowerBoundsOnOrthogonalComponentOfDeterminantNearHypersurface}.
\medskip
{\bf\noindent 4:\ }Suppose that $x\in \Pi^{-1}(X\minter U(K)\minter (\partial K_1)\minter (\partial K_2))\setminus K$ and that $\msf{N}_1(\Pi(x))\neq\msf{N}_2(\Pi(x))$. Denote $y:=\Pi(x)$. Let $s$ and $l$ be as in Lemma \procref{NormalLiesBetweenTheNormalsOfEachHypersurface}. Let $C_1>1$ be such that $(1/C_1)\opId\leq A_1(y)\leq C_1\opId$ and $(1/C_1)\opId\leq A_2(y)\leq C_1\opId$ for all $y$ in $X\minter U(K_1)$ and $X\minter U(K_2)$ respectively. Let $C_2\geq 0$ be as in Lemma \procref{DerivativeOfClosestPointProjectionIsLikeNormalProjection}. Define $r:=1/(2C_1^2C_2)$. Then if $d(x)<r$, for all vectors $V$ normal to $\msf{N}_1(y)$ and $\msf{N}_2(y)$,
$$\eqalign{
\langle A_1(y)D\Pi(x)V,V\rangle &= \langle A_1(y)V,V\rangle + \langle A_1(y)(D\Pi(x) - \pi^{1,2})(V),V\rangle\cr
&\geq \frac{1}{C_1}\|V\|^2 - \frac{1}{2C_1}\|V\|^2\cr
&=\frac{1}{2C_1}\|V\|^2.\cr}$$
Likewise, for all such $x$ and $V$,
$$
\langle A_2(y)D\Pi(x)V,V\rangle \geq \frac{1}{2C_1}\|V\|^2.
$$
Thus, if $s$ and $l$ are as in Lemma \procref{NormalLiesBetweenTheNormalsOfEachHypersurface}, by Lemma \procref{SecondDerivativesCaseIVA}, for all such $x$ and $V$,
$$\eqalign{
D^2d(x)(V,V)&= \frac{1-s}{l}\langle A_1(y)D\Pi(x)V,V\rangle + \frac{s}{l}\langle A_2(y)D\Pi(x)V,V\rangle\cr
&\geq \frac{1-s}{2C_1l}\|V\|^2 + \frac{s}{2C_1l}\|V\|^2\cr
&\geq \frac{1}{2C_1}\|V\|^2.\cr}
$$
Upon applying an isometry, we may suppose that the plane spanned by $e_n$ and $e_{n+1}$ coincides with the plane spanned by $\msf{N}_1(y)$ and $\msf{N}_2(y)$ and furthermore that $e_{n+1}=Dd(x)$. We denote by $M$ the restriction of $D^2d(x)$ to $\langle e_1,...,e_{n-1}\rangle$. By the preceeding discussion, $M\geq (1/2C_1)\opId$. However, by Lemma \procref{SecondDerivativesCaseIVB}, for all $i$,
$$
D^2d(x)(e_i,e_n) = D^2d(x)(e_n,e_i) = \frac{1}{d(x)}\delta_{in}.
$$
Reducing $r$ if necessary, we may suppose that $r<(2C_1)^{1-n}(k-\epsilon)^{-n}$ so that, if $d(x)<r$, then
$$
\opDet(D^2d(x),\langle Dd(x)\rangle^\perp)\geq (k-\epsilon)^n,
$$
as desired.\qed
\newsubhead{Smoothing functions and convexity}
Let $\chi\in C_0^\infty(\Bbb{R}^{n+1})$ be a smooth, non-negative function such that $\chi=0$ outside the unit ball $B_1(0)$, and
\subheadlabel{SmoothingOperators}
$$
\int_{\Bbb{R}^{n+1}}\chi(x)\opdVol_x = 1
$$
For all $s>0$, we define $\chi_s\in C_0^\infty(\Bbb{R}^{n+1})$ by
$$
\chi_s(x) := s^{-(n+1)}\chi(x/s).
$$
Let $E$ be a finite-dimensional vector space. For any function $f\in L^1_\oploc(\Bbb{R}^{n+1},E)$, and for all $s>0$, we define the function $f_s:\Bbb{R}^{n+1}\rightarrow E$ by,
$$
f_s(x) := \int_{\Bbb{R}^{n+1}}f(x-y)\chi_s(y)\opdVol_y.
$$
We recall the following properties of smoothing functions.
\proclaim{Lemma \nextprocno}
\noindent For all $f\in L^1_\oploc(\Bbb{R}^{n+1})$ and for all $s>0$, $f_s$ is continuous.
\endproclaim
\proclabel{MollifiedFunctionIsContinuous}
{\sl\noindent{\bf Remark:\ }In fact, as is well known, $f_s$ is smooth.}
\medskip
\proof Choose $x\in\Bbb{R}^{n+1}$ and $s>0$. By local uniform continuity, there exists $\delta>0$ such that if $\|y\|<s$ and $\|z-y\|<\delta$, then $\left|\chi_s(y) - \chi_s(z)\right|<\epsilon$. Thus, if $\|z-x\|<\delta$, using a change of variable, we obtain
$$\eqalign{
\|f_s(z) - f_s(x)\| &= \|\int_{\Bbb{R}^{n+1}}f(z-y)\chi_s(y) - f(x-y)\chi_s(y)\opdVol_y\|\cr
&=\|\int_{\Bbb{R}^{n+1}}f(x-y)(\chi_s(y + (z-x)) - \chi_s(y))\opdVol_y\|\cr
&\leq \int_{\Bbb{R}^{n+1}}\|f(x-y)\|\left|\chi_s(y + (z-x)) - \chi_s(y)\right|\opdVol_y\cr
&\leq \epsilon\int_{B_{R+\delta}(0)}\|f(x-y)\|\opdVol_y.\cr}
$$
Since $\epsilon$ may be chosen arbitrarily small, continuity of $f_s$ at $x$ follows. Since $x\in\Bbb{R}^{n+1}$ is arbitrary, it follows that $f_s$ is continuous, as desired.\qed
\proclaim{Lemma \nextprocno}
\noindent Choose $f\in L^1_\oploc(\Bbb{R}^{n+1})$. If $f$ is continuous, then $(f_s)_{s>0}$ converges to $f$ locally uniformly as $r$ tends to $0$.
\endproclaim
\proclabel{MollifiedFunctionsConvergeUniformly}
\proof Choose $\epsilon>0$ and $R>0$. By uniform continuity, there exists $\delta>0$ such that if $\|x\|<R$ and if $\|x-y\|<\delta$, then $\|f(x) - f(y)\|<\epsilon$. Then, for $s<\delta$ and for $\|x\|<R$, bearing in mind that $\chi_s$ is non-negative and has integral equal to $1$,
$$\eqalign{
\|f(x) - f_s(x)\| &= \| f(x) - \int_{\Bbb{R}^{n+1}}f(x-y)\chi_s(y)\opdVol_y\|\cr
&= \|\int_{\Bbb{R}^{n+1}}(f(x) - f(x-y))\chi_s(y)\opdVol_y\|\cr
&\leq \int_{\Bbb{R}^{n+1}}\|f(x) - f(x-y)\|\chi_s(y)\opdVol_y\cr
&\leq \epsilon\int_{\Bbb{R}^{n+1}}\chi_s(y)\opdVol_y\cr
&=\epsilon.\cr}
$$
Since $R,\epsilon>0$ are arbitrary, we conclude that $(f_s)_{s>0}$ converges locally uniformly to $f$ over $\Bbb{R}^{n+1}$ as $s$ tends to $0$, as desired.\qed
\proclaim{Lemma \nextprocno}
\noindent Choose $f\in L^1_\oploc(\Bbb{R}^{n+1})$. If $f$ has $L^1_\oploc$ distributional derivatives, then, for all $s>0$, $f_s$ is differentiable and $D(f_s) = (Df)_s$.
\endproclaim
\proclabel{MollifierOfDerivative}
\proof Choose $x\in\Bbb{R}^{n+1}$ and $s>0$. Choose $\epsilon>0$. Since $\chi_s$ is smooth, there exists $\eta>0$ with the property that for all $y$ and for all vectors $V$ such that $\|V\|\leq\eta$,
$$
\left|\chi_s(y + V) - \chi_s(y) - D\chi_s(y)V\right| \leq \epsilon\|V\|.
$$
Thus, using the definition of the distributional derivative and a change of variable, for all $V$ such that $\|V\|\leq\eta$, we obtain
$$\eqalign{
\|f_s(x + V) - f_s(x) - (Df)_s(x)V\|
&= \|\int_{\Bbb{R}^{n+1}}f(x + V - y)\chi_s(y)  - f(x-y)\chi_s(y)\cr
&\qquad - Df(x-y)V\chi_s(y)\opdVol_y\|\cr
&=\|\int_{\Bbb{R}^{n+1}}f(x-y)(\chi_s(y+V)\cr
&\qquad- \chi_s(y) - D\chi_s(y)V)\opdVol_y\|\cr
&\leq \epsilon\|V\|\int_{B_s(x)}\|f(x-y)\|\opdVol_y.\cr}
$$
Since $V$ is arbitrary, and since $\epsilon$ may be chosen arbitrarily small, we conclude that $f_s$ is differentiable at $x$ with derivative equal to $(Df)_s(x)$, as desired.\qed
\medskip
\noindent Combining these results yields
\proclaim{Theorem \nextprocno}
\noindent If $f\in L^1_\oploc(\Bbb{R}^{n+1})$ is $C^k$, then for all $s>0$, $J^k(f_s) = (J^kf)_s$, and $(f_s)_{s>0}$ converges to $f$ in the $C^k_\oploc$ sense as $s$ tends to $0$.
\endproclaim
\proof We work by induction on $k$. By Lemmas \procref{MollifiedFunctionIsContinuous} and \procref{MollifiedFunctionsConvergeUniformly}, the result holds when $k=0$. Suppose that the result holds for $k=l$. Choose $f\in L^1_\oploc(\Bbb{R}^{n+1})$ such that $f$ is $C^{l+1}$. In particular, $Df$ is $C^l$. Using the induction hypothesis together with Lemma \procref{MollifierOfDerivative}, we obtain, for all $s$,
$$
J^{l+1}(f_s)=J^l(D(f_s))=J^l((Df)_s) = (J^l(Df))_s = (J^{l+1}f)_s.
$$
Moreover, $(f_s)_{s>0}$ and $(J^l(Df)_s)_{s>0}$ converge locally uniformly to $f$ and $J^l(Df)$ respectively as $s$ tends to $0$. $(J^{l+1}(f_s))_{s>0}$ therefore converges locally uniformly to $J^{l+1}f$, and the result follows by induction.\qed
\medskip
Importantly, the smoothing operation preserves preserves convexity.
\proclaim{Lemma \nextprocno}
\noindent If $f:\Bbb{R}^{n+1}\rightarrow\Bbb{R}$ is convex, then so too is $f_s$ for all $s>0$.
\endproclaim
\proclabel{SmoothingPreservesConvexity}
\proof Fix $s>0$. Using the convexity of $f$ and the positivity of $\chi_s$, for all $x,y\in\Bbb{R}^{n+1}$ and for all $t\in[0,1]$, we obtain
$$\eqalign{
f_s(tx+(1-t)y)
&=\int_{B_s(0)}f(tx + (1-t)y + z)\chi_s(z)\opdVol_z\cr
&=\int_{B_s(0)}f(t(x+z) + (1-t)(y+z))\chi_s(z)\opdVol_z\cr
&\geq\int_{B_s(0)}\left(tf(x+z) + (1-t)f(y+z)\right)\chi_s(z)\opdVol_z\cr
&=tf_s(x) + (1-t)f_s(y),\cr}
$$
so that $f_s$ is convex, as desired.\qed
\medskip
\noindent Of particular use to us is
\proclaim{Lemma \nextprocno}
\noindent Let $E$ be a finite dimensional vector space. Let $K$ be a closed convex subset of $E$. Let $U$ be an open subset of $\Bbb{R}^{n+1}$. Let $f\in L^1_\oploc(\Bbb{R}^{n+1})$ be such that for almost all $x\in U$, $f(x)\in K$. Then for all $s>0$ and for all $x$ with the property that $B_s(x)\subseteq U$, we have $f_s(x)\in K$.
\endproclaim
\proclabel{MollifiedFunctionRemainsWithinConvexSet}
\proof We use the terminology of Section \headref{WeakBarriers}. Let $H(\msf{N},t)$ be an open half-space of $E$ containing $K$. In particular, $\langle z,\msf{N}\rangle<t$ for all $z\in K$. Choose $s>0$ and $x\in\Bbb{R}^{n+1}$ such that $B_s(x)\subseteq U$. Then, bearing in mind that $\chi$ is positive,
$$\eqalign{
\langle f_s(x),\msf{N}\rangle &= \langle \int_{B_s(0)}f(x-y)\chi_s(y)\opdVol_y,\msf{N}\rangle\cr
&=\int_{B_s(0)}\langle f(x-y),\msf{N}\rangle\chi_s(y)\opdVol_y\cr
&<\int_{B_s(0)}t\chi_s(y)\opdVol_y\cr
&=t,\cr}
$$
so that $f_s(x)\in H(\msf{N},t)$. Since $H(\msf{N},t)$ is an arbitrary open half-space containing $K$, we conclude that $f_s(x)\in \opConv(K)=K$, and the result follows.\qed
\proclaim{Corollary \nextprocno}
\noindent If $f:\Bbb{R}^{n+1}\rightarrow\Bbb{R}$ is $1$-Lipschitz, then for all $s>0$,
$$
\|f-f_s\|_0 \leq s.
$$
\endproclaim
\proclabel{MollifierOfDistanceFunction}
\proof Choose $x\in\Bbb{R}^{n+1}$. Since $f$ is $1$-Lipschitz, for all $y\in B_s(x)$, $f(y)\in[f(x)-s,f(x)+s]$. Since this interval is convex, by Lemma \procref{MollifiedFunctionRemainsWithinConvexSet}, $f_s(x)\in[f(x)-s,f(x)+s]$ so that $\left|f(x)-f_s(x)\right|\leq s$. Since $x\in\Bbb{R}^{n+1}$ is arbitrary, it follows that $\|f-f_s\|_0\leq s$, as desired.\qed
\newsubhead{Smoothing the intersection}
We return to the situation discussed in Section \subheadref{IntersectingConvexSets}. Thus, let $K_1$ and $K_2$ be compact, convex subsets of $\Bbb{R}^{n+1}$ whose intersection has non-trivial interior. Let $U$ be an open subset of $\Bbb{R}^{n+1}$ and suppose that both $U(K_1)$ and $U(K_2)$ are smooth of gaussian curvature at least $k$. As before, we denote $K:=K_1\minter K_2$ and we denote $d:=d_{K_1\minter K_2}$, $d_1:=d_{K_1}$ and $d_2:=d_{K_2}$. We recall the following version of the submersion theorem.
\subheadlabel{SmoothingTheIntersection}
\proclaim{Lemma \nextprocno}
\noindent Let $U\subseteq\Bbb{R}^{n+1}$ be an open set. Let $f:U\rightarrow\Bbb{R}$ be a smooth mapping and denote $\Sigma=f^{-1}(\left\{0\right\})$. If $0$ is a regular value of $f$, then $\Sigma$ is a smooth, embedded submanifold. Moreover, for all $x\in\Sigma$, $Df(x)/\|Df(x)\|$ is a unit, normal vector field over $\Sigma$, and if we denote by $A$ the shape operator of $\Sigma$ with respect to this normal, then for all $x\in\Sigma$ and for all $X,Y$ tangent to $\Sigma$ at $x$,
$$
\langle A(x)X,Y\rangle = \frac{1}{\|Df(x)\|}D^2f(x)(X,Y).
$$
\endproclaim
\proclabel{SubmersionTheoremWithShapeOperator}
\proof If $0$ is a regular value of $f$, then it follows by the submersion theorem (c.f. \cite{GuillemanPollack}) that $\Sigma$ is a smooth, embedded submanifold of $U$. Choose $x\in U$ and let $X$ be a tangent vector to $\Sigma$ at $x$. Let $\gamma:]-\epsilon,\epsilon[\rightarrow\Sigma$ be a smooth curve such that $\gamma(0)=x$ and $\gamma'(0)=X$. In particular, $(f\circ\gamma)(t)=0$ for all $t$. Thus, by the chain rule,
$$
\langle Df(x),X\rangle=\langle Df(x),\gamma'(0)\rangle = (f\circ\gamma)'(0)=0.
$$
Since $X$ is an arbitrary vector tangent to $\Sigma$ at $x$, it follows that $Df(x)$ is normal to $\Sigma$ at $x$. Since, furthermore, $\|Df(x)\|\neq 0$, we conclude that $Df(x)/\|Df(x)\|$ is a unit normal vector to $\Sigma$ at $x$, as desired. Now let $X$ and $Y$ be tangent vectors to $\Sigma$ at $x$. We denote $\msf{N}=Df/\|Df\|$. By definition of $A$, and using the chain and product rules,
$$\eqalign{
\langle A(x)X,Y\rangle &= \langle D\msf{N}(x)X,Y\rangle\cr
&= \langle D(Df/\|Df\|)(x)X,Y\rangle\cr
&= \frac{1}{\|Df(x)\|} D^2f(x)(X,Y) - \frac{1}{\|Df(x)\|^3}\langle Df(x),Y\rangle\langle Df(x),X\rangle.\cr}
$$
However, by the previous discussion, $Df(x)$ is normal to $\Sigma$ at $x$, and so,
$$
\langle A(x)X,Y\rangle = \frac{1}{\|Df(x)\|}D^2f(x)(X,Y),
$$
as desired.\qed
\medskip
For all $k,B>0$, and for all $\msf{N}\in\Sigma^n$, we define the set $\kappa(k,B,\msf{N})\subseteq\opSymm(2,\Bbb{R}^{n+1})$ by
$$
\kappa(k,B,\msf{N}) := \left\{ A\ |\ \|A\|\leq B,\ A\geq 0,\ \opDet(A;\langle\msf{N}\rangle^\perp)\geq k^n\right\}.
$$
\proclaim{Lemma \nextprocno}
\noindent For all $k,B>0$ and for all $\msf{N}\in\Sigma^n$, $\kappa(k,B,\msf{N})$ is compact and convex.
\endproclaim
\proclabel{ThatSetKappaIsConvex}
\proof The set of all matrices of norm no greater than $B$ is compact. Since $\kappa(k,B,\msf{N})$ is a closed subset of this set, it too is compact. Observe that the space of positive-definite matrices is convex. Furthermore, by Lemma \procref{ConcavityOfF}, the function $\opDet(M;\langle\msf{N}\rangle^\perp)^\frac{1}{n}$ is convex over this space. Since the norm is also convex, we conclude that $\kappa(k,B,N)$ is convex, as desired.\qed
\proclaim{Lemma \nextprocno}
\noindent For every compact subset $X$ of $U$ and for all $\epsilon>0$, there exists $\rho>0$ with the property that if $x\in\Pi^{-1}(X\minter U(K))\setminus K$, if $d(x)<\rho$ and if $D^2d(x)$ is defined and is symmetric, then
$$
D^2d(x) \in \kappa(k-\epsilon,2/d(x),Dd(x)).
$$
\endproclaim
\proclabel{DistanceFunctionHasHessianInCorrectSet}
\proof By Lemma \procref{DIsTwiceDifferentiable}, $\|D^2d(x)\|\leq 2/d(x)$ and the result follows by Lemma \procref{SecondDerivativeIsInTheCorrectSet}.\qed
\medskip
We now consider smoothings of $d$ as described in Section \subheadref{SmoothingOperators}.
\proclaim{Lemma \nextprocno}
\noindent For every compact subset $X$ of $U$ and for all $\epsilon>0$, there exists $\rho>0$ with the property that for all $r\in]0,\rho[$, there exists $\sigma>0$ such that if $s<\sigma$, if $x\in X$ and if $d_s(x)=r$, then $0<\|Dd_s(x)\|\leq 1$ and
$$
D^2d_s(x) \in \kappa(k-\epsilon,4/r,Dd_s(x)/\|Dd_s(x)\|).
$$
\endproclaim
\proclabel{DerivativeOfMollifiedFunctionIsInCorrectSet}
\proof Choose $\sigma_1>0$ such that $X_1:=\overline{B}_{\sigma_1}(X)\subseteq U$. Since $X$ is compact, so too is $X_1$. We first claim that there exists a compact subset $X_2$ of $U$ and $\rho_1>0$ such that
$$
X_1\minter d^{-1}(]0,\rho_1[)\subseteq\Pi^{-1}(X_2\minter U(K))\setminus K.
$$
Indeed, suppose the contrary. There exists a sequence $(x_m)_\minn$ in $X_1$ with the properties that $d(x_m)>0$ for all $m$, $(d(x_m))_\minn$ converges to $0$ and $(\Pi(x_m))_\minn$ is not contained in any compact subset of $U$. For all $m$, denote $y_m:=\Pi(x_m)$ and $\msf{N}_m:=Dd(x_m)$. Since $K$ is compact, there exists $y_\infty\in K$ towards which $(y_m)_\minn$ subconverges. By hypothesis, $y_\infty$ lies in the boundary of $U$. By Lemma \procref{FirstDerivativeOfDistance}, for all $m$, $x_m=y_m + d(x_m)\msf{N}_m$. In particular, since $(d(x_m))_\minn$ converges to $0$ and since $\msf{N}_m$ has unit length for all $m$, it follows that $(x_m)_\minn$ also subconverges to $y_\infty$. By compactness, $y_\infty$ is also an element of $X_1$, which is absurd, and the assertion follows.
\par
By Lemma \procref{DistanceFunctionHasHessianInCorrectSet}, there exists $\rho_2<\rho_1$ with the property that if $x\in\Pi^{-1}(X_2\minter U(K))\setminus K$, if $d(x)<\rho_2$ and if $D^2d(x)$ is defined and is symmetric, then
$$
D^2d(x) \in \kappa(k-\epsilon/4,2/d(x),Dd(x)).
$$
Let $\delta\in]0,1[$ be such that if $\msf{N}$ is any vector in $\Sigma^n$ and if $V\in\overline{B}_\delta(\msf{N})$, then
$$
\kappa(k-\epsilon/4,4/r,V/\|V\|) \subseteq \kappa(k-\epsilon/2,4/r,\msf{N}) \subseteq \kappa(k-\epsilon,4/r,V/\|V\|).
$$\par
Choose $r\in[0,\rho_2[$. Fix $\eta>0$ such that $2\eta<\opMin(r/2,\rho_2-r)$. Since $K$ is compact, so too is $d^{-1}([r-\eta,r+\eta])$. By continuity, there therefore exists $\sigma_2<\opMin(\eta,\sigma_1)$ such that if $x\in d^{-1}([r-\eta,r+\eta])$ and if $y\in B_{\sigma_2}(x)$, then $Dd(y)\subseteq\overline{B}_\delta(Dd(x))$. Now choose $s<\sigma_2$. Fix $x\in d_s^{-1}(\left\{r\right\})\minter X$. By Corollary \procref{MollifierOfDistanceFunction}, $x\in d^{-1}([r-\eta,r+\eta])$. Thus, if $y\in B_{\sigma_2}(x)$, then $d(y)\in]r/2,\rho_2[$ and $Dd(y)\in\overline{B}_\delta(Dd(x))$. Furthermore, every such $y$ is an element of $X_1\minter d^{-1}(]0,\rho_2[)\subseteq\Pi^{-1}(X_2\minter U(K))\setminus K$, so that if $D^2d(y)$ is defined and symmetric, then
$$
D^2d(y) \in \kappa(k-\epsilon/4,4/r,Dd(y)) \subseteq \kappa(k-\epsilon/2,4/r,Dd(x)).
$$
Since $\overline{B}_\delta(Dd(x))$ is compact and convex, it follows by Lemma \procref{MollifiedFunctionRemainsWithinConvexSet} that
$$
Dd_s(x) \in \overline{B}_{\delta}(Dd(x)),
$$
and, in particular, $Dd_s(x)\neq 0$. Likewise, since $\kappa(k-\epsilon/4,4/r,Dd(x))$ is compact and convex, by Lemma \procref{MollifiedFunctionRemainsWithinConvexSet} again,
$$
D^2d_s(x) \in \kappa(k-\epsilon/2,4/r,Dd(x)) \subseteq \kappa(k-\epsilon,4/r,Dd_s(x)/\|Dd_s(x)\|).
$$
Finally, since $Dd(y)\in\overline{B}_1(0)$ at every point where it is defined, and since $\overline{B}_1(0)$ is closed and convex, it follows by Lemma \procref{MollifiedFunctionRemainsWithinConvexSet} again that $\|Dd_s(x)\|\leq 1$, and this completes the proof.\qed
\proclaim{Theorem \nextprocno}
\noindent For every compact subset $X$ of $U$ and for all $\epsilon>0$, there exists $\rho>0$ with the property that for all $r<\rho$, there exists $\sigma>0$ such that if $s<\sigma$, if $x\in X$ and if $d_s(x)=r$, then $d_s^{-1}(\left\{r\right\})$ is smooth near $x$ and has gaussian curvature at least $k-\epsilon$ at $x$.
\endproclaim
\proclabel{SmoothIntersectionHasHighCurvature}
\proof Let $\rho$ be as in Lemma \procref{DerivativeOfMollifiedFunctionIsInCorrectSet}. Choose $r<\rho$. Let $\sigma$ be as in Lemma \procref{DerivativeOfMollifiedFunctionIsInCorrectSet}. Choose $s<\sigma$. We denote $\Sigma_{r,s}=d_s^{-1}(\left\{r\right\})$. Choose $x\in X\minter\Sigma_{r,s}$. By Lemma \procref{DerivativeOfMollifiedFunctionIsInCorrectSet}, $Dd_s(x)\neq 0$ and $\|Dd_s(x)\|\leq 1$. Thus, by Lemma \procref{SubmersionTheoremWithShapeOperator}, $\Sigma_{r,s}$ is smooth near $x$ and $Dd_s(x)/\|Dd_s(x)\|$ is the normal to $\Sigma_{r,s}$ at $x$. Moreover, if we denote by $A(x)$ the shape operator of $\Sigma_{r,s}$ at $x$ with respect to this normal, then, for all vectors $X$ and $Y$ tangent to $\Sigma_{r,s}$ at $x$,
$$
A(x)(X,Y) = \frac{1}{\|Dd_s(x)\|}D^2d_s(x)(X,Y).
$$
Thus, bearing in mind that $\|Dd_s(x)\|\leq 1$, if we denote by $\kappa(x)$ the gaussian curvature of $\Sigma$ at $x$, then
$$
\kappa(x) = \opDet(A(x))^{1/n} \geq \opDet(D^2d(x);\langle Dd_s(x)\rangle^\perp)^{1/n}.
$$
However, by Lemma \procref{DerivativeOfMollifiedFunctionIsInCorrectSet},
$$
D^2d(x) \in \kappa(k-\epsilon,4/r,Dd_s(x)/\|Dd_s(x)\|),
$$
so that $\kappa(x)\geq k-\epsilon$, as desired.\qed
\newsubhead{Weak barriers}
Let $U\subseteq\Bbb{R}^{n+1}$ be an open set. Let $k>0$ be a positive real number. Let $K$ be a compact, convex subset of $\Bbb{R}^{n+1}$. We say that $K$ is a \gloss{strong barrier} of gaussian curvature at least $k$ inside $U$ whenever $(\partial K)\minter U$ is smooth and has gaussian curvature at least $k$ at every point. We say that $K$ is a \gloss{weak barrier} of gaussian curvature at least $k$ inside $U$ whenever there exists a sequence $(\epsilon_m)_{\minn}>0$ converging to $0$, an increasing sequence $(V_m)_\minn$ of open sets and a sequence $(K_m)_\minn$ of convex sets converging to $K$ in the Hausdorff sense with the properties that $U=\munion_{m\in\Bbb{N}}V_m$ and, for all $m$, $K_m$ is a strong barrier of gaussian curvature at least $k-\epsilon_m$ inside $V_m$. That is, weak barriers are Hausdorff limits of strong barriers.
\par
We first show that the set of weak barriers is closed in the Hausdorff topology.
\proclaim{Lemma \nextprocno}
\noindent Let $U\subseteq\Bbb{R}^{n+1}$ be an open set. Let $k>0$ be a positive real number. Let $(U_m)_\minn$ be an increasing sequence of open sets such that $U=\munion_\minn U_m$. Let $(k_m)_\minn$ be a sequence of positive real numbers converging to $k$. Let $(K_m)_\minn,K_\infty$ be compact, convex subsets of $\Bbb{R}^{n+1}$ and suppose that $(K_m)_\minn$ converges to $K_\infty$ in the Hausdorff sense. If $K_m$ is a weak barrier of gaussian curvature at least $k_m$ inside $U_m$ for all $m$, then $K_\infty$ is a weak barrier of gaussian curvature at least $k$ inside $U$.
\endproclaim
\proclabel{WeakBarriersClosedUnderLimits}
\proof Upon extracting a subsequence, we may suppose that for all $m$, $d_H(K_m,K_\infty)\leq 1/m$ and that $k_m\geq k-1/m$. For all $m$, let $(\epsilon_{m,p})_{p\in\Bbb{N}}>0$ be a sequence converging to $0$, let $(V_{m,p})_{p\in\Bbb{N}}$ be an increasing sequence of open subsets of $U_m$ such that $U_m=\munion_{p\in\Bbb{N}}V_{m,p}$ and let $(K_{m,p})_{p\in\Bbb{N}}$ be a sequence of convex sets converging to $K_m$ in the Hausdorff sense such that, for all $m$, $K_{m,p}$ is a strong barrier of gaussian curvature at least $k_m-\epsilon_{m,p}$ inside $V_{m,p}$. Upon extracting subsequences, we may suppose, in addition, that for all $m$ and for all $p$, $\epsilon_{m,p}\leq 1/p$ and $d_H(K_{m,p},K_m)\leq 1/p$. Let $(V_m)_{\minn}$ be an increasing sequence of relatively compact open subsets of $U$ such that $U=\munion_{m\in\Bbb{N}}V_m$. We may suppose that $V_p\subseteq V_{m,p}$ for all $m$ and for all $p$. For all $m$, define $K'_m:=K_{m,m}$. Then, for all $m$, $d_H(K'_m,K_\infty)\leq 2/m$ and $K'_m$ is a strong barrier of gaussian curvature at least $k_m-1/m\geq k-2/m$ inside $V_m$. In particular, $(K_m')_\minn$ converges to $K_\infty$ in the Hausdorff sense and we conclude that $K_\infty$ is a weak barrier of gaussian curvature at least $k$ over $U$ as desired.\qed
\medskip
\noindent We now show that the set of weak barriers in closed under intersection.
\proclaim{Lemma \nextprocno}
\noindent Let $U\subseteq\Bbb{R}^{n+1}$ be an open set. Let $k>0$ be a positive real number. Let $K_1$ and $K_2$ be compact, convex subsets of $\Bbb{R}^{n+1}$. If $K_1$ and $K_2$ are both weak barriers of gaussian curvature at least $k$ inside $U$, and if $K_1\minter K_2$ has non-trivial interior, then $K_1\minter K_2$ is also a weak barrier of gaussian curvature at least $k$ inside $U$.
\endproclaim
\proclabel{WeakBarriersClosedUnderIntersection}
\proof By definition, for each $i\in\left\{1,2\right\}$, there exists an increasing sequence $(V_{i,m})_\minn$ of open subsets of $U$, a sequence $(\epsilon_{i,m})_\minn$ of positive real numbers converging to $0$, and a sequence $(K_{i,m})_\minn$ of compact, convex subsets of $\Bbb{R}^{n+1}$ converging to $K_i$ in the Hausdorff sense with the properties that $U=\munion_{\minn}V_{i,m}$ and, for all $m$ and for all $x\in(\partial K_{i,m})\minter V_{i,m}$, $(\partial K_{i,m})$ is smooth near $x$ and has gaussian curvature at least $k-\epsilon_{i,m}$ at $x$. Let $V_m$ be an increasing sequence of relatively compact, open subsets of $U$ such that $U=\munion_{\minn} V_m$. Upon extracting subsequences, we may suppose that for all $m$, $V_m\subseteq V_{1,m},V_{2,m}$. For all $m$, denote $\epsilon_m:=\opMax(\epsilon_{1,m},\epsilon_{2,m})$, so that $(\epsilon_m)_\minn$ also converges to $0$, and that, for all $m$, for each $i$, and for all $x\in(\partial K_{i,m})\minter V_m$, $(\partial K_{i,m})$ is smooth with gaussian curvature at least $k-\epsilon_m$ at $x$.
\par
Let $(W_m)_\minn$ be a sequence of relatively compact open subsets of $U$ such that $U=\munion_{\minn}W_m$ and that, for all $m$, $\overline{W}_m\subseteq V_m$. For all $m$, denote $K_m:=K_{1,m}\minter K_{2,m}$. Observe that $(K_m)_\minn$ converges to $K_1\minter K_2$ in the Hausdorff sense, and we may therefore suppose that $d_H(K_m,K_1\minter K_2)\leq 1/m$ for all $m$. Choose $m\in\Bbb{N}$. Denote by $d_m$ the distance in $\Bbb{R}^{n+1}$ to $K_m$. By Theorem \procref{SmoothIntersectionHasHighCurvature}, there exists $r<1/2m$ and $\sigma<r$ such that if $s<\sigma$, if $x\in\overline{W}_m$ and if $d_{m,s}(x)=r$, then $d_{m,s}^{-1}(\left\{r\right\})$ is smooth near $x$ and has gaussian curvature at least $k-2\epsilon_m$ at $x$. In particular, if we denote $K_m':=d_{m,s}^{-1}(]-\infty,r])$ then for all $m$, $K_m'$ is a strong barrier of gaussian curvature at least $k-2\epsilon_m$ in $W_m$.
\par
By Corollary \procref{MollifierOfDistanceFunction}, for all $s<\sigma<r$, $\|d_{m,s}-d_m\|_0\leq r$, and so
$$
K_m = d_m^{-1}(]-\infty,0]) \subseteq d_{m,s}^{-1}(]-\infty,r]) = K_m',
$$
and
$$
K_m' = d_{m,s}^{-1}(]-\infty,r]) \subseteq d_m^{-1}(]-\infty,2r[) = \overline{B}_{2r}(K_m),
$$
so that $d_H(K_m,K_m')\leq 2r < 1/m$. It follows that $d_H(K_m',K_1\minter K_2)< 2/m$, so that $(K_m')_\minn$ converges to $K_1\minter K_2$ in the Hausdorff sense, and $K_1\minter K_2$ is therefore a weak barrier of gaussian curvature at least $k$ in $U$, as desired.\qed
\medskip
We refine Lemma \procref{WeakBarriersClosedUnderIntersection} in order to construct a local excision operation which allows us to obtain regularity for extremal weak barriers, as we shall see presently.
\proclaim{Lemma \nextprocno}
\noindent Let $K$ be a compact, convex subset of $\Bbb{R}^{n+1}$. Let $V$ be an open, convex subset of $\Bbb{R}^{n+1}$. Let $L$ be a compact, convex subset of $\overline{V}$. If $K\minter(\partial V)\subseteq L$, then $(K\setminus\overline{V})\munion(K\minter L)$ is compact and convex.
\endproclaim
\proof Denote $K':=(K\setminus\overline{V})\munion (K\minter L)$. Since $K\minter(\partial V)\subseteq L$, $K'=(K\setminus V)\munion(K\minter L)$, and since both $K\setminus V$ and $K\minter L$ are compact, so too is $K'$. Choose $x,x'\in K'$. For all $t\in[0,1]$, denote $x_t:=(1-t)x+tx'$. We claim that $x_t\in K'$ for all $t$. Indeed, since $x$ and $x'$ are both elements of $K$, by convexity, $x_t\in K$ for all $t$. Let $I$ be the set of all $t$ such that $x_t\in\overline{V}$. Observe that $I$ is a closed subinterval of $[0,1]$. Consider $t\in\partial I$. If $t\in]0,1[$, then $x_t$ is an element of $K\minter\partial V\subseteq L$. Otherwise, if $t\in\left\{0,1\right\}$, then $x_t\in K'\minter\overline{V}=K\minter L\subseteq L$. In each case, $x_t\in L$ for each $t\in\partial I$, and so, by convexity, $x_t\in L$ for all $t\in I$. That is, for all such $t$, $x_t\in K\minter L\subseteq K'$. However, for all $t\in[0,1]\setminus I$, $x_t\in K\setminus\overline{V}\subseteq K'$, so that $x_t\in K'$ for all $t$. Since $x,x'\in K'$ are arbitrary, we conclude that $K'$ is convex, as desired.\qed
\proclaim{Lemma \nextprocno}
\noindent Choose $k>0$. Let $K$ be a compact, convex subset of $\Bbb{R}^{n+1}$. Let $U$ be an open subset of $\Bbb{R}^{n+1}$ and suppose that $K$ is a strong barrier of gaussian curvature at least $k$ in $U$. Let $V$ be an open, convex subset of $\Bbb{R}^{n+1}$ whose closure is contained in $U$, and let $L$ be a compact, convex subset of $\overline{V}$. If $L$ is a strong barrier of gaussian curvature at least $k$ in $V$ and if $K\minter(\partial V)$ is contained in the relative interior of $L$ in $\overline{V}$, then $(K\setminus\overline{V})\munion(K\minter L)$ is a weak barrier of gaussian curvature at least $k$ in $U$.
\endproclaim
\proclabel{WeakBarriersClosedUnderExcision}
\proof Denote $K':=(K\setminus\overline{V})\munion(K\minter L)$. Let $d'$, $d_K$ and $d_L$ be the respective distances in $\Bbb{R}^{n+1}$ to $K'$, $K$ and $L$. Likewise, let $\Pi'$, $\Pi_K$ and $\Pi_L$ be their respective closest point projections. Since $K\setminus V$ is compact, and since $K\minter(\partial V)$ is contained in the relative interior of $L$ in $\overline{V}$, there exists $\delta>0$ such that for all $x\in K\setminus V$, $(\overline{B}_{\delta}(x)\minter \overline{V})\subseteq L$. We claim that for all $x\in K\setminus V$,
$$
K\minter\overline{B}_\delta(x) = K'\minter\overline{B}_\delta(x).\eqno{\nexteqnno}
$$
Indeed, for all $x\in K\setminus V$,
\eqnlabel{EqnKLooksTheSameInSmallBalls}
$$
(K\minter\overline{V})\minter\overline{B}_\delta(x)
\subseteq K\minter L\minter\overline{V}\minter\overline{B}_\delta(x)
= (K'\minter\overline{V})\minter\overline{B}_\delta(x),
$$
However, for all $x$,
$$
(K\setminus\overline{V})\minter\overline{B}_\delta(x) = (K'\setminus\overline{V})\minter\overline{B}_\delta(x),
$$
so that, for all $x\in K\setminus V$,
$$
K\minter\overline{B}_\delta(x) \subseteq K'\minter\overline{B}_\delta(x).
$$
On the other hand, since $K'\subseteq K$, $K'\minter\overline{B}_\delta(x)\subseteq K\minter\overline{B}_\delta(x)$, and the two sets therefore coincide, as desired.
\par
Define
$$
X := K'\setminus \munion_{x\in K\setminus V}B_\delta(x),
$$
and observe that $X$ is a compact subset of $V$. Choose $\rho_1>0$ such that $X_1:=\overline{B}_{2\rho_1}(X)\subseteq V$. We now claim that for all $x\in\Bbb{R}^{n+1}\setminus\overline{B}_{\rho_1}(X)$ such that $d'(x)<\rho_1$,
$$
d'(x)=d_K(x).
$$
Indeed, for such an $x$, denote $y:=\Pi'(x)$ and $\msf{N}:=(x-y)/\|x-y\|$. In particular, $y$ is the closest point in $K'$ to $x$ and by Lemma \procref{CharacterisationOfClosestPoints}, $\msf{N}$ is a supporting normal to $K'$ at this point. However, since $y\notin X$, there exists $y'\in K\setminus V=K'\setminus V$ such that $y\in B_\delta(y')$. In particular, $\msf{N}$ is a supporting normal to $\overline{B}_\delta(y')\minter K'$ at $y$. However, by \eqnref{EqnKLooksTheSameInSmallBalls}, $\overline{B}_\delta(y')\minter K'=\overline{B}_\delta(y')\minter K$, so that, by Lemma \procref{SupportingNormalIsLocalProperty}, $\msf{N}$ is a supporting normal to $K$ at $y$. In particular, by Lemma \procref{CharacterisationOfClosestPoints} again, $y$ is also the closest point in $K$ to $x$, so that
$$
d_K(x)=\|y-z\|=d'(x),
$$
as asserted.
\par
We now claim that for all $r<\rho_1$, there exists $\sigma_1:=\sigma_1(r)<r$ such that if $s<\sigma_1$, if $x\notin X_1$ and if $d'_s(x)=r$, then for all $y$ near $x$,
$$
d'_s(y) = d_{K,s}(y).
$$
Indeed, choose $r<\rho_1$. Fix $\eta>0$ such that $3\eta<\opMin(r,\rho_1-r)$ and fix $\sigma_1<\eta$. Choose $s<\sigma_1$ and $y\notin X_1$ such that $d'_s(y)\in[r-\eta,r+\eta]$. By Corollary \procref{MollifierOfDistanceFunction}, $d'(y)\in[r-2\eta,r+2\eta]$. Thus, if $z\in B_{\sigma_1}(y)$, then $d'(z)\in]0,\rho_1[$ and $z\notin\overline{B}_{\rho_1}(X)$, so that, by the discussion of the preceeding paragraph, $d'(z)-d_K(z)=0$, and it follows that $d'_s(y)=d_{K,s}(y)$, as desired.
\par
Let $(W_m)_\minn$ be an increasing family of relatively compact open subsets of $U$ with closure contained in $U$ such that $U=\munion_\minn W_m$. Suppose furthermore that $X_1\subseteq W_m$ for all $m$. Fix $m\in\Bbb{N}$. Choose $R>0$ such that $K\subseteq\overline{B}_R(0)$. By Theorem \procref{SmoothIntersectionHasHighCurvature} (with $K_2=\overline{B}_R(0)$), there exists $\rho_2<\opMin(\rho_1,\frac{1}{2n})$ with the property that for all $r<\rho_2$, there exists $\sigma_2:=\sigma_2(r)<\sigma_1(r)$ such that if $s<\sigma_2$, if $x\in\overline{W}_m$ and if $d_{K,s}(x)=r$, then $(d_{K,s})^{-1}(\left\{r\right\})$ is smooth near $x$ and has gaussian curvature at least $k-1/m$ at $x$. In particular, by the discussion of the preceeding paragraph that if $r<\rho_2$, if $s<\sigma_2$ and if $x\in\overline{W}_m\setminus X_1$ is such that $d'_s(x)=r$, then $(d_s')^{-1}(\left\{r\right\})$ is also smooth near $x$ and has gaussian curvature at least $k-1/m$ at $x$. On the other hand, by Theorem \procref{SmoothIntersectionHasHighCurvature} again, there exists $\rho_3<\rho_2$ with the property that for all $r<\rho_3$, there exists $\sigma_3:=\sigma_3(r)<\sigma_2(r)$ such that if $s<\sigma_3$, if $x\in X_1$ and if $d'_s(x)=r$, then $(d'_s)^{-1}(\left\{r\right\})$ is smooth near $x$ and has gaussian curvature at least $k-1/m$ at $x$. It follows that if $K_m:=(d_s')^{-1}(]0,r])$, then $K_m$ is a strong barrier of gaussian curvature at least $k-1/m$ in $W_m$. Furthermore, by Corollary \procref{MollifierOfDistanceFunction}, $\|d'-d'_s\|_0<\sigma_3<r$ so that
$$
K' = (d')^{-1}(]-\infty,0]) \subseteq (d_s')^{-1}(]-\infty,r])=K_m,
$$
and
$$
K_m = (d_s')^{-1}(]-\infty,r]) \subseteq d^{-1}(]-\infty,2r]) =\overline{B}_{2r}(K').
$$
Since $r<1/2m$, it follows that $d_H(K',K_m)\leq 1/m$, and since $m$ is arbitrary, we conclude that $(K_m)_\minn$ converges to $K'$ in the Hausdorff sense, so that $K'$ is a weak barrier of gaussian curvature at least $k$ over $U$, as desired.\qed
\newsubhead{The Plateau problem}
The machinery developed in the preceeding sections allows us to solve a general version of the Plateau Problem. Let $K$ be a compact, convex subset of $\Bbb{R}^{n+1}$ with smooth boundary and non-trivial interior. Let $X$ be a closed subset of the boundary of $K$ such that $\opConv(X)$ also has non-trivial interior. Choose $k>0$, and suppose that $\partial K$ has gaussian curvature at least $k$ at every point of $(\partial K)\setminus X$. Observe that, using the terminology of the preceeding section, this means that $K$ is a strong barrier of gaussian curvature at least $k$ in $\Bbb{R}^{n+1}\setminus X$. We define the family $\Cal{B}(k,K,X)$ to be the set of all compact, convex subsets $K'$ of $\Bbb{R}^{n+1}$ such that $X\subseteq K'\subseteq K$ and $K'$ is a weak barrier of gaussian curvature at least $k$ in $\Bbb{R}^{n+1}\setminus X$. Since strong barriers are also weak barriers, we see that $K$ itself is an element of $\Cal{B}(k,K,X)$ so that this family is non-empty.
\proclaim{Lemma \nextprocno}
\noindent If $L$ is an element of $\Cal{B}(k,K,X)$, then $L$ has non-trivial interior.
\endproclaim
\proclabel{EveryElementHasNonTrivialInterior}
\proof By definition, $L$ is compact and convex. Since $X\subseteq L$, using Lemma \procref{ConvexHullOfConvexIsTheSame}, we have $\opConv(X)\subseteq\opConv(L)=L$. Since $\opConv(X)$ has non-trivial interior, it follows that $L$ too has non-trivial interior, as desired.\qed
\proclaim{Lemma \nextprocno}
\noindent Let $L$ be an element of $\Cal{B}(k,K,X)$. If $\Sigma$ is a smooth embedded hypersurface (without boundary) such that $\Sigma\subseteq L$, then $\Sigma$ has gaussian curvature at least $k$ at every point of $(\Sigma\minter\partial L)\setminus X$.
\endproclaim
\proclabel{SmoothPointsHaveCurvatureBoundedBelow}
{\sl{\bf\noindent Remark:\ }In other words, every element of $\Cal{B}(k,K,X)$ is a viscosity supersolution of the gauss curvature equation (c.f. \cite{CrandhallIshiiLyons}.}
\medskip
\proof Consider $x\in(\Sigma\minter\partial L)\setminus X$. Without loss of generality, we may suppose that $x=0$. Since every supporting normal to $L$ at $0$ is also normal to $\Sigma$, $L$ has only one supporting normal at this point, which we may take to be $-e_{n+1}$. By Theorem \procref{ExistenceOfGraphFunction}, there exist $C,\rho>0$ and a convex, $C$-Lipschitz function $\omega:B_\rho'(0)\rightarrow]-C\rho,C\rho[$ such that $\partial L\minter (B_\rho'(0)\times]-2C\rho,2C\rho[)$ coincides with the graph of $\omega$. Upon reducing $\rho$ if necessary, we may suppose furthermore that there exists a smooth function $f:B_\rho'(0)\rightarrow]-C\rho,C\rho[$ such that $\Sigma\minter(B_\rho'(0)\times]-2C\rho,2C\rho[)$ coincides with the graph of $f$. In particular, since $\Sigma\subseteq L$, $f\geq\omega$.
\par
Fix $r<\rho$. For $0<t<C/2r$ and for $\left|s\right|<t$, denote $f_{s,t}(x'):=f(x') + t\|x'\|^2 + sr^2$ and let $\Sigma_{s,t}$ be the graph of $f_{s,t}$ over $\overline{B}'_r(0)$. Observe that for all $\left|s\right|<t$, $\partial\Sigma_{s,t}$ lies in the interior of $L$, and for all $s>0$, the whole of $\Sigma_{s,t}$ lies in the interior of $L$.
\par
Fix $0<t<C/2r$. Let $(\epsilon_m)_{\minn}$ be a sequence of positive numbers converging to $0$, let $(V_m)_{m\in\Bbb{N}}$ be an increasing sequence of open subsets of $\Bbb{R}^{n+1}\setminus X$ and let $(L_m)_{m\in\Bbb{N}}$ be a sequence of convex sets converging to $L$ in the Hausdorff sense such that $\Bbb{R}^{n+1}\setminus X=\munion_{\minn} V_m$ and, for all $m$, $L_m$ is a strong barrier of gaussian curvature at least $k-\epsilon_m$ inside $V_m$. Fix $s_0<0<s_1$ such that $\left|s_0\right|,\left|s_1\right|<t$. For sufficiently large $m$, the whole of $\Sigma_{s_1,t}$ is contained in the interior of $L_m$, $\partial\Sigma_{s,t}$ is contained in the interior of $L_m$ for all $s\in[-s_0,s_1]$, but $\Sigma_{-s_0,t}$ intersects the complement of $L_m$ non-trivially. There therefore exists $s\in]-s_0,s_1[$ such that $\Sigma_{s,t}$ is an interior tangent to $\partial L_m$ at some point, $(x_m',f_{s,t}(x_m'))$, say. By the maximum principal, $\Sigma_{s,t}$ has gaussian curvature at least $k-\epsilon_m$ at this point. By compactness, letting $t$ tend to $0$, we conclude that there exists $x'\in\overline{B}'_r(0)$ such that $\Sigma$ has gaussian curvature at least $k$ at $x'$, and since $r>0$ is arbitrary, we conclude that $\Sigma$ has gaussian curvature at least $k$ at $0$, as desired.\qed
\medskip
For any Borel measurable subset $X$ of $\Bbb{R}^{n+1}$, we define the \gloss{volume} of $X$ to be its $(n+1)$-dimensional Lebesgue measure, and we denote it by $\opVol(X)$. We denote
$$
V_0 := \minf_{L\in\Cal{B}(k,K,X)}\opVol(L).
$$
\proclaim{Lemma \nextprocno}
\noindent $V_0>0$.
\endproclaim
\proclabel{InfimumOfVolumeIsPositive}
\proof Choose $L\in\Cal{B}(k,K,X)$. By definition, $L$ is compact and convex. Since $X\subseteq L$, and bearing in mind Lemma \procref{ConvexHullOfConvexIsTheSame}, $\opConv(X)\subseteq\opConv(L)=L$. Thus, by monotonicity of Lebesgue measure, $\opVol(L)\geq\opVol(\opConv(X))$. However, since $\opConv(X)$ has non-trivial interior, $\opVol(\opConv(X))>0$, and so
$$
V_0 = \minf_{L\in\Cal{B}(k,K,X)}\opVol(L) \geq \opVol(\opConv(X)) > 0,
$$
\noindent as desired.\qed
\proclaim{Lemma \nextprocno}
\noindent Let $K_0\subseteq K_1$ be compact, convex subsets of $\Bbb{R}^{n+1}$ with non-trivial interiors. If $K_0\neq K_1$, then $\opVol(K_0)<\opVol(K_1)$.
\endproclaim
\proclabel{ContainedMeansLessVolume}
\proof Choose $x\in K_1\setminus K_0$. Since $K_0$ is compact, there exists $\delta_1>0$ such that $B_{\delta_1}(x)\minter K_0=\emptyset$. Let $y$ be an interior point of $K_0$. There exists $\delta_2>0$ such that $B_{\delta_2}(y)\subseteq K_0$. By convexity, for all $t\in]0,1]$, and for all $z\in B_{t\delta_2}(0)$,
$$
(1-t)x + ty + z = (1-t)x + t(y + z/t) \in K_1.
$$
That is, for all $t\in]0,1]$, $B_{t\delta_2}((1-t)x+ty)\subseteq K_1$. Choose $t>0$ such that $t(\|y-x\| + \delta_2)<\delta_1$. In particular $B_{t\delta_2}((1-t)x+ty)\minter K_0=\emptyset$, and so, by additivity and monotonicity of Lebesgue measure,
$$
\opVol(K_1) \geq \opVol(K_0) + \opVol(B_{t\delta_2}((1-t)x+ty)) > \opVol(K_0),
$$
as desired.\qed
\proclaim{Theorem \nextprocno}
\noindent There exists a unique element $K_0\in\Cal{B}(k,K,X)$ such that
$$
\opVol(K_0) = V_0.
$$
\endproclaim
\proclabel{VolumeMinimiserExists}
\proof We first show uniqueness. Indeed, suppose that there exists $K_0\neq K_0'\in\Cal{B}(k,K,X)$ such that, $\opVol(K_0)=\opVol(K_0')=V_0$. Since $K_0\neq K_0'$, without loss of generality, we may assume that $K_0\minter K_0'\neq K_0$. Since $X$ is contained in each of $K_0$ and $K_0'$, it is also contained in $K_0\minter K_0'$. Moreover, since both $K_0$ and $K_0'$ are contained in $K$, so too is $K_0\minter K_0'$. Finally, by Lemma \procref{WeakBarriersClosedUnderIntersection}, $K_0\minter K_0'$ is a weak barrier of gaussian curvature at least $k$ over $\Bbb{R}^{n+1}\setminus X$ and we conclude that $K_0\minter K_0'$ is an element of $\Cal{B}(k,K,X)$. However, by Lemma \procref{ContainedMeansLessVolume}, $\opVol(K_0\minter K_0')<\opVol(K_0)$. This contradicts minimality of $K_0$, and uniqueness follows.
\par
Let $(L_m)_\minn\in\Cal{B}(k,K,X)$ be a sequence such that $(\opVol(L_m))_\minn$ converges to $V_0$. For all $m$, define $K_m := L_1\minter...\minter L_m$. For all $m$, $X\subseteq K_m\subseteq K$, and, by Lemma \procref{WeakBarriersClosedUnderIntersection}, $K_m$ is also a weak barrier of gaussian curvature at least $k$ over $\Bbb{R}^{n+1}\setminus X$, so that $K_m\in\Cal{B}(k,K,X)$. Moreover, by monotonicity of the Lebesgue measure, for all $m$, $V_0\leq\opVol(K_m)\leq \opVol(L_m)$. In particular, $(\opVol(K_m))_\minn$ also converges to $V_0$. Define
$$
K_\infty := \minter_{m\in\Bbb{N}}K_m.
$$
Since $\opVol(K_\infty)\leq\opVol(K_m)$ for all $m$, we have $\opVol(K_\infty)\leq V_0$. However $X\subseteq K_\infty\subseteq K$, and, by Lemma \procref{NestedIntersectionIsHausdorffLimit}, $(K_m)_\minn$ converges to $K_\infty$ in the Hausdorff sense. It follows by Lemma \procref{WeakBarriersClosedUnderLimits} that $K_\infty$ is a weak barrier of gaussian curvature at least $k$ over $\Bbb{R}^{n+1}\setminus X$. That is, $K_\infty\in\Cal{B}(k,K,X)$ and so $V_0\leq\opVol(K_\infty)$. We conclude that $\opVol(K_\infty)=V_0$, and this completes the proof.\qed
\newsubhead{Singularities and smoothness}
Continuing to use the notation of the preceeding section, we now show that the volume minimiser solves the Plateau problem modulo singularities of a type that are now well understood.
\subheadlabel{SingularitiesAndSmoothness}
\proclaim{Theorem \nextprocno}
\noindent Let $K_0\in\Cal{B}(k,K,X)$ be the volume minimiser. Then $(\partial K_0)\minter K^o$ has gaussian curvature equal to $k$ in the viscosity sense. Furthermore, if $x\in(\partial K_0)\setminus X$, then either
\medskip
\myitem{(1)} $(\partial K_0)$ is smooth near $x$ and has gaussian curvature equal to $k$ at $x$; or
\medskip
\myitem{(2)} $K_0$ satisfies the local geodesic property at $x$.
\endproclaim
\proclabel{VolumeMinimiserIsSmoothWhenNotLGP}
\proof By Lemma \procref{SmoothPointsHaveCurvatureBoundedBelow}, $(\partial K_0)\minter K^o$ has gaussian curvature at least $k$ in the viscosity sense. Now choose $x\in(\partial K_0)\setminus X$. Suppose that $K_0$ satisfies the local geodesic property at $x$. Then, if $\Sigma$ is a smooth, embedded surface (without boundary) contained in $\overline{K_o^c}$ and if $x\in\Sigma$, then, provided $\Sigma$ is oriented such that its normal points outwards from $K_o$, this surface has non-positive curvature at $x$. In particular, $(\partial K_0)\minter K^o$, has curvature at most $0\leq k$ in the viscosity sense at $x$.
\par
Now suppose that $K_0$ does not satisfy the local geodesic property at $x$. Let $U$ be a relatively compact neighbourhood of $x$ whose closure is contained in $\Bbb{R}^{n+1}\setminus X$. Let $(K_m)_\minn$ be a sequence of compact, convex subsets of $\Bbb{R}^{n+1}$ with the properties that $(K_m)_\minn$ converges to $K_0$ in the Hausdorff sense and, for all $m$, $K_m$ is a strong barrier of gaussian curvature at least $k-1/m$ in $U$. Denote $K_\infty:=K_0$ and $x_\infty=x_0$.
\par
Let $(x_m)_\minn$ be a sequence converging to $x_\infty$ such that $x_m\in\partial K_m$ for all $m$. Upon applying a convergent sequence of affine isometries, we may suppose that $x_m=0$ for all $m$. Since $K_\infty$ has non-trivial interior, by Lemma \procref{NonTrivialInteriorMeansSupportingNormalsContainedInHypersphere}, $\Cal{N}(0;K_\infty)$ is strictly contained in a hemisphere. By Lemma \procref{OptimalSupportingNormal}, there exists $\msf{N}\in\Cal{N}(0;K_\infty)$ such that $\langle\msf{N},\msf{M}\rangle>0$ for all $\msf{M}\in\Cal{N}(0;K_\infty)$. By compactness of $\Cal{N}(0;K_\infty)$, there exists $\theta\in[0,\pi/2[$ such that $\langle\msf{N},\msf{M}\rangle>3\opCos(\theta)$ for all $\msf{M}\in\Cal{N}(0;K_\infty)$. Denote $C:=\opTan(\theta)$.
\par
By Lemma \procref{ConvergenceOfSupportingNormals}, upon extracting a subsequence, there exists $r>0$ such that $\overline{B}_r(0)\subseteq U$ and, for all $m$, for all $x\in(\partial K_m)\minter B_r(0)$ and for all $\msf{M}\in\Cal{N}(x;K_m)$, $\langle\msf{N},\msf{M}\rangle>2\opCos(\theta)$. We denote $\rho=r/\sqrt{1+4C^2}$. By Lemma \procref{NoLGPMeansPointyBit}, there exists $\msf{N}'$, which we may choose as close to $\msf{N}$ as we wish such that for all $x\in K_\infty\setminus B_{\rho/2}(0)$, $\langle x,\msf{N}'\rangle<0$. Moreover, we may assume that for all $m$, for all $x\in(\partial K_m)\minter B_r(0)$ and for all $\msf{M}\in\Cal{N}(x;K_m)$, $\langle \msf{N}',\msf{M}\rangle > \opCos(\theta)$.
\par
Upon applying a rotation, we may suppose that $\msf{N}'=-e_{n+1}$. By Theorem \procref{ExistenceOfGraphFunction}, for all $m$, there exists a convex, $C$-Lipschitz function $\hat{f}_m:B_\rho'(0)\rightarrow]-C\rho,C\rho[$ such that $\hat{f}_m(0)=0$ and $(\partial K_m)\minter (B_\rho'(0)\times]-2C\rho,2C\rho[)$ coincides with the graph of $\hat{f}_m$ over $B_\rho'(0)$. By the Arzela-Ascoli theorem, every subsequence of $(\hat{f}_m)_\minn$ has a subsubsequence converging in the local uniform sense over $B_\rho'(0)$ to some limit $\hat{f}_\infty'$ say. Furthermore, since $(K_m)_\minn$ converges to $K_\infty$ in the Hausdorff sense, we conclude that $\hat{f}_\infty'=\hat{f}_\infty$. It follows that $(\hat{f}_m)_\minn$ converges in the local uniform sense over $B_\rho'(0)$ to $\hat{f}_\infty$.
\par
By construction, $\hat{f}_\infty(x')>2\delta>0$ for all $x'\in\partial B_{\rho/2}'(0)$ and for some $\delta>0$. Since $(\hat{f}_m)_\minn$ converges locally uniformly to $\hat{f}_\infty$ over $B_\rho'(0)$, we may suppose that $\hat{f}_m(x')>\delta$ for all $m$ and for all $x'\partial B'_{\rho/2}(0)$.
\par
Choose $m<\infty$. Observe that $\hat{f}_m$ is smooth and strictly convex. Denote $\overline{\Omega}_m:=\hat{f}_m^{-1}(]-\infty,\delta])$ and observe that $\overline{\Omega}_m$ is a compact, convex subset of $B_{\rho/2}'(0)$. By strict convexity, $D\hat{f}_m$ only vanishes at the unique absolute minimum of $\hat{f}_m$ over $B'_{\rho/2}(0)$. However, since $\hat{f}_m(0)=0$, this absolute minimum is contained in the interior of $\overline{\Omega}_m$. In particular, $D\hat{f}_m$ does not vanish at any boundary point of $\overline{\Omega}_m$, so that $\overline{\Omega}_m$ has smooth boundary. Thus, by Theorem \procref{FirstExistenceTheorem}, there exists a unique, smooth, strictly convex function $f_m:\overline{\Omega}_m\rightarrow]-\infty,\delta]$ such that $f_m(x')=\delta$ for all $x'\in\partial\Omega_m$ and the graph of $f_m$ has constant gaussian curvature equal to $k$. By convexity, $f_m\leq \delta$, and by Lemma \procref{StrongMaximumPrinciple}, $f_m\geq\hat{f}_m$.
\par
Define $V_m:=\Omega_m\times]-2C\rho,2C\rho[$. Observe that $V_m$ is open and convex. Moreover, $\overline{V}_m\subseteq\overline{B}_r(0)\subseteq U$. We define the subset $L_m$ of $\overline{V}_m$ by
$$
L_m := \left\{ (x',t)\ |\ x'\in\overline{\Omega}_m,\ f_m(x')\leq t\leq 2C\rho\right\}.
$$
Observe that $L_m$ is compact and convex, $K_m\minter\partial V_m\subseteq L_m$ and $L_m\subseteq K_m\minter\overline{V}_m$. Define $K'_m:=(K_m\setminus\overline{V})\munion(K_m\minter L_m)=(K_m\setminus\overline{V})\munion L_m$. We claim that $K'_m$ is a weak barrier of gaussian curvature at least $k$ in $\Bbb{R}^{n+1}\setminus X$. Indeed, for all $s\in[0,C\rho[$, define the subset $L_{m,s}$ of $\overline{V}_m$ by
$$
L_{m,s} := \left\{ (x',t)\ |\ x'\in\overline{\Omega}_m,\ f_m(x')-s\leq t\leq 2C\rho\right\},
$$
and define $K'_{m,s}:=(K_m\setminus\overline{V})\munion(K_m\minter L_{m,s})$. For all $s>0$, $K_m\minter\partial V_m$ is contained in the relative interior of $L_m$ in $V_m$. Thus, by Lemma \procref{WeakBarriersClosedUnderExcision}, $K'_{m,s}$ is a weak barrier of gaussian curvature at least $k$ over $\Bbb{R}^{n+1}\setminus X$. Thus, since $(K'_{m,s})_{s\in[0,(1/2)C\rho[}$ converges to $K'_{m,0}=K'_m$ in the Hausdorff sense as $s$ tends to $0$, it follows by Lemma \procref{WeakBarriersClosedUnderLimits} that $K'_m$ is also a weak barrier of gaussian curvature at least $k$ in $\Bbb{R}^{n+1}$, as asserted.
\par
By Lemmas \procref{CompactnessOfCompactSets} and \procref{ConvexSetsIsAClosedSet}, we may suppose that $(K'_m)_\minn$ converges towards a compact, convex subset, $K'_\infty$, say, of $\Bbb{R}^{n+1}$. We claim that $K_\infty'=K_\infty$. Indeed, by Lemma \procref{WeakBarriersClosedUnderLimits}, $K'_\infty$ is a weak barrier of gaussian curvature at least $k$ in $\Bbb{R}^{n+1}\setminus X$. Moreover, for all $m$, $X\subseteq K_m\setminus B_r(0)\subseteq K_m\setminus V_m\subseteq K_m'$ and so $X\subseteq K'_\infty$. Finally, for all $m$, $K'_m\subseteq K_m\subseteq K$, so that $K'_\infty\subseteq K_\infty\subseteq K$. We conclude that $K_\infty'$ is an element of $\Cal{B}(k,K,X)$. In particular, $V_0\leq\opVol(K'_\infty)$. However, since $K'_\infty\subseteq K_\infty$, $\opVol(K'_\infty)\leq\opVol(K_\infty)\leq V_0$, so that $\opVol(K'_\infty)=V_0$. It follows by uniqueness that $K'_\infty=K_\infty$, as asserted.
\par
By continuity, there exists $\rho'<\rho$ such that for all $x'\in \overline{B}_{\rho'}'(0)$, $\hat{f}_\infty(x')<\delta/2$. Since $(\hat{f}_m)_\minn$ converges to $\hat{f}_\infty$ uniformly over $\overline{B}_{\rho'}'(0)$, we may suppose that for all $m$ and for all $x'\in\overline{B}_{\rho'}'(0)$, $\hat{f}_m(x')\leq\delta$. In particular, for all $m$, $\overline{B}_{\rho'}'(0)\subseteq\overline{\Omega}_m$. We therefore define $W:=B_{\rho'}'(0)\times]-(3/2)C\rho,(3/2)C\rho[$, and, for all $m<\infty$, $(\partial K_m)\minter W=(\partial L_m)\minter W$ is smooth with constant gaussian curvature equal to $k$. However, since $K_\infty$ does not satisfy the local geodesic property at $x$, it follows by Theorem \procref{StructureOfSingularities}, that $(\partial K_\infty)\minter W$ is smooth with constant gaussian curvature equal to $k$, and this completes the proof.\qed
\medskip
The boundary of the volume minimiser, $K_0$, therefore solves the Plateau problem in the very general setting where $X$ is any closed subset of $\partial K$. We say that a point $x\in\partial K_0$ is \gloss{regular} if $\partial K_0$ is smooth near that point. We define the \gloss{singular set}, $\opSing(K_0)$, to be the set of all points of $\partial K_0$ that are not regular. We obtain the following characterisation.
\proclaim{Theorem \nextprocno}
\noindent There exists a family $(X_\alpha)_{\alpha\in A}$ of subsets of $X$ such that
$$
\opSing(K_0) = \munion_{\alpha\in A}\opConv(X_\alpha).
$$
\endproclaim
\proclabel{ThmDescriptionOfSingularSet}
\proof By definition, $\opSing(K_0)$ is closed. Furthermore, by Theorem \procref{VolumeMinimiserIsSmoothWhenNotLGP}, $\opSing(K_0)\setminus X$ consists of all points of $\partial K_0\setminus X$ satisfying the local geodesic property, so that, by Theorem \procref{GBPBdyPointsAreConvexHull}, $\opSing(K_0)\subseteq\opConv(X)$.
\par
Now choose $x\in\opSing(K_0)$. Let $H$ be a supporting tangent hyperplane to $K_0$ at $x$. Since $\opConv(X)\subseteq K_0$, $H$ is also a supporting tangent hyperplane to $\opConv(X)$ at $x$. Denote $X_x:=X\minter H$. Since $H\minter\opConv(X)=\opConv(X_x)$, it follows that $x\in\opConv(X_x)$. Furthermore, since $\opConv(X_x)\subseteq K_0\minter H$, every point of $\opConv(X_x)$ is a boundary point of $K_0$. However, by Theorem \procref{ConvexHullHasLGP}, the set $\opConv(X_x)$ satisfies the local geodesic property at every point of $\opConv(X_x)\setminus X_x$. In particular, $K_0$ also satisfies the local geodesic property at every point of this subset, so that, by Theorem \procref{VolumeMinimiserIsSmoothWhenNotLGP}, $\opConv(X_x)\setminus X$ is contained in $\opSing(K_0)$. Since $\opSing(K_0)$ is closed, we conclude that the whole of $\opConv(X_x)$ is contained in $\opSing(K_0)$, and since $x\in\opSing(K_0)$ is arbitrary, we conclude that
$$
\opSing(K_0) = \munion_{x\in\opSing(K_0)}\opConv(X_x),
$$
as desired.\qed
\medskip
Various ad-hoc arguments can now be used to eliminate singularities. For example, by Lemma \procref{SmoothPointsHaveCurvatureBoundedBelow}, if the boundary of $\opConv(X)$ is smooth at some point, then that point must lie in the interior of $K_0$. In the particular case at hand, however, singularites are removed as follows.
\proclaim{Lemma \nextprocno}
\noindent Suppose that for every point $x$ of $\partial X$, there exists a $C^2$ function $f:\partial K\rightarrow\Bbb{R}$ such that $f(x)=0$, $Df(x)\neq 0$ and $f^{-1}(]-\infty,0])\subseteq X$. Then, $\opSing(K_0)\subseteq X$.
\endproclaim
\proclabel{InteriorBallCondition}
\proof Suppose the contrary. Let $x$ be a point of $\opSing(K_0)\setminus X$. By Theorem \procref{ThmDescriptionOfSingularSet}, there exists a subset $X'\subseteq X$ such that $x\in\opConv(X')\subseteq\partial K_0$. Furthermore, since $x\in X$, $X'$ contains at least two distinct points, $y_1$ and $y_2$, say, and, without loss of generality, $x$ lies along the straight line, $\Gamma$, passing through these two points. Let $\msf{N}$ be a supporting normal to $K_0$ at $x$. In particular, $\msf{N}$ is normal to $\Gamma$. For $\epsilon>0$, denote
$$
C_\epsilon := \munion_{z\in\Gamma}\overline{B}_\epsilon(x-\epsilon\msf{N}),
$$
so that, for all $\epsilon$, $C_\epsilon$ is the closed cylinder about the straight line, $\Gamma_\epsilon$, obtained by displacing $\Gamma$ a distance $\epsilon$ in the $-\msf{N}$ direction. We claim that for all sufficiently small $\epsilon$, $C_\epsilon\minter K\subseteq\opConv(X)$.
\par
It suffices to show that for sufficiently small $\epsilon$, $C_\epsilon\minter\partial K\subseteq X$ near $y_1$ and $y_2$. Without loss of generality, we may suppose that $y_1=0$, that $\msf{N}=e_n$ and that $x$ lies on the positive $x_{n+1}$ axis. Since $x$ is an interior point of $K$, for all $\msf{N}\in\Cal{N}(y_1;K)$, $\langle\msf{N},e_{n+1}\rangle=\|x-y\|^{-1}\langle x - y_1,\msf{N}\rangle<0$. By compactness of $\Cal{N}(y_1;K)$, we may suppose that there exists $\theta\in]0,\pi/2[$ such that $\langle\msf{N},e_{n+1}\rangle>2\opCos(\theta)$ for all $\msf{N}\in\Cal{N}(y_1)$. By Lemma \procref{ContinuityOfSupportingNormals}, there exists $r>0$ such that for all $y\in\partial K\minter B_r(y_1)$ and for all $\msf{N}\in\Cal{N}(y;K)$, $\langle\msf{N},e_{n+1}\rangle>\opCos(\theta)$. Denote $C:=\opTan(\theta)$ and $\rho=r/\sqrt{1+4C^2}$. By Theorem \procref{ExistenceOfGraphFunction}, there exists a convex, $C$-Lipschitz function $\omega:B'_\rho(0)\rightarrow]-C\rho,C\rho[$ such that $\partial K\minter (B'_\rho(0)\times]-2C\rho,2C\rho[)$ coincides with the graph of $\omega$. Furthermore, since $\partial K$ is smooth, so too is $\omega$.
\par
Let $f:\partial K\rightarrow\Bbb{R}$ be a $C^2$ function such that $f(y_1)=0$, $Df(y_1)\neq 0$ and $f^{-1}(]-\infty,0])\subseteq X$. Define $g:B'_\rho(0)\rightarrow\Bbb{R}$ by $g(x'):=f(x',\omega(x'))$. Observe that $g$ is $C^2$, $g(0)=0$ and $Dg(0)\neq 0$. Furthermore, if $g(z')\leq 0$, then $(z',\omega(z'))\in X$, and so, recalling that $\msf{N}$ is a supporting normal to $\opConv(X)$ at $y_1$, $\langle z',e_n\rangle = \langle (z',\omega(z')),\msf{N}\rangle \leq 0$. It follows that $Dg(0)=\lambda e_n$ for some $\lambda>0$. Thus, since $g$ is $C^2$, for sufficiently small $\epsilon>0$, $\overline{B}_\epsilon(-\epsilon e_n)\subseteq g^{-1}(]-\infty,0])$, so that
$$
C_\epsilon\minter\partial K\minter (B'_\rho(0)\times]-2C\rho,2C\rho[)
=\left\{ (z',\omega(z'))\ |\ z'\in\overline{B}_\epsilon(-\epsilon e_n)\right\} \subseteq X.
$$
That is, $C_\epsilon\minter\partial K\subseteq X$ near $y_1$. In like manner, we show that $C_\epsilon\minter\partial K\subseteq X$ also near $y_2$ so that, for sufficiently small $\epsilon$, $C_\epsilon\subseteq\opConv(X)$, as desired. However, $\partial C_\epsilon$ has zero curvature at every point. This is absurd, by Lemma \procref{SmoothPointsHaveCurvatureBoundedBelow}, and we conclude that $\opSing(K_0)$ is empty, as desired.\qed
\medskip
\noindent In particular, the classical existence result follows as an immediate corollary.
\proclaim{Theorem \procref{MostGeneralPlateau}}
\noindent Choose $k>0$. Let $K$ be a compact, convex subset of $\Bbb{R}^{n+1}$ with smooth boundary. Let $X$ be a closed subset of $\partial K$ with $C^2$ boundary $C=\partial X$. If $\partial K$ has gaussian curvature bounded below by $k$ at every point of $(\partial K)\setminus X$, then there exists a compact, strictly convex, $C^{0,1}$ embedded hypersurface $S\subseteq\Bbb{R}^{n+1}$ with the properties that
\medskip
\myitem{(1)} $S\subseteq K$;
\medskip
\myitem{(2)} $\partial S = C$; and
\medskip
\myitem{(3)} $S\setminus\partial S$ is smooth and has constant gaussian curvature equal to $k$.
\endproclaim
\line{\hfill\it Barcelona-Granada, May-June, 2012}
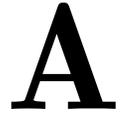
\global\headno=0
\inappendicestrue
\newhead{Terminology}
{\bf\noindent Derivatives:\ }For any vector spaces $E,F$, let $\opSymm(n,E)\otimes F$ denote the space of symmetric multilinear forms from $E$ into $F$. When $F=\Bbb{R}$, we denote simply $\opSymm(n,E)=\opSymm(n,E)\otimes\Bbb{R}$. For any open subset $U\subseteq E$ and for any $k$-times differentiable function $f:U\rightarrow F$, we denote the $k$'th total derivative by $D^kf:U\rightarrow\opSymm(k,E)\otimes F$. For any point $p\in U$ and for $k$ vectors $V_1,...,V_k\in E$, we denote $D^kf(p)(V_1,...,V_k)\in F$ the image of the $k$-tuplet $(V_1,...,V_k)$ under the action of $D^kf$ at the point $P$.
\headlabel{Terminology}
\medskip
\noindent For any vector spaces $E$ and $F$, for any open subset $U$ of $E$, and for all $k\in\Bbb{N}$, we denote by $C^k(U,F)$ the space of $k$-times continuously differentiable functions from $U$ into $F$. We denote by $C^\infty(U,F)$ the space of functions from $U$ into $F$ which have continuous derivatives of arbitrarily high order. When $F=\Bbb{R}$, we denote simply $C^k(U)=C^k(U,\Bbb{R})$ and $C^\infty(U)=C^\infty(U,\Bbb{R})$.
\medskip
\noindent For any $k\in\Bbb{N}$ and for any $f\in C^k(U)$, we define $J^k(f)\in C^0(U,\oplus_{k=0}^m\opSymm(n,E))$ by:
$$
J^k(f)(x) = (f(x),Df(x),...,D^kf(x)).
$$
\noindent We refer to $J^k(f)$ as the $k$-jet of $f$.
\medskip
{\bf\noindent Canonical Basis of Euclidean Space:\ }For all $n$, we denote by $\Bbb{R}^n$, $n$-dimensional, real space and by $e_1,...,e_n$ its canonical basis. We denote by $\langle\cdot,\cdot\rangle$ the Euclidean inner product and by $\|\cdot\|$ the Euclidean norm. For any open subset $U\subseteq\Bbb{R}^n$, for any $k$-times differentiable function $f:U\rightarrow\Bbb{R}$, and for any $k$-tuple of indices $1\leq i_1,...,i_k\leq n$, we define the function $(\partial_{i_1}...\partial_{i_k}f)$ such that for all $x\in U$:
$$
(\partial_{i_1}...\partial_{i_k}f)(x) = D^kf(x)(e_{i_1},...,e_{i_k}).
$$
\noindent We will also use the more concise notation:
$$
f_{i_1...i_k} := (\partial_{i_1}...\partial_{i_k}f).
$$
{\bf\noindent Distributional Derivatives:\ }Let $E$ be a vector space furnished with a volume form $\opdVol$. Let $U$ be an open subset of $E$ and let $f:U\rightarrow \Bbb{R}$ be a real valued function which is locally $L^1$. Let $g=(g_0,g_1,...,g_k):U\rightarrow \oplus_{k=0}^m\opSymm(n,E)$ be locally $L^1$. We say that $g$ is the $k$'th order distributional derivative of $f$ whenever it has the property that for any smooth function $\phi$ with compact support, for all $1\leq k\leq n$, and for all vectors $V_1,...,V_k$:
$$
\int_E f(x) (D^k\phi)(x)(V_1,...V_k)\opdVol = (-1)^k\int_E g_k(x)(V_1,...,V_k)\phi(x)\opdVol.
$$
{\bf\noindent Smooth Functions on Sets with Boundary:\ }$\Omega$ will always represent a bounded, strictly convex, open subset of $\Bbb{R}^n$. Given any vector space $E$, a function $f:\overline{\Omega}\rightarrow E$ is said to be $C^k$ whenever there exists an extension $\hat{f}$ of $f$ to $\Bbb{R}^n$ which is $k$-times continuously diferentiable. By Whitney's Extension Theorem (c.f. \cite{Simon}), the extension $\hat{f}$ can be chosen such that for all $k\leq l$:
$$
\|D^k\hat{f}\|_{L^\infty} = \|D^k\hat{f}|_{\overline{\Omega}}\|_{L^\infty} = \|D^kf\|_{L^\infty}.
$$
\noindent We say that $f$ is smooth whenever it is $C^k$ for all finite $k$. Given any open subset $U$ of $E$, we denote by $C^\infty(\overline{\Omega},U)$ the set of all smooth functions from $\overline{\Omega}$ into $E$ taking values in $U$. In particular, when $U=E=\Bbb{R}$, we denote $C^\infty(\overline{\Omega}) = C^\infty(\overline{\Omega},\Bbb{R})$.
\medskip
{\bf\noindent Non-linear Operators:\ }Given open subsets $U\subseteq\Bbb{R}^n$ and $V\subseteq\opSymm(2,\Bbb{R}^n)$ and a smooth function $F:\Bbb{R}\times\opSymm(1,\Bbb{R}^n)\times V\rightarrow \Bbb{R}$, for any function $f:U\rightarrow\Bbb{R}$ with the property that $D^2f(x)\in V$ for all $x\in U$, we define the function $F(f,Df,D^2f)$ such that, for all $x\in U$:
$$
F(f,Df,D^2f)(x) = F(f(x),Df(x),D^2f(x)).
$$
\noindent $F$ thus represents the most general second-order, non-linear partial differential operator acting on functions over $U$ which is homogeneous in the spatial variables.
\medskip
{\bf\noindent Decomposition of Euclidean Space:\ }We often decompose $\Bbb{R}^{n+1}$ as $\Bbb{R}^n\times\Bbb{R}$. For all $r>0$ and for all $x\in\Bbb{R}^{n+1}$, we denote by $B_r(x)$ the open ball of radius $r$ about $x$ in $\Bbb{R}^{n+1}$. For all $r>0$ and for all $x'\in\Bbb{R}^n$, we denote by $B'_r(x)$ the open ball of radius $r$ about $x$ in $\Bbb{R}^n$.
\medskip
{\bf\noindent Metrics:\ }Let $X$ and $Y$ be two compact subsets of $\Bbb{R}^{n+1}$. We define the Hausdorff distance between $X$ and $Y$ by:
$$
d_H(X,Y) = \msup_{x\in X}\minf_{y\in Y}\|x-y\| + \msup_{y\in Y}\minf_{x\in X}\|x-y\|.
$$
\noindent We denote by $\Sigma^n$ the sphere of unit radius in $\Bbb{R}^{n+1}$. We define the spherical distance $d_\Sigma:\Sigma^n\times\Sigma^n\rightarrow\Bbb{R}$ by:
$$
d_\Sigma(\msf{N},\msf{M}) = \opCos^{-1}\langle\msf{N},\msf{M}\rangle.
$$
\noindent The spherical distance thus measures the angle between two points in the sphere. Let $X$ and $Y$ be two compact subsets of $\Sigma^n$. We define the spherical-Hausdorff distance between $X$ and $Y$ by:
$$
d_{H,\Sigma}(X,Y) = \msup_{x\in X}\minf_{y\in Y}d_\Sigma(x,y) + \msup_{y\in Y}\minf_{x\in X}d_\Sigma(x,y).
$$
{\bf\noindent Miscellaneous:\ } If $X$ is any subset of $\Bbb{R}^n$, we denote its closure by $\overline{X}$, its interior by $X^o$ and its boundary by $\partial X$.
\medskip
\noindent Let $E$ be a vector space furnished with an inner product. For vectors $X$ and $Y$ in $E$, we denote by $\langle X,Y\rangle$ the inner product of $X$ with $Y$.
\medskip
\noindent Let $E$ be any vector space. For vectors $X_1,...,X_n$, we denote by $\langle X_1,...,X_n\rangle$ the linear subspace of $E$ spanned by $X_1,...,X_n$. This should not be confused with the inner product. It will be clear from the context which is meant. 
\newhead{Index}
\noindent\hbox{\vbox{%
\hsize = 6.5cm%
\raggedright%
\noindent affine projection\qquad 61\par
\noindent Banach space\qquad 29\par
\noindent convex\qquad 61\par
\noindent convex hull\qquad 47\par
\noindent derivative\qquad 28\par
\noindent differentiable\qquad 28\par
\noindent dual\qquad 63\par
\noindent elliptic\qquad 24\par
\noindent elliptic function\qquad 34\par
\noindent extrinsic curvature\qquad 2\par
\noindent Fredholm\qquad 30\par
\noindent Hausdorff distance\qquad 39\par
\noindent H\"older norm\qquad 32\par
\noindent H\"older semi-norm\qquad 31\par
\noindent index\qquad 30\par
\noindent link\qquad 66\par
\noindent local geodesic property\qquad 48\par
\noindent lower barrier\qquad 8\par
\noindent gaussian curvature\qquad 2\par
\noindent graph\qquad 44\par
\noindent great-circular arc\qquad 61\par
}
\vbox{%
\hsize=6.5cm%
\raggedright%
\noindent non-antipodal\qquad 61\par
\noindent normal\qquad 74\par
\noindent open half-space\qquad 59\par
\noindent open hemisphere\qquad 61\par
\noindent outer barrier\qquad 3\par
\noindent regular\qquad 97\par
\noindent shape operator\qquad 1, 74\par
\noindent singular set\qquad 97\par
\noindent smooth\qquad 28\par
\noindent strictly contained in a hemisphere\qquad 50\par
\noindent strictly convex\qquad 2\par
\noindent strong barrier\qquad 89\par
\noindent solution space\qquad 30, 36\par
\noindent southern hemisphere\qquad 61\par
\noindent submanifold\qquad 29\par
\noindent supporting normal\qquad 41\par
\noindent trivialising chart\qquad 29\par
\noindent volume\qquad 94\par
\noindent weak barrier\qquad 89\par
\noindent Weingarten operator\qquad 1\par
\noindent \null\par
}}\par
\newhead{Bibliography}
{\leftskip = 5ex \parindent = -5ex
\leavevmode\hbox to 4ex{\hfil \cite{Brezis}}\hskip 1ex{Brezis H., {\sl Functional analysis, Sobolev spaces and partial differential equations}, Universitext, Springer, New York, (2011)}
\medskip
\leavevmode\hbox to 4ex{\hfil \cite{CaffNirSpruckI}}\hskip 1ex{Caffarelli L., Nirenberg L., Spruck J., The Dirichlet problem for nonlinear second-order elliptic equations. I. Monge Amp\`ere equation, {\it Comm. Pure Appl. Math.} {\bf 37} (1984), no. 3, 369--402}
\medskip
\leavevmode\hbox to 4ex{\hfil \cite{CaffNirSpruckII}}\hskip 1ex{Caffarelli L., Kohn J. J., Nirenberg L., Spruck J., The Dirichlet problem for nonlinear second-order elliptic equations. II. Complex Monge Amp\`ere, and uniformly elliptic, equations, {\it Comm. Pure Appl. Math.} {\bf 38} (1985), no. 2, 209--252}
\medskip
\leavevmode\hbox to 4ex{\hfil \cite{CaffNirSpruckIII}}\hskip 1ex{Caffarelli L., Nirenberg L., Spruck J., Nonlinear second-order elliptic equations. V. The Dirichlet problem for Weingarten hypersurfaces, {\sl Comm. Pure Appl. Math.} {\bf 41} (1988), no. 1, 47--70}
\medskip
\leavevmode\hbox to 4ex{\hfil \cite{Calabi}}\hskip 1ex{Calabi E., Improper affine hyperspheres of convex type and a generalization of a theorem by K. J\"orgens, {\it Michigan Math. J.} {\bf 5} (1958), 105--126}
\medskip
\leavevmode\hbox to 4ex{\hfil \cite{doCarmo}}\hskip 1ex{do Carmo M. P., {\it Differential Geometry of Curves and Surfaces}, Pearson, (1976)}
\medskip
\leavevmode\hbox to 4ex{\hfil \cite{doCarmoII}}\hskip 1ex{do Carmo M. P., {\it Riemannian Geometry}, Birkha\"user, Boston-Basel-Berlin, (1992)}
\medskip
\leavevmode\hbox to 4ex{\hfil \cite{CrandhallIshiiLyons}}\hskip 1ex{Crandall M. G., Ishii H., Lions P.-L., User's guide to viscosity solutions of second order partial differential equations, {\it Bull. Amer. Math. Soc.} {\bf 27} (1992), no. 1, 1--67}
\medskip
\leavevmode\hbox to 4ex{\hfil \cite{Dold}}\hskip 1ex{Dold A., {\it Lectures on Algebraic Topology}, Classics in Mathematics, Springer-Verlag, Berlin-Heidelberg-New York, (1995)}
\medskip
\leavevmode\hbox to 4ex{\hfil \cite{FriedlanderJoshi}}\hskip 1ex{Friedlander F. G., Joshi M., {\it Introduction to the Theory of Distributions}, Cambridge University Press, Cambridge, (1998)}
\medskip
\leavevmode\hbox to 4ex{\hfil \cite{GilbTrud}}\hskip 1ex{Gilbarg D., Trudinger N. S., {\it Elliptic partical differential equations of second order}, Die Grundlehren der mathemathischen Wissenschagten, {\bf 224}, Springer-Verlag, Berlin, New York (1977)}
\medskip
\leavevmode\hbox to 4ex{\hfil \cite{GuanSpruckII}}\hskip 1ex{Guan B., Spruck J., The existence of hypersurfaces of constant Gauss curvature with
prescribed boundary, {\it J. Differential Geom.} {\bf 62} (2002), no. 2, 259--287}
\medskip
\leavevmode\hbox to 4ex{\hfil \cite{GuillemanPollack}}\hskip 1ex{Guillemin V., Pollack A., {\it Differential Topology}, Prentice-Hall, Englewood Cliffs, N.J.,
(1974)}
\medskip
\leavevmode\hbox to 4ex{\hfil \cite{Guttierez}}\hskip 1ex{Gutierrez C. E., {\it The Monge-Ampere Equation}, Progress in Nonlinear Differential Equations and Their Applications, Birkha\"user, (2001)}
\medskip
\leavevmode\hbox to 4ex{\hfil \cite{Milnor}}\hskip 1ex{Milnor J. W., {\it Topology from the differential viewpoint}, Princeton Landmarks in Mathematics, Princeton, (1997)}
\medskip
\leavevmode\hbox to 4ex{\hfil \cite{Nirenberg}}\hskip 1ex{Nirenberg L., {\it Topics in non-linear functional analysis}, Courant Lecture Notes, {\bf 6}, AMS, (2001)}
\medskip
\leavevmode\hbox to 4ex{\hfil \cite{Pog}}\hskip 1ex{Pogorelov A. V., On the improper convex affine hyperspheres, {\it Geometriae Dedicata} {\bf 1} (1972), no. 1, 33--46}
\medskip
\leavevmode\hbox to 4ex{\hfil \cite{RudinOld}}\hskip 1ex{Rudin W., {\it Principles of Mathematical Analysis}, International Series in Pure \& Applied Mathematics, McGraw-Hill, (1976)}
\medskip
\leavevmode\hbox to 4ex{\hfil \cite{Rudin}}\hskip 1ex{Rudin W., {\it Real \& Complex Analysis}, McGraw-Hill, (1987)}
\medskip
\leavevmode\hbox to 4ex{\hfil \cite{RudinII}}\hskip 1ex{Rudin W., {\it Functional Analysis}, International Series in Pure \& Applied Mathematics, McGraw-Hill, (1990)}
\medskip
\leavevmode\hbox to 4ex{\hfil \cite{ShengUrbasWang}}\hskip 1ex{Sheng W., Urbas J., Wang X., Interior curvature bounds for a class of curvature equations. (English summary), {\it Duke Math. J.} {\bf 123} (2004), no. 2, 235--264}
\medskip
\leavevmode\hbox to 4ex{\hfil \cite{Simon}}\hskip 1ex{Simon L., {\it Lectures on geometric measure theory}, Centre for Mathematical Analysis, Australian National University (1984)}
\medskip
\leavevmode\hbox to 4ex{\hfil \cite{Smale}}\hskip 1ex{Smale S., An infinite dimensional version of Sard's theorem, {\it Amer. J. Math.}, {\bf 87}, (1965), 861--866}
\medskip
\leavevmode\hbox to 4ex{\hfil \cite{SmithCGC}}\hskip 1ex{Smith G., Compactness results for immersions of prescribed Gaussian curvature II - geometric aspects, {\it Geom. Dedicata},{\bf 172}, no. 1, (2014), 303--350}
\medskip
\leavevmode\hbox to 4ex{\hfil \cite{SmithPPG}}\hskip 1ex{Smith G., The Plateau problem for convex curvature functions, arXiv:1008.3545}
\medskip
\leavevmode\hbox to 4ex{\hfil \cite{TrudWang}}\hskip 1ex{Trudinger N. S., Wang X., On locally locally convex hypersurfaces with boundary, {\it J.
Reine Angew. Math.} {\bf 551} (2002), 11--32}
\par}
\enddocument